 \documentclass[12pt]{article}
\usepackage{xcolor} 
\RequirePackage[OT1]{fontenc}
\RequirePackage{amsthm,amsmath}
\usepackage[square,numbers]{natbib}
\usepackage{algorithm,algpseudocode}
\usepackage{xr}
\usepackage[figuresright]{rotating}
\makeatletter
\newcommand*{\addFileDependency}[1]{
	\typeout{(#1)}
	\@addtofilelist{#1}
	\IfFileExists{#1}{}{\typeout{No file #1.}}
}
\makeatother

\usepackage{pdfsync}
\usepackage{color}
\usepackage{multirow}
 \usepackage{rotating}
\usepackage{amsmath, amssymb, amsfonts, amsthm}
\usepackage{graphicx}
\usepackage{epstopdf}
\usepackage{geometry}
\usepackage{color}

\usepackage{caption}
\usepackage{subcaption}
\usepackage{mathtools}
\mathtoolsset{showonlyrefs=true}
\newcommand{\field}[1]{\mathbb{#1}}

\newcommand{\p}{\field{P}}

\newcommand{\E}{\field{E}}
\newcommand{\pp}{\mathcal{P}}
\newcommand{\FF}{\mathcal{F}}

\def\argmax{\mathop{\mbox{argmax}}}
\def\argmin{\mathop{\mbox{argmin}}}
\theoremstyle{remark}
\theoremstyle{definition}
\newtheorem{theorem}{Theorem}[section]

\newtheorem{remark}{Remark}[section]
\newtheorem{lemma}{Lemma}[section]
\newtheorem{definition}{Definition}[section]
\newtheorem{corol}{Corollary}[section]

\newtheorem{proposition}{Proposition}[section]

\def\lf{\lfloor}
\def\rf{\rfloor}

\newcommand{\bea}{\begin{eqnarray*}}
	\newcommand{\eea}{\end{eqnarray*}}
\newcommand{\be}{\begin{eqnarray}}
\newcommand{\ee}{\end{eqnarray}}

\usepackage{booktabs}
\numberwithin{equation}{section}
\theoremstyle{plain}

\begin{document}

{	\title{Tractably Modelling Dependence in Networks Beyond Exchangeability}

		\author{\small Weichi Wu \hspace{.2cm} \\
		\small Center for Statistical Sciences,
		Department of Industrial Engineering, Tsinghua University\\
		{\small email: wuweichi@mail.tsinghua.edu.cn }\\
		\small Sofia Olhede \\	\small Institute of Mathematics, Chair of Statistical Data Science, EPFL\\
		{\small email:sofia.olhede@epfl.ch }\\
			\small Patrick Wolfe \\	\small Departments of Statistics and Computer Science, Purdue University\\
			{\small email: patrick@purdue.edu}\\
	}\maketitle
}	
	\begin{abstract}
	We propose a general framework for modelling network data that is designed to describe aspects of non-exchangeable networks. Conditional on latent (unobserved) variables, the edges of the network are generated by their finite growth history (with latent orders) while the marginal probabilities of the adjacency matrix are modeled by a generalization of a graph limit function (or a graphon). In particular, we study the estimation, clustering and degree behavior of the network in our setting. We determine (i) the minimax estimator of a composite graphon with respect to squared error loss; (ii) that spectral clustering is able to consistently detect the latent membership when the block-wise constant composite graphon is considered under additional conditions; and (iii) we are able to construct models with heavy-tailed empirical degrees under specific scenarios and parameter choices. This explores why and under which general conditions non-exchangeable network data can be described by a stochastic block model. The new modelling framework is able to capture empirically important characteristics of network data such as sparsity combined with heavy tailed degree distribution, and add understanding as to what generative mechanisms will make them arise.
	\end{abstract}
	MSC2010 subject classifications: 62G05,62R07,62E20, 62G20,  secondary 53C20.

\section{Introduction}  
The major problem facing modern network analysis is representing sufficient network heterogeneity. Classically heterogeneity is not incorporated in popular models because they are assumed exchangeable~\citep{BiCh2009,OlhedeWolfe2014}; e.g., the models are invariant to permutations, and thus have no nodes that are ``too extreme''. To capture additional heterogeneity research has therefore focused on relaxing models away from standard forms of exchangeability~\citep{veitch2015class,Borgs2014,BCJHV2017}, often modelling edge variables instead of formulating models in terms of the network nodes.

Capturing more network heterogeneity requires us to pose a mechanism for the generation of non-exchangeable networks. 
In this paper we propose a mechanism that can mimic temporal network growth, and is based on the popular graphon model~\citep{BiCh2009,OlhedeWolfe2014}, but is still able to capture additional variability. We shall call our model the `composite graphon model', a special case of `latent order non-anticipatory graphs', a dependent network model class we introduce and describe in detail in Section~\ref{Sec:notation}. Our understanding is encapsulated by using models of latent dependence, and we explore the performance of standard network algorithms with data produced from such a generative model.  These networks exhibit power-law degrees and significant data heterogeneity, typical observed features of non-exchangeability. The common types of data that could require such models include for example citation networks~\citep{ji2016coauthorship}, ecological networks~\citep{sole2001complexity}, technological networks such as the powergrid~\citep{pagani2013power} or communications networks~\citep{broido2019scale}.

Mimicking the mechanism of network growth, or network evolution, to produce an output network, is a very general idea. This idea can be said to be the genesis of other very popular frameworks such as the Barabasi--Albert network scheme~\cite{BraAlb1999}, and temporal evolution underpins the graph processes of Borgs et al.~\cite{borgs2017sparse}.  Borgs et al.~\cite{borgs2017sparse} clarify the complementary relationship of their models to those of~\cite{veitch2015class}, whose modelling framework is nearly identical, even if their achieved results are naturally complementary. Our aim is different: by introducing a model class that has a simple dependence parameter $\chi_n$ (rather than a latent Poisson process) which then drives the degree of exchangeability through controlling the probability of a run in consecutively generated edges, and then understanding the effects on estimation in this setting for networks with $n$ nodes. 

Furthermore, it is important to study the application of standard network algorithms to non-standard network data, as in practice we often cannot check if conditions of exchangeability are satisfied. Motivated by usage of composite likelihood in classical inference, we develop a model with a parameter that quantifies departure from exchangeability.  This is the parameter $\chi_n>0$, quantifying the dependence strength of the latent variables. By defining and studying the impact of $\chi_n$ on estimation, clustering and degree patterns we conclude that when $\chi_n$ is small, the exchangability assumption is adequate even if the true model is non-exchangeable, while the assumption is inadequate exactly when $\chi_n$ approaches $1$.

To be more concrete and granular, the contributions of this paper are fourfold. First, we obtain the minimax estimator of the composite graphon model (as well as the composite version of the stochastic block model) with the $\mathcal L_2$ loss function. Nonparametric regression with stationary/non-stationary time series has already attracted increasing research attention; for example, see \cite{Fanyao2003, ZhouWu2010, WuZhou2015}. Our result can be considered as a network counterpart of \cite{GaoZhou2015}.

Second, we investigate spectral clustering for the composite stochastic block model,  which is the non-exchangeable counterpart of the existing results such as \cite{Rohe2011}. We find that the spectral clustering algorithm is robust to certain dependence structure of edges, which answers the question why the algorithm works well when the assumption of conditional independence for the  Stochastic Block Model (SBM) fails, for example see \cite{SaYu2015}. In addition to in network science, spectral clustering has been applied in many scientific fields including image analysis, data mining and speech recognition, see for instance \cite{gong1995} and \cite{Berkhin2006}, where we do not know of any latent dependence strength.

Third, we construct a model with a  heavy-tailed degree distribution by considering an unobserved latent order $\omega(\cdot,\cdot)$ and certain dependence structure of the edge variables. This shows a new mechanism  resulted from the missing information of correlation between edges  that can produce a power law degree distribution that is not preferential attachment model \cite{BraAlb1999}, \cite{Pa2012} nor the inhomogeneous edge connection probability model \cite{OlhedeWolfe2012}, or dropping the assumption of array exchangeability \cite{veitch2019sampling}, \cite{veitch2015class}.

Fourth, we establish a theoretical framework for the analysis of network data with latent time--order non--anticipatory edge structure by developing a concept closely related to the notion of dependence measure (\cite{Wu2005}) in the literature of time series analysis in which area a similar framework has been successfully set up to accommodate non-stationarity (see for example \cite{zhou2014}). This motivates us to develop many useful mathematical tools in this paper.

The paper is organized as follows. Section \ref{Sec:notation} introduces our notation and basic model structure. Section \ref{Sec:dependence} introduces a dependence structure for the edges. The composite graphon and its minimax estimator are investigated in Section \ref{Chap:graphon}. The composite SBM and the spectral clustering algorithm are studied in Section \ref{Chap:spectral}. In Section \ref{simpleedge} we investigate via an example the basic behaviour of Latent time-order graphs, especially in the case when the dependence is strong. Finally, the proofs of most results are relegated to the supplementary material.

\section{Notation and the Composite Graphon Model}\label{Sec:notation}
For any set $A$, let $|A|$ denote the cardinality of $A$. For a positive integer $n$, we write $[ n ]=\{1,2,...,n\}$. For $n$--dimensional random vectors $\mathbf v=(v_i,1\leq i\leq n)$ and $\mathbf u=(u_i,1\leq i\leq n)$, write $(v_i,1\leq i\leq n)\overset{d}{=}(u_i,1\leq i\leq n)$ if $\mathbf v$ and ${\mathbf{u}}$ have the same distribution. For two numbers $i,j$, denote by $\{i,j\}$ the collection of $i$ and $j$, e.g. $\{i,j\}=\{j,i\}$ and  $\{i,i\}=\{i\}$. Whenever the notation $\{i,j\}$ appears, by default we assume that $i\neq j$. Let $\{[n],[n]\}$ denote the set $\{\{i,j\}, 1\leq i\leq n, 1\leq j\leq n\}$. Denote by $\mathbf 1(\cdot)$ the usual indicator function which is one if the corresponding event is true and zero otherwise. Write $a\wedge b$ for $\min(a,b)$, and $a\vee b$ for $\max(a,b)$. For a graph with adjacency matrix $A_{i,j}$, its marginal probability is  the collection of $\{\p(A_{i,j}=1),\{i,j\}\in \{[n],[n]\}\}$. Each $A_{i,j}$ indicates the presence of an edge between node $i$ and $j$, and we refer to it as an edge variable. It is only an edge if $A_{i,j}=1$. The joint probability of the graph adjacency matrix is $\p(A_{i,j}=a_{i,j},\{i,j\}\in\{[n],[n]\})$ for $\{a_{i,j}, \{i,j\}\in\{[n],[n]\}\}\in \{0,1\}^{N}$ with $N=n(n-1)/2$.  For fixed $(i,j)$, we say $A_{i,j}$ is a edge variable between nodes $i$ and $j$. When $A_{i,j}=1$, we say that there is a linked edge between $i$ and $j$. Let $0^0=1$ as per usual. For any vector $\mathbf V=(v_1,...,v_d)\in \mathbb R^d$, let $\|\mathbf V\|=\sqrt{\sum_{i=1}^dv_i^2}$, and $\|\mathbf V\|_{\mathcal L_p}=
\E((\sum_{i=1}^dv_i^p)^{1/p})$.

With this notation, we will describe the technical framework that we use to quantify a network's departure from exchangeability. We note that lack of exchangeability is a non-property. Non-properties are notoriously hard to quantify, and are in fact not uniquely quantifiable. A network can be non-exchangeable in more than one possible manner. Our method of quantifying departure from exchangeability is merely a possible choice; a non-unique and possibly an imperfect choice. It allows us to quantify that under mild forms of non-exchangeability standard network analysis tools are still applicable and useful.
We shall start by proposing a model which allows us to adjust the degree of departure from exchangeability.

\begin{definition}{Composite Graphon Model.}\label{defcompgraph}
We say a network with adjacency matrix {\it A} is generated by the composite graphon model $f(\cdot,\cdot)$ with respect to a series of latent $i.i.d.$ random variables $\{\xi_i,1\leq i\leq n\}$ if
\begin{description}
\item (a) $\p(A_{ij}=1|\xi_i,\xi_j)=f(\xi_i,\xi_j)$ for some symmetric integrable function $f(\cdot,\cdot)\in [0,1]$; 
\item (b) There exists a bijective map $\omega(\cdot,\cdot):{[n],[n]}\rightarrow [\frac{n(n-1)}{2}]$, such that conditioning on latent variables $\{\xi_i\}$, $B_{s}=A_{\omega^{-1}(s)}$ forms an order $l$ Markovian chain for some $l\geq 0$;
\end{description}
where we call $l$ the long memory parameter, and we include the parameter dependence strength $\chi_n$ which we discuss in Proposition \ref{dependenceprop} in detail.
\end{definition} 

The complete and rigorous definition of a composite graphon is provided in Section \ref{Chap:graphon}, but we give this intuitive definition here to motivate further developments.
When $l=0$, the composite graphon reduces to the usual graphon model. We say a network follows a composite stochastic block model (composite SBM or CSBM) if the composite graphon $f(\cdot,\cdot)$ of its adjacency matrix is block-wise constant, just like composite likelihood ignores correlations.
In fact, the likelihood of composite graphon/SBM is the composite likelihood \citep{VaRe2011} of graphon/SBM, which motives the name of model. By varying the parameter $\chi_n$ continuously we go from a standard exchangeable network, to one exhibiting increasing dependence between the edge variables.

\section{A Graph Sequence Model and a Graph Dependence Measure}\label{Sec:dependence}

Consider a sequence of graphs $\{G_n\}$ with adjacency matrices $\{A_n\}$. 
For a series of $1-1$ corresponding mappings $\omega_n(\{\cdot,\cdot\})$: $\{[n],[n]\}\rightarrow [N]$, 
define the \underline{dependence measure} of the adjacency matrices w.r.t $\omega_n(\{\cdot,\cdot\})$ to be
\begin{align}\label{dependencemeasure}
\ \ \ \ \Delta_n(k):=\max_{s\in\{0,1\}, i,j}|\p(A_{\omega_n(\{i,j\})}=s|\FF_{\omega_n(\{i,j\})-k})-\p(A_{\omega_n(\{i,j\})}=s)|,
\end{align}
 where $\FF_{i}=(B_{-\infty},...,B_i)$, and $B_{\omega_n(i,j)}=A_{i,j,n}$.
 
 For convenience, we let $B_s,s\leq 0$ follow an $i.i.d.$ ${\mathrm{Bernoulli}}(1/2)$ law and be independent of $\{B_s,1\leq s\leq  N\}$. In our paper, we call $\{B_s\}_{1\leq s\leq N}$ the ``ordered edge variables \textcolor{black}{with respect to $\omega$}".  Note that for a sequence of graphs, their adjacency matrices form an array of dependent Bernoulli random variables, with the $n_{th}$ row of the array corresponding to edge variables of a size $n$ graph ordered by $\omega_n(\{\cdot,\cdot\})$. We call $A_{i,j}$ edge variables, and $B_s$ the ordered edge variables.
   For each graph $G_n$, its edge variables behave as a time series indexed by $\omega_n(\{i,j\})$. The quantity produced by \eqref{dependencemeasure} is closely related to the physical  and predictive dependence measure introduced by \cite{Wu2005} that quantifies the degree of dependence of outputs on inputs in (nonlinear) physical systems. It is easy to compute for many stochastic process and  has been used to quantify the strength of dependence in both stationary and non-stationary time series,  see for example \cite{ZhouWu2010}, \cite{WuZhou2015} among many others. We both introduce the dependence measure to networks, and use it to characterise dependence in our network sequence.
   The dependence measure \eqref{dependencemeasure} can be tailored to network data and is easy to calculate due to  Bernoulli random variables being bounded by unity.

In Corollary \ref{Coroldependence} we show that for the non-exchangeable network models built in this paper,   a uniform $M\geq 0$ and a series $\chi_n$ exist, such that $\Delta_n(k)\leq M\chi_n^{|k|}$,
i.e. for each $G_n$, the dependence measure for its edge variables is geometrically decaying with respect to $\omega_n(\{\cdot,\cdot\})$, which we refer to as the Geometric Convergence (GC) assumption. We refer to $\chi_n$ as the {\em dependence parameter}.

Graph sequence models have been well studied in the literature. Among others, for example, \cite{BoRi2009} studied the metrics for sparse graphs via a graph sequence model; \cite{BiCh2009} established a graph sequence model with a scaling parameter $\rho_n$ to address the sparsity issues for exchangeable graph model; and \cite{OlhedeWolfe2014} approximated the graphon model. The minimax rate of the estimation of sparse graphon sequence model was studied by  \cite{Klopp2015} and others. 

\section{Latent Non-Anticipatory Order Graphs}\label{Chap:graphon}
Subsequent to this section, we omit the subscript $n$ if this omission produces no potential for ambiguity. We formulate the Composite Graphon Model (CGM) in this section, which corresponds to a special case of a Latent Non--Anticipatory Order Graph, defined as per below. We take inspiration from 
the non-linear Wold representation used in~\cite{Wu2005} to
build graphs that have latent dependence structure. We can generalize the composite graphon model by the following framework.

 \begin{definition}\label{defmarkovgraph}
We say $G$ is a undirected, edge-based and finite memory Latent Non-Anticipatory Order Graphs w.r.t. latent variables $\xi=(\xi_1,...,\xi_n)$ and mapping $\omega_n(\cdot,\cdot)$ if conditionally on the latent variables $\xi$, there exists a finite $k$, such that for all $i$, the conditional distributions of the edge variables $B_{\omega(i,j)}=A_{i,j}, i\neq j$  
\begin{align}\label{Markov-edge}\p(B_i|B_{i-1},...,B_{-\infty},\xi)=\p(B_i|B_{i-1},...,B_{i-k+1},\xi),\end{align} 
where edge variables $\{B_j,j\leq 0\}$ correspond to the burn-in process, which could be chosen as $i.i.d.$ $Bernoulli(1/2)$ independent of $B_s, 1\leq s\leq N$. Let $l$ be the smallest $k$ such that \eqref{Markov-edge} holds. Then we say $G$ has a memory parameter $l$.
  \end{definition}
  The ``burn-in process'' (e.g. the parameter $l$) has little impact on the network, rather like the starting values for a time series autoregressive (AR) model. The Latent Non-Anticipatory Order Graph has memory parameter $l$ and has the property that the probability of linking an edge variable relies on the past $l-1$  network edge variables. Note that there are network models such that the linkage probability of every edge variables depends on {\em all} the edge variable generated before it. Such models are not Latent Non-Anticipatory Order graphs. An important example is the preferential attachment model (\cite{BraAlb1999}), which we will further discuss in Remark \ref{Newremark5}. Recently
 the asymptotic normality of the affine preferential attachment network models has been studied by \cite{Gaovan2017}.

  Consider a latent time-order graph with memory parameter $l$ and w.r.t. the map $\omega_n(\{\cdot,\cdot\})$. Define $U_i=(B_{i},...,B_{i-l+1})^T$. Denote by $\mathcal X\in \mathbb R^{l}$ the set of $l$-dimensional binary vectors with all entries $0$ or $1$. The following proposition shows that  the dependence between $A_{i,j}$ and $A_{k,l}$ decreases as $|\omega_n(\{i,j\})-\omega_n(\{k,l\})|$ increases.
  \begin{proposition}\label{dependenceprop} Consider a 
  	Latent time-order graph with memory parameter $l$ and mapping $\omega_n(\cdot,\cdot)$. Assume $\min_{a,b\in \mathcal X\times \mathcal X } \p(U_{i}=a|U_{i-l}=b,\xi)\geq \alpha'>0$.
  	Then uniformly for $i,k$ and for all $u_i, u_{i-k}\in \mathcal X$, we have that for $k>l$,
  	\begin{align}
  	\p(U_i=u_i|U_{i-k}=u_{i-k},\xi)-\p(U_i=u_i|\xi)
  	=O(\chi^k\tilde p),
  	\end{align}
	where $\chi=(1-2\alpha')^{1/l}$, and the latent variable $\xi$ is defined in Definition  \ref{defmarkovgraph}. Take
  	\begin{align}
  	\tilde p=\max_{1\leq s\leq {l-1}}\max_{1\leq i\leq n}\max_a|\max_b\p(U_i=a|U_{i-s}=b,\xi)-\min_b\p(U_i=a|U_{i-s}=b,\xi)|.
  	\end{align}
  	Notice that if given $\xi$, $l=0$ and $\{B_i\}$ is an independent series, then $\tilde p=0$.  
  \end{proposition}
  This proposition shows the implication of Definition \ref{defmarkovgraph}.  When a graph sequence model is considered, $\alpha$ may in fact depend on $n$. In that case we shall write $U_{i,n}$ for $U_i$.\ For each $n$, we assume $U_{i,n}$ is an order one Markov process. At this time, $\chi$ depends on $n$ and we denote it by $\chi_n$ when it is important to emphasise the size-dependent relationship. In this article we assume that  $\chi_n<1$. When $\chi_n\rightarrow 1$, $\tilde p$ will be strictly bounded away from zero. Hence $\chi_n$ can be regarded as a proxy of dependence strength.  
 By a straightforward argument using the Markov property, we have the following corollary:
\begin{corol}\label{Coroldependence}
Under the conditions of Proposition \ref{dependenceprop}, we have that 
\begin{align}\label{new4.Fuzhou}
\Delta_n(k)=\max_{\substack{1\leq i\leq N,\\b_s\in \{0,1\},\\ 1-k\leq s\leq i-k}}|\p(B_i=b_i|B_{s}=b_s,s\leq i-k,\xi)-\p(B_i=b_i|\xi)|=O(\chi^k),
\end{align}
where $\chi$ is defined in Proposition \ref{dependenceprop}.
\end{corol}
 If $|\chi|\leq 1-\epsilon$ for some $\epsilon>0$, then we say the latent time-order graph sequence is short-range dependent w.r.t. $\omega(\{\cdot,\cdot\})$. In this case equation \eqref{new4.Fuzhou} implies a geometric decay of $\p(B_i=b_i|B_{s}=b_s,s\leq i-k,\xi)-\p(B_i=b_i|\xi)$ in $k$. This shows a stronger link with an AR(1) process where the  term $\p(B_i=b_i|B_{s}=b_s,s\leq i-k,\xi)$ plays the role of conditional expectation of an observation on another observation $k$ steps ahead  from an $AR(1)$ process. 
 \begin{definition}\label{latent_Markov_def}(Composite graphon model)
A Latent time-order graph is a composite graphon model with respect to $i.i.d.$ latent variables $\xi_i$ if
\begin{align}\label{Compositegraphon}
\theta_{i,j}:=\p(A_{i,j}=1|\xi_1,...,\xi_n)=f_n(\xi_i,\xi_j),\  i\neq j
\end{align}
for some symmetric function $f_n(\cdot,\cdot)$. If in addition $f_n(\cdot,\cdot)$ is block-wise constant in the $\mathbb {R}^2$ plane, 
 then we call \eqref{Compositegraphon} a composite stochastic block model (composite SBM). We shall omit the subscript $n$ when no confusion can arise.
\end{definition}

The parameter of interest in \eqref{Compositegraphon} is independent of the mapping $\omega(\cdot,\cdot)$. This fact is crucial for estimating the composite graphon model without estimating $\omega(\cdot,\cdot)$. Indeed, model \eqref{Compositegraphon} is quite flexible, including the usual graphon model as its special case.
 We then present a general pseudo algorithm for constructing the composite graphon model with memory parameter $l$. For $\omega(\{\cdot,\cdot\})$, define $\omega^{-1}:[N]\rightarrow \{[n]\times [n]\}$ as $\omega^{-1}(k)=\{\{i,j\}, \omega(\{i,j\})=k\}$. For any series $\xi_i,i\in \mathbb Z$, denote by $\xi_{\{i,j\}}=\{\xi_i,\xi_j\}$ for short.
 	\begin{algorithm}
 		\caption{Generative Pseudo Algorithm }
 		\begin{algorithmic}[1]
 			\State Generate $\xi_1,...,\xi_n$. Calculate $f(\xi_i,\xi_j)$.
 			
 			\State For given $\omega$, we first generate $B_1$ as $\p(B_1=1|\xi)=f(\xi_{\omega^{-1}(1)})$.
 			\State  Generate $B_2$ by two parameters
 			$\p(B_2=1|B_1=1,\xi)$ and $\p(B_2=1|B_1=0,\xi)$, which satisfies that $f(\xi_{\omega^{-1}(2)})=\p(B_2=1,B_1=1)+\p(B_2=1,B_1=0)$, where
 			\begin{align}
 			\p(B_2=1,B_1=1|\xi)=f(\xi_{\omega^{-1}(1)})\p(B_2=1|B_1=1,\xi),\notag\\ \p(B_2=1,B_1=0|\xi)=(1-f(\xi_{\omega^{-1}(1)}))\p(B_2=1|B_1=0,\xi).
 			\end{align}
 			By using $\p(B_2=0|B_1,\xi)=1-\p(B_2=1|B_1,\xi)$, in step 2 we have constructed a two-dimensional multivariate Bernoulli model $(B_1,B_2)$.
 			\For{$i = 1$ to $N-1$}
 			\State  Generate $B_{i+1}$ by parameters $\p(B_{i+1}=1|U_i=u_i,\xi)$ which satisfy the following constraints:
 			\begin{align}\label{constraintsmarginal}
 			\p(B_{i+1}=1|\xi)=f(z_{\omega^{-1}(i+1)})=\sum_{u_i}\p(B_{i+1}=1|U_i=u_i,\xi)\p(U_i=u_i|\xi),
 			\end{align}
 		  for $u_i\in \{0,1\}^{i-((i-l)\vee 1)}$,	where $\p(U_i=u_i|\xi)$ could be obtained by the $(i+1-((i-l+1)\vee 1))$-dimensional multivariate Bernoulli model generated in previous iteration.
 			\EndFor \ with output $B_1, B_2,...,B_N$. 
 			
 		\end{algorithmic}\label{Algorithm_4.1}
 	\end{algorithm}
From the generative pseudo--algorithm Algorithm \ref{Algorithm_4.1}, we see that the joint distribution of the edge variables of the composite graphon model, as well as dependence strength $\chi_n$, is fully determined by the following (infinite dimensional) parameters:
 \begin{description}\item (i) $f(\omega^{-1}(i)), 1\leq i\leq N$, 
 \item(ii) $\p(B_i=1|U_i=u_i,\xi)$ for $2\leq i\leq N$, $u_i\in \{0,1\}^{i-(i-l+1)\vee 1}$ with constraints \eqref{constraintsmarginal} for $2\leq i\leq N$.
 \end{description}
 In the classic case, the first  of these specifications is solely via the graphon function, while specification in (ii) breaks the model exchangeability, and so makes the model more flexible. 
{\color{black} We recover the classic graphon model when $\p(B_i=1|U_i=u_i)=\p(B_i=1)=f(\omega^{-1}(i))$}. For any composite graphon model \eqref{Compositegraphon}, we define its associated composite graph as follows.
 \begin{definition}\label{compositegraphdef}
We say $\tilde G$ is a composite graph with respect to a composite graphon model $G$ from \eqref{Compositegraphon}  and with respect to $i.i.d.$ latent variables $z_i$ following $U(0,1)$, if (i) $V(G)=V(\tilde G)$ where $V(G)$ (or $V(\tilde G)$) is the vertex set of $G$ (or $\tilde G$), and if (ii) $\p(\tilde A_{i,j}=1|\xi)=f(\xi_i,\xi_j)$ where $\tilde A_{i,j}$ is the edge variables of $\tilde G$ and (iii) conditioning on $\xi_i, 1\leq i\leq n$, the  edge variables $\{\tilde A_{i,j}, 1\leq i<j\leq n\}$ are independently distributed.
\end{definition}
 The composite graphon $f(\cdot,\cdot)$ is the limit of the composite graph in the sense of \cite{Lov2012}. For any simple graph $F$ with vertex set $V(F)$, the integral of composite graphon $f(\cdot,\cdot)$ on $V(F)$ corresponds to a homomorphism density (see the definition in Section 5.2.2 of \cite{Lov2012}) of $F$ in the composite graph. 

Our model connects to the graph limit and convergence in the language of \cite{Lov2012} in a marginal way similar to the ``composite'' concept in the classic statistics literature, see for instance \cite{VaRe2011}. Furthermore, \cite{BCJHV2017} introduces a latent birth time concept, which is similar to our latent order concept in that it is also temporal. The differences lie in the fact that their latent birth time is for each vertex, while our latent order is for each edge variable, and more fundamentally, lie in the procedure that they drop this latent birth time in their final step to make the model exchangeable such that the labels carry no information while our model is not exchangeable by assuming the information of the labels is missing. As a consequence, \cite{BCJHV2017} involves additional cost for the exchangeability. 

{\color{black}Throughout the paper, we shall focus on the composite graphon model given in \eqref{latent_Markov_def}. For $\alpha\in (0,1]$ and a sufficiently large constant $M$, we define the H\"{o}lder  class
\begin{align*}
\mathcal H_\alpha(M):&=\left\{f:|f(x,y)-f(x',y')|\leq M(|x-x'|+|y-y'|)^\alpha, x\geq y,x'\geq y',\right. \\
&\left.  f\  \text{symmetric}\right\},
\end{align*}
for all $x\geq y$, $x'\geq y'$.
We consider the following scenarios:
\begin{description}
\item(A) $f(\cdot,\cdot)$ is a block-wise constant symmetric function, or
\item(B) $f\in \mathcal{F}_{\alpha}(M)$, where
$\FF_\alpha(M)=\{0\leq f\leq 1:f\in \mathcal H_\alpha(M)\}$.
\end{description}
Under (A), our model reduces to the composite stochastic block model (composite SBM). Under (B), the composite graphon $f(\cdot,\cdot)$ is smooth and estimable. The smoothness is assumed by for example \cite{GaoZhou2015}, \cite{Klopp2015}, \cite{OlhedeWolfe2014}, \cite{Airoldi2013} (which assumes $\alpha=1$) among others.} 

\subsection{Inhomogeneity of the Composite Graphon Model}\label{twostates}
In this subsection, we explore the inhomogeneity introduced by conditional dependence via studying examples of composite SBM with memory parameter $1$. The  memory parameter is given by Definition \ref{defmarkovgraph}. Note that the memory parameter $0$ corresponds to the classical stochastic block model. 
\subsubsection{Inhomogeneity introduced by communities}\label{notationtheta}
We first consider the composite SBM with fixed $k$ communities constructed as follows. Define the map $\gamma: [n]\rightarrow [k]$, which assigns $n$ nodes into $k$ different groups. Define $\tilde \gamma:\{[n]\times[n]\}\rightarrow \{[k]\times [k]\}$ as $\tilde \gamma(i,j)=\{\gamma(i),\gamma(j)\}$. We shall construct a composite SBM such that the edge variable connection probability depends on its previous edge variable with respect to latent (and unobservable) map $\omega$. For this purpose, let $q,l$ be the numbers such that  $\{q,l\}=\tilde \gamma(\omega^{-1}(\omega(\{i,j\})-1))$, i.e., $q,l$ are the communities of vertices of $B_{\omega(\{i,j\})-1}$. Assume
\begin{align}\label{multigenerate}
A_{i,j}=B_{w(\{i,j\})}\sim \left\{ \begin{array}{ll}
\text{Bernoulli($\varrho_{\gamma(i),\gamma(j)}^{0,q,l}$)} &  \text{if}\  B_{\omega(i,j)-1}=0\\
\text{Bernoulli($\varrho_{\gamma(i),\gamma(j)}^{1,q,l}$)} &  \text{if}\  B_{\omega(i,j)-1}=1
\end{array}.
\right.
\end{align}
The probability of linking $(i,j)$ depends on its ``parent'' edge variable ${B_{\omega(i,j)-1}}$, and the communities of the four nodes $q$, $l$, $\tilde \gamma(i,j)$. Conditioning on the latent memberships, \eqref{multigenerate} reduces to an inhomogeneous two-state Markov process. Recently in time series analysis, researchers have developed certain inhomogeneous  models to characterise non-stationarity of integer-valued and categorical data, see for example \cite{truquet2019local}.   
{\color{black}To define a composite SBM with $k$ groups such that $f(\xi_i,\xi_j)=\theta_{\gamma(i)\gamma(j)}$ for $\frac{k(k+1)}{2}$ connection probabilities $\{\theta_{i,j}, i,j\in [k],\theta_{i,j}=\theta_{j,i}\}$}, using Algorithm \ref{Algorithm_4.1}, we specify a composite SBM with the parameters $\{\varrho_{a,b}^{u,c,d},u\in\{0,1\}, \{a,b,c,d\}\in[k]^4 \}$ satisfying the following constraints:
\begin{description}
	\item(a) For $1\leq i\leq j\leq k$,
	\begin{align}\label{Fuzhou.56}
	\varrho_{i,j}=\frac{\varrho_{i,j}^{0,i,j}}{1+\varrho_{i,j}^{0,i,j}-\varrho_{i,j}^{1,i,j}}.
	\end{align}
	\item(b) For $1\leq i\leq j\leq k$, $1\leq s\leq l\leq k$, $\{s,l\}\neq \{i,j\}$, $\varrho_{i,j}^{0,s,l}$ and $\varrho_{i,j}^{1,s,l}$ satisfy
	\begin{align}\label{Fuzhou.57}
	\varrho_{i,j}=\varrho_{i,j}^{0,s,l}(1-\varrho_{s,l})+\varrho_{i,j}^{1,s,l}\varrho_{s,l}.
	\end{align}
\end{description}
In fact, 
each  sub-chain that maps $(i,j)\rightarrow(i,j)$ describes a homogeneous Markov chain, i.e., if we consider any consecutively generated  edge variables which connect the vertices that belong to the same pair of groups $(i,j)$, then these edge variables form a homogeneous Markov Chain with stationary probability ($\varrho_{i,j},1-\varrho_{i,j}$). If $k=1$ (corresponding to the scenario of only one group), constraints (a), (b) degenerate to a strictly stationary 2-states Markov process. From this point of view, the inhomogeneity is introduced by the specification of communities with stationary probabilities ($\varrho_{1,1},1-\varrho_{1,1}$). 
\subsubsection{Inhomogeneity Introduced by Individuals}
{\color{black}Another source of inhomogeneity is  {\color{black} due to the dependence introduced by the latent position $\omega_n(i,j)$}. Consider the single group, or $k=1$ case. As described by the algorithm under definition  \ref{latent_Markov_def}, another composite SBM could be specified by
	\begin{align}\label{inhomomarkv-singlegroup}
	\varrho_{1,1}=\frac{q_{0,i}}{1-q_{1,i}+q_{0,i}}=\frac{q_{0,j}}{1-q_{1,j}+q_{0,j}} \ \ \text{for $1\leq i<j\leq N$},
	\end{align}
	where $q_{0,i}=\p(B_{i}=1|B_{i-1}=0)$ and  $q_{1,i}=\p(B_{i}=1|B_{i-1}=1)$. 
	It is not hard to see that \eqref{inhomomarkv-singlegroup} defines an inhomogeneous Markov chain. The conditional connection probabilities of the edge variables depend on the  its ordered edge variables}' positions in the history of the Markov chain. Since all nodes belong to the same community, the inhomogeneity is evident only at the individual level. This is very different from its Erd\H{o}s-R\'{e}nyi counterpart, in which each node is stochastically equivalent. 
\section{Minimax Rate Estimator of the Composite Graphon}\label{Sec:Minimax}
In this section, we discuss the minimax estimator of the composite graphon with respect to squared error loss.
 Let $\mathcal Z_{n,k}=\{z:[n]\rightarrow [k]\}$ be the collection of all possible mappings from $[n]$ to $[k]$. Then for any $\bar{z}\in \mathcal Z_{n,k}$, $\{\bar{z}^{-1}(a):a\in [k]\}$ forms a partition of $[n]$, or equivalently: $\cup_{a\in [k]}\bar{z}^{-1}(a)=[n]$ and $\bar{z}^{-1}(a)\cap \bar{z}^{-1}(b)=\emptyset$ for any $a\neq b\in [k]$.
In the following, we adopt the notation of \cite{GaoZhou2015}. Define
\begin{align}
\bar \eta_{ab}(\bar{z})=\frac{1}{|\bar{z}^{-1}(a)||\bar{z}^{-1}(b)|}\sum_{i\in \bar{z}^{-1}(a)}\sum_{j\in \bar{z}^{-1}(b)}\eta_{ij} \text{ for $a\neq b\in [k]$},\\
\bar \eta_{aa}(\bar{z})=\frac{2}{|\bar{z}^{-1}(a)|(|\bar{z}^{-1}(a)|-1)}\sum_{i\in \bar{z}^{-1}(a),j\in \bar{z}^{-1}(a),i< j}\eta_{ij} \text{ for $a\in [k]$, $|\bar{z}^{-1}(a)|>1$},
\end{align}
where $A=A_{i,j}$ is the adjacency matrix. In this section, let $\theta_{i,j}=\p(A_{i,j}=1|\xi)=f(\xi_i,\xi_j)$.
Define the estimate $\hat \theta_{ij}=\hat Q_{\hat z(i)\hat z(j)}$ where
\begin{align}\label{estimate1}
(\hat Q,\hat z)=\argmin_{Q\in \mathbb R^{k\times k}_{sym},\bar z\in\mathcal Z_{n,k}} L(Q,\bar z),\\
L(Q,\bar z)=\sum_{a,b\in [k]}\sum_{(i,j)\in \bar{z}^{-1}(a)\times \bar{z}^{-1}(b),i< j}(A_{ij}-Q_{ab})^2.\label{estimate2}
\end{align}\\
This procedure \eqref{estimate1} is referred to as minimizing combinatorial least squares (see \cite{GaoZhou2015}). The word combinatorial is inserted, as we need to determine group membership $(z)$ which is a combinatorial problem, rather than solely estimating a parameter by weighted averaging.
Straightforward calculations show that $\hat Q_{ab}=\bar A_{ab}(\hat z)$ for all $a,b\in [k]$.
Therefore similarly to \cite{wolfe2013nonparametric} we propose a block constant estimates for the composite graphon model. For this purpose, \textcolor{black}{we define \begin{align}
\Theta_k=\{\{\theta^\diamond_{i,j}\}\in [0,1]^{n\times n}:\theta^\diamond_{i,i}=0,\theta^\diamond_{i,j}=Q_{ab}=Q_{ba}, \forall (i,j)\in \bar{z}^{-1}(a)\times \bar{z}^{-1}(b),\notag\\ \text{for some}\  Q_{ab}\in [0,1], \bar{z}\in \mathcal Z_{n,k}\}.
\end{align}}Define the true value on each block by $\{Q^*_{ab}\}\in [0,1]^{k\times k}$, and the oracle assignment $z^*\in\mathcal Z_{n,k}$,
writing $\theta_{i,j}=Q^*_{z^*(i)z^*(j)}$. For {each estimate} $\hat z$, define $\tilde Q_{ab}\in [0,1]^{k\times k}$ by $\tilde Q_{ab}=\bar \theta_{ab}(\hat z)$ and $\tilde \theta_{i,j}=\tilde Q_{\hat z(i)\hat z(j)}$ for $i\neq j$. For all $i\in \lf n\rf$, let $\hat \theta_{i,i}=\tilde \theta_{i,i}=\theta_{i,i}=0$ (as we have assumed no self loops). Define $n_a=|\bar{z}^{-1}(a)|$. 
We first consider the composite SBM model.
 Recall the definition of $\chi_n$ from Corollary \ref{Coroldependence}. In this section, we assume that the memory parameter is bounded. 
\begin{theorem}\label{cconsistent}
Considering the composite SBM model $G_n$ with $k$ groups. Assume the conditions of (A) hold. 
For any constant $C'>0$,
there is a constant $C>0$ which only depends on $C'$, such that
\begin{align}\label{new.61}
\frac{1}{n^2}\sum_{i,j\in [n]\times [n]}(\hat \theta_{i,j}-\theta_{i,j})^2\leq C\left(\frac{k^2}{n^2}+\frac{\log k }{n}\right)(1-\chi)^{-2}
\end{align}
with probability at least $1-\exp(-C'n\log k)$ uniformly over $\theta\in\Theta_k$, and
\begin{align}
\sup_{\theta\in\Theta_k}\E\left\{(\hat \theta_{i,j}-\theta_{i,j})^2\right\}\leq C_1\left(\frac{k^2}{n^2}+\frac{\log k }{n}\right)(1-\chi)^{-2}
\end{align}
for all $k\in [n]$ with some universal constant $C_1>0$.
\end{theorem}
{\it Proof.}
This proof proceeds along the lines of \cite{GaoZhou2015}. By using the fact that $L(\hat Q,\hat z)\leq L(Q^*,z^*)$, we have that
\begin{align}
\|\hat \theta-\theta\|^2\leq 2\langle \hat \theta-\theta,A-\theta\rangle.
\end{align}
Direct calculations show that
\begin{align}\label{new.21-2017}
\langle \hat \theta-\theta,A-\theta\rangle\leq \|\tilde \theta-\hat \theta\|\left|\left\langle \frac{\hat \theta-\tilde \theta}{\|\hat \theta-\tilde \theta\|},A-\theta\right\rangle\right|+(\|\tilde \theta-\hat \theta\|\notag\\+\|\hat \theta-\theta\|)\left|\left\langle \frac{\tilde \theta-\theta}{\|\tilde \theta- \theta\|},A-\theta\right\rangle\right|.
\end{align}
Define for any $\bar{z}\in \mathcal Z_{n,k}$,  
\[\hat \theta(\bar{z})=\argmin _{Q\in \mathbb R_{sym}^{k\times k}}L(Q,\bar{z}),\] 
and $\tilde \theta(\bar{z})=\argmin _{Q\in \mathbb R_{sym}^{k\times k}}\tilde L(Q,\bar{z})$, where
$$\tilde L(Q,\bar{z})=\sum_{a,b\in [k]}\sum_{(i,j)\in \bar{z}^{-1}(a)\times \bar{z}^{-1}(b),i< j}(\theta_{ij}-Q_{ab})^2.$$
Note that $\tilde \theta=\tilde \theta(\hat z).$ Thus,
by the property of least squares estimator, we have that $\|\tilde \theta-\theta\|\leq \|\hat \theta -\theta\|$. As a result, we have
\begin{align}
\|\tilde \theta-\hat \theta\|\leq 2\|\hat \theta-\theta\|.
\end{align}
It follows from Lemmas \ref{New.lemma_7} and \ref{New.lemma_8} in the supplementary material that the terms
\begin{align}\label{new.23-2017}
\left|\left\langle \frac{\hat \theta-\tilde \theta}{\|\hat \theta-\tilde \theta\|},A-\theta\right\rangle\right|,\
\left|\left\langle \frac{\tilde \theta-\theta}{\|\tilde \theta-\theta\|},A-\theta\right\rangle\right|,
\end{align}
could be bounded by $\sqrt{k^2 +n\log k}(1-\chi)^{-1}$ with probability at least $1-\exp(-C'n\log k)$. Finally, the theorem follows from combining \eqref{new.21-2017}--\eqref{new.23-2017}.
\hfill $\Box$\\

 Note that when a graph sequence model is considered, the factor $\chi=\chi_n$ is allowed to depend on $n$. This is discussed further in Remark \ref{explainp1p0}. Regarding the convergence rate, the term $\frac{k^2}{n^2}$  corresponds to the estimation of $k^2$ unknown parameters with an order of $n^2$ observations (edge variables), and the term $\frac{\log k}{n}$ corresponds to clustering rate, see for example \cite{GaoZhou2015} and \cite{Klopp2015}. 
Meanwhile, the term $(1-\chi)^{-2}$ is the effect of the non-exchangeability due to the latent order $\omega(\{\cdot,\cdot\})$ and the conditional  dependence between edges given latent  variables $(\xi_i)$. 
When $p_0=p_1$, our model reduces to the usual SBM. In this situation, the second part of the convergence rate degenerates to a constant, while the first part agrees with the rate in \cite{GaoZhou2015}, which has been shown to be rate optimal.

{\color{black} Consider the composite graphon model generated from Algorithm \ref{Algorithm_4.1} and specified by $\{\xi_1,...,\xi_n, A,\omega(\{\cdot,\cdot\}),l\}$ with symmetric composite grahon} $f(\cdot,\cdot)$, i.e., $\theta_{i,j}=f(\xi_i,\xi_j)$, where latent variables $\{\xi_i,1\leq i\leq n\}$  are $i.i.d$ $U(0,1)$, $\omega$ is the latent order of  edge variables and $l$ is the memory parameter.  
When $f\in \mathcal F_\alpha(M)$, which is {\color{black} the bounded H\"{o}lder's class} defined in (B) of Section \ref{Chap:graphon},
arguments of Gao et al. (2015) show that there exists an {\em oracle} (\cite{wolfe2013nonparametric}) $z^+_k\in \mathcal Z_{n,k}$ such that for some universal constant $C>0$,
\begin{align}\label{thetastar}
\frac{1}{n^2}\sum_{a,b\in [k]}\sum_{i\neq j:z^+_k(i)=a, z^+_k(j)=b}(\theta_{i,j}-\bar \theta_{a,b}(z^+_k))^2\leq CM^2(\frac{1}{k^2})^{\alpha \wedge 1}.
\end{align}
Consider
\begin{align}\label{Est0111}
(\theta^*,z^*)=\argmin_{Q\in \mathbb R^{k\times k}_{sym},z\in\mathcal Z_{n,k}} \tilde L(Q,z),
\tilde L(Q,z)=\sum_{a,b\in [k]}\sum_{(i,j)\in \bar{z}^{-1}(a)\times \bar{z}^{-1}(b),i< j}(\theta_{i,j}-Q_{ab})^2.
\end{align} By choosing $k=n^{\frac{1}{(1+\alpha \wedge 1)}}$, we have the following theorem:
\begin{theorem} \label{Marginalgraphon}
Consider a composite graphon model $G=\{\xi_1,...,\xi_n, A,\omega(\{\cdot,\cdot\}),l\}$ 
Assume the conditions of (B) hold. Then there exist constant $C,C'$
\begin{align}
\frac{1}{n^2}\sum_{i,j\in [ n]}(\hat \theta_{i,j}-\theta_{i,j})^2\leq C n^{\frac{-2(\alpha \wedge 1)}{1+\alpha\wedge 1}}(1-\chi)^{-2}\log n
\end{align}
with probability at least $1-\exp(-C'n)$, uniformly over $f\in \mathcal F_\alpha(M)$. Furthermore,
\begin{align}\label{expupperbound}
\sup_{f\in \mathcal {F}_\alpha(M)}\E\left\{\frac{1}{n^2}\sum_{i,j\in \lf n\rf}(\hat \theta_{i,j}-\theta_{i,j})^2\right\}\leq C_1 n^{\frac{-2(\alpha \wedge 1)}{1+\alpha\wedge 1}} (1-\chi)^{-2}\log n,
\end{align}
for some constant $C_1>0$.
\end{theorem}
{\it Proof.} By similar arguments to those of \cite{GaoZhou2015},
we have that
\begin{align}
\label{thetathetastar}
\|\hat \theta-\theta^*\|\leq \max\{16 (D+C)^2, 4(B+D)E\},
\end{align}
where \begin{align}
E=\|\tilde \theta-\hat \theta\|\leq 2\|\hat \theta-\theta\|,\notag\\ B=\left|\left\langle \frac{\hat \theta-\tilde \theta}{\|\hat \theta-\tilde \theta\|}, A-\theta\right\rangle\right|,C=\|\theta-\theta^*\|, D=\left|\left\langle \frac{\theta^*-\tilde \theta}{\|\theta^*-\tilde \theta\|}, A-\theta\right\rangle\right|.
\end{align}
 Direct calculations show that Lemma \ref{New.lemma_7} in the supplementary material still holds when replacing $\left|\left\langle \frac{\tilde \theta-\theta}{\|\tilde \theta-\theta\|},A-\theta\right\rangle\right|$ with $\left|\left\langle \frac{\theta^*-\tilde \theta}{\|\theta^*-\tilde \theta\|}, A-\theta\right\rangle\right|$ .
By 
Lemma \ref{New.lemma_7}, Lemma \ref{New.lemma_8} in the supplementary material and \eqref{thetastar},
the theorem follows. \hfill $\Box$

\begin{remark}\label{explainp1p0}
Comparing with the results of the usual SBM and the graphon model as discussed in \cite{GaoZhou2015}, our Theorems \ref{cconsistent} and \ref{Marginalgraphon} introduce an extra factor of $(1-\chi_n)^{-2}$. 
The convergence rate is therefore slow when $\chi_n$ is close to $1$,  of which the situation indicates the strong conditional dependence between the  edge variables on latent variables, see Proposition \ref{dependenceprop}. A straightforward calculation using Proposition \ref{dependenceprop} shows that if given $\xi_n$, $l=0$ and $\{B_i\}$ is an independent series, 
the rate is fully consistent with previous results in the sense that both the model and the rate recovers the known optimal rate of \cite{GaoZhou2015}.
\end{remark}
\begin{remark}\label{sparsityexp1}
Consider the composite graphon model defined in Theorem \ref{Marginalgraphon} with  ordered edge variables $B_s,1\leq s\leq n$. We construct a sparse composite graphon model $\tilde G=\{\xi_1,...,\xi_n,\tilde A,\omega(\{\cdot,\cdot\}),l\}$ by $\p(\tilde U_i=u_i|\tilde U_{i-1}=u_{i-1},\xi)=\rho_n\p( U_i=u_i| U_{i-1}=u_{i-1},\xi)$, where $U_i=(B_i,...,B_{i-l+1})^T$,  $\tilde U_i=(\tilde B_i,...,\tilde B_{i-l+1})^T$, $\tilde B_i, 1\leq i\leq N$ are  ordered  edge variables of $\tilde A$, and $\rho_n$ is a positive sequence that converge to $0$.
We therefore represent the ``sparsity'' by the parameter $\rho_n>0$. This parameter was used by \cite{BoRi2009} to uniformly control the success probability across all nodes and uniformly controls the number of edges present.  
 {\color{black} Straightforward calculations} show that the upper bound of RHS of \eqref{expupperbound} could be lowered to the order of $\min_k\{\rho_n^2(\frac{1}{k^2})^{\alpha\wedge 1}+\rho_n(\frac{k^2}{n^2}+\frac{\log k}{n})(1-\chi)^2\}$, which coincides with the upper bound of that in \cite{Klopp2015}.  In Section \ref{simpleedge}, we shall see  scenarios of homogeneity that $\frac{\max_{i,j} \p(A_{i,j}=1)}{\min_{i,j} \p(A_{i,j}=1)}\rightarrow \infty$ which is able to produce power law degree distribution. This scenario cannot be captured by scale parameter $\rho_n$.  As a result, we do not focus on the scaled sparse model in detail. 
\end{remark}
{\color{black}\begin{remark}\label{Consiscondition}
Assume the setting of the composite SBM sequence such that $\chi=\chi_n$ is regulated by $n$. Assume that the conditions of Theorem \ref{cconsistent} hold. By Theorem \ref{cconsistent}, when the number of communities $k$ is fixed, a sufficient condition for the consistency of the $\mathcal L_2$ estimator \eqref{estimate1} is $(\sqrt n(1-\chi_n))^{-1}=o(1)$. Theorem \ref{Marginalgraphon} indicates that the estimator will be inconsistent under strong dependence such that $\chi_n$ approaches $1$ at a rate faster than $\frac{1}{\sqrt n}$. Similarly, Theorem \ref{Marginalgraphon} implies that when the composite graphon $f\in \mathcal F_\alpha(M) $ and $k=\lf n^\frac{\alpha\wedge 1}{1+\alpha\wedge 1}\rf$, a sufficient condition for the consistency of the $\mathcal L_2$ estimator \eqref{estimate1}  is $(n^{\frac{\alpha\wedge1}{1+\alpha\wedge1}}(1-\chi_n))^{-1}=o(1)$.
\end{remark}}
{\bf Interpretation:} In other words, our estimator is consistent for 
 $1-\chi_n=\Omega(n^{-a})$, $0<a<1/2$, where we write $a_n=\Omega(b_n)$ for series $a_n$, $b_n$ if $b_n=O(a_n)$. By proposition \ref{dependenceprop} this means $\max_{a,b\in \mathcal X\times \mathcal X } \p(U_{i}=a|U_{i-l}=b,\xi)\leq 1-Cn^{-a}$ for some constant $C$.

\begin{remark}\label{matrixorfunction}
Theorem~\ref{Marginalgraphon} determines that matrix estimation can be done for this problem, i.e.  the sampled graphon can be estimated from an observed adjacency matrix. It does not necessarily relate to the underlying graphon function, unless we derive further results. As noted in~\cite{wolfe2013nonparametric} the mean square error of the estimate of $f(x,y)$ can be directly related to the matrix mean square error of~\eqref{thetathetastar}. An issue with this statement is that the discretized $p_{ij}=f(\xi_i,\xi_j)$ is still random as $\xi_i$ is random and so statements can be made either marginally or conditionally on $\xi_1,\dots, \xi_n$. The added problem of estimating $f(x,y)$, the function, is to use an appropriate metric, and factor out measure preserving transformations. 
\end{remark}

It is discussed in {\cite{GaoZhou2015}} and~\cite{wolfe2013nonparametric} that the  graphon model is closely related to non-parametric regression with unknown design and $i.i.d.$ errors. Consider the one-dimensional regression problem $y_i=f(z_i)+e_i$, where $z_i,1\leq i\leq n$ are $i.i.d.$ samples, and $e_i$ are zero mean errors. When $f\in \mathcal H_\alpha(M)$ and $e_i$ are $i.i.d$ normals, the local polynomial estimator achieves the minimax rate $n^{-\frac{2\alpha}{1+\alpha}}$ under the squared error loss $\frac{1}{n}\sum_{i\in [n]}\{\hat f(z_i)-f(z_i)\}^2$.
When $e_i$ is a short range dependent non-stationary time series for example the piecewise locally stationary time series in \cite{zhou2013}, Lemma 5 in {\cite{ZhouWu2010}} shows that $\E(e_ie_j)=O(\eta^{|i-j|})$ for some $\eta\in (0,1)$. It follows from this fact and Proposition 1.13 in \cite{Tsybakov08}, that the convergence rate of the local polynomial estimator with non-stationary time series error has the same order as with $i.i.d.$ error. 

However, under the situation that the design is unknown, an additional difference arises between the time series error and the $i.i.d$ error due to the unknown chronological order. 
Indeed, missing chronological order affects time series but not the $i.i.d.$ errors. 
Surprisingly, for the time series error, the impact of the missing chronological order on the estimation is negligible in terms of order under certain situation. %
Suppose the edge variables of the graph are short range dependent with respect to the mapping $\omega(\{\cdot,\cdot\})$. 

Recall $L(Q,\bar z)$ in \eqref{estimate2}, and define the new objection $L^0(Q,\bar z)$ by replacing $A$ with $\theta=\E(A|\xi)$. For given $\bar z$, let $\tilde Q(\bar z)$ and $Q^0(\bar z)$ be the minimizer of $L(Q,\bar z)$ and  $L^0(Q,\bar z)$, respectively. In fact, $\tilde Q(\bar z)$ and $Q^0(\bar z)$  are the average among partitions of adjacency matrix and of true but unknown conditional linkage probability matrix, respectively. 
Thus $\E(\tilde Q(\bar{z}))=Q^0(\bar{z})$. Since  Bernoulli random variables are bounded, we show in our paper that under mild conditions and similar to the time series counterpart, 
the deviations between the average and the mean of the edge variables are bounded uniformly over all possible partitions $z$ 
as if the  edge variables are conditional independent.

 The time series structure of our model is important for many real applications. For example, \cite{GaoZhou2015}  relates link prediction to the graphon model. In real application, links can be modeled by time series in dynamic network, see for example \cite{Sarkar2014}. Despite the convergence rate of the proposed estimator for the composite graphon model is similar to that under usual graphon model, its moment behavior is different under two scenarios. This will further impact on any estimation and hypothesis testing procedure, see for example \cite{Bickle2011} and \cite{Bickle2015}. 
\section{Spectral Clustering Algorithm for Composite Stochastic Block Model}\label{Chap:spectral}

In the previous sections we investigate the estimation of the composite graphon model and the composite stochastic block model. In addition to estimation of the linkage probabilities, community detection is another research topic in the network analysis. The connection between estimation and spectral clustering is complicated, and they are not identical problems. A good estimation result for the block heights of a stochastic blockmodel does not necessarily guarantee a good community detection result. For a more detailed discussion of the link between parameter estimation and spectral clustering, we refer to \cite{GaoZhou2015}.  
In the area of community detection, spectral clustering and its variants have already been widely applied (\cite{Von2007}).  The consistency of spectral clustering for certain exchangeable network models has been studied by for example \cite{Rohe2011}, \cite{QinRohe2013}, \cite{zhao2012consistency}, \cite{li2018hierarchical} among others. In the following, we shall study the performance of spectral clustering for estimating the composite SBM.
\subsection{Re-parameterization of the Composite Stochastic Block Model}\label{MarginalSBM}
 By choosing a block-wise constant symmetric $f(\cdot,\cdot)$ in Definition \ref{latent_Markov_def}, the composite SBM has the form of $\theta_{i,j}=\p(A_{i,j}=1|\xi_i,\xi_j)$, where $\{\xi_i,1\leq i\leq n \in \mathbb R^k\}$ are $i.i.d.$ latent vectors, with one entry equal to one and all other entries equal to zero. Let $\theta$ be an $n\times n$ matrix, with $\theta_{i,j}$ its $(i,j)_{th}$ entry. Then $\theta$ could be parameterized as 
\begin{align}
\theta=Z B^\dag Z^T,
\end{align}
where $B^\dag \in [0,1]^{k\times k}$ is full rank and symmetric, and $Z\in \mathbb R^{n\times k}$ is a matrix with $i{th}$ row $\xi_i$ such that it has one $1$ in each row and at least one $1$ in each column. For each node $i$, we say it belongs to  group $j$ if $\xi_{i,j}$, which is the $j_{th}$ element of $\xi_{i}$, equals $1$. 
With the re-parameterization, we are able to define the graph Laplacian, which is essential for the spectral clustering algorithm.  Define diagonal matrices $D$ and $\bar D$ with diagonal elements $D_{i,i} $ and $\{\bar D_{i,i}\}, i=1...n$, respectively, where  
\begin{align}
 D_{i,i}=\sum_{k=1}^n  A_{i,k},\ \ \ \bar D_{i,i}=\sum_{k=1}^n \theta_{i,k}.\end{align}
 Define $L$ and $\bar L$ for the Laplacian of $A$ and $\theta$, respectively, as 
 \begin{align}
 L= D^{-1/2} A  D^{-1/2},\bar L=\bar D^{-1/2}\theta \bar D^{-1/2}.
\end{align}

 Note that  $\bar L$ is the population version of $L$ since the former is the Laplacian of $\theta$ and the latter is the Laplacian of adjacency matrix $A$. Both $L$ and $\bar L$  depend on the number of nodes $n$.  Let $c_i=\bar D_{i,i}/n$
 and $\tau_n=\min_{i=1,...,n}c_i$.  We shall write
$L$ as $L^{(n)}$, $\bar L$ as $\bar L^{(n)}$ and $\tau$ as $\tau_n$ when we need to emphasise the sample size. In the remainder of this section, we assume $Z$ is unknown but fixed (unless specified).
After obtaining $L$, the spectral clustering algorithm is given by:
\begin{enumerate}
	\item Compute the eigenvectors $u_1,..., u_k$ w.r.t. the first $k$ largest eigenvalues of $L$.
	\item Run a $k$-means algorithm on vectors $y_1,...,y_n$, $\{y_i\}_{1\leq i\leq n}\in \mathbb R^{1\times k}$ to cluster them into clusters $C_1,...,C_k$, where $y_i$ is the $ith$ row of matrix $U$, an $n\times k$ matrix such that the $jth$ column of $U$ is $u_k$.
\end{enumerate}
Then node $i$ is in class $g$ if $y_i$ is assigned to $C_g$.
\subsection{Properties of Mis-clustered Nodes}
For simplicity of exploration and the ease of comparison, we will use the notion of~\cite{Rohe2011}.
In order to discuss the property of mis-clustered nodes, we first introduce the following notation.   
In addition we define
\begin{align}
P_n=\max_{j=1,...,k}(Z^TZ)_{j,j}.
\end{align}
 We then give two properties of composite SBM. The validity of the properties could be shown similarly to Lemma 3.1 and Lemma 3.2 of {\color{black}\cite{Rohe2011}}, and so we omit the proof for the sake of brevity.
\begin{description}
\item(a) There exists a matrix $V_1\in \mathbb R^{k\times k}$ such that the columns of $ ZV_1$ are the eigenvectors of $\bar L$ which correspond to the nonzero eigenvalues. In addition, $z_iV_1=z_jV_1$ if and only if $z_i=z_j$.
\item(b) Let $V_2\in \mathbb R^{n\times k}$ be a matrix whose orthonormal columns are the eigenvectors which correspond to the ordered largest $k$ eigenvalues of $L$ (in absolute value). Let $c_i,1\leq i\leq n$ be the centroid corresponding to the $i_{th}$ row of $V_2$. Let the columns of $U,\bar U \in \mathbb R^{n\times k}$ be $k$ orthonormal eigenvectors of $LL$ and $\bar L\bar L$ (recall $L$ and $\bar L$ are symmetric matrix) which correspond to the first $k$ largest eigenvalues of the two matrices in absolute value, respectively. Define matrices $O_{1}$ and $O_{2}$ with the singular decomposition $\bar U^TU=O_{1}\Sigma O_{2}^T$, where $O_1, O_2$ are orthonormal matrices and $\Sigma$ is a diagonal matrix. Let $O=O_{1}O_{2}^T$. 
    Then $\|c_i-z_iV_1O\|<\frac{1}{\sqrt{2P_n}}$ if and only if  $\|c_i-z_iV_1O\|<\|c_i-z_jV_{1}O\| $ for any $z_i\neq z_j$.
\end{description}

Under conditions of Theorem \ref{Spectralthm3} below, the Davis-Kahan Theorem \cite{davis1970rotation} shows that $\|V_2-ZV_1O\|_F=o(1)$ almost surely, which leads to that the corresponding eigenvectors of the observed graph Laplacian $L$ is close to that of the population graph Laplacian $\bar L$; { see \cite{Rohe2011} for a detailed introduction of the Davis-Kahan Theorem.} 
As a result, by (a), (b), we define the set of mis-clustered nodes as \begin{align}
\mathcal M=\left\{i:\|c_i-z_iV_1O\|\geq\frac{1}{\sqrt{2P_n}}\right\},
\end{align}
since similarly to the argument in \cite{Rohe2011}, we can show that if any node $i\not \in \mathcal M$, then $i$ will be correctly clustered by spectral clustering algorithm.
For any symmetric matrix $M$, define $\lambda(M)$ to be the eigenvalues of $M$. For any interval $S\in \mathbb R$,
define $\lambda_S(M)=\{\lambda(M)\cap S\}$,  Let $\bar \iota_1^{}\geq...\geq \bar \iota_n^{}$ be the elements of $\lambda(\bar L^{}\bar L^{})$, and $ \iota_1^{}\geq...\geq \iota_n^{}$ be the elements of $\lambda( L^{} L^{})$. Define $G(\chi,N,u)=\sum_{r=0}^{N}r^u\chi^{r/2}$, where $\chi$ is defined in Proposition \ref{dependenceprop}.

Before stating Theorem \ref{Spectralthm3} regarding the performance spectral clustering for estimating the composite SBM, we present the following Proposition \ref{lemmaconsistent} which studies the tail probability of $\|LL-\bar L\bar L\|_F$. The latter is the difference between the population version of and the usual graph Laplacian. The graph Laplacian plays a central role in the spectral clustering, therefore the difference $\|LL-\bar L\bar L\|_F$ is key to study the asymptotic behavior of the corresponding clusters.  The proof of Proposition \ref{newprop6} and Proposition \ref{lemmaconsistent} are inspired by \cite{zhou2014}, as well as \cite{Rohe2011} and \cite{zhou2014}, respectively.  The proof of Theorem \ref{Spectralthm3} rests on the the following Proposition \ref{lemmaconsistent}.
\begin{proposition}\label{lemmaconsistent}
	Under conditions of Theorem \ref{Spectralthm3}, there exist sufficiently large positive constants $\eta_0,\eta_1, M'$ such that if $n\geq M'$
	\begin{align}
	\p(\|LL-\bar L\bar L\|_F\geq \frac{\log n}{\tau^2n^{1/2}}G^{1/2}(\chi,N,3)(1-\chi)^{-1/2})\leq \zeta(n),\label{summable100}
	\end{align}
	where \begin{align}
	\zeta(n)=\frac{\eta_0}{n\log ^4n}+\eta_1n^{-2}.
	\end{align}
\end{proposition}
\begin{theorem}\label{Spectralthm3}
Consider a size $n$ composite SBM with a fixed unknown mapping $\omega(\{\cdot,\cdot\})$. Denote by $k_n$ the number of groups of nodes, and by $n_{k_n}$ the corresponding group size. Let $|\lambda_1|>...>|\lambda_{k_n}|$ be the absolute values of ordered  $k_n$ largest absolute and non-zero eigenvalues of $\bar L$. Assume that $n^{-1/2}(\log n)^2G^{1/2}(\chi,N,3)(1-\chi)^{-1/2}=O(\lambda_{k_n}^2)$, and $\tau_n^2>M/\log n$ for a sufficiently large constant $M$.  
Then we have that the number of miss-specified nodes has the order of \begin{align}
|\mathcal M|=o\left(\frac{P_n\log^2 n}{\lambda_{k_n}^4\tau_n^4n}G^{}(\chi,N,3)(1-\chi)^{-1}
\right),\quad a.s.
\end{align}
provided that $G(\chi,N,3)(1-\chi)\log n$ is sufficiently large such that $\zeta(n)$ (defined as in Proposition \ref{lemmaconsistent} below) is summable.
\end{theorem}

The conditions on the eigenvalues and on $\tau$ are similar to those of \cite{Rohe2011} that ensure the eigengap of $\bar L\bar L$ and the smallest nonzero eigenvalues of  $\bar L$ cannot be too small.  Hence we omit the discussion here for the sake of brevity.

{\bf {Proof of Theorem \ref{Spectralthm3}. }} Since $\zeta(n)$ in Proposition \ref{lemmaconsistent} is summable, by the Borel-Cantelli Lemma, we have that
\begin{align}
\|L^{}L^{}-\bar L^{}\bar L^{}\|_F=o\left( \frac{\log n}{\tau^2n^{1/2}}G^{1/2}(\chi,N,3)(1-\chi)^{-1/2}\right),\quad a.s. ,
\end{align}
where  $\|\cdot\|_F$ represents the Frobenius norm.
Thus we have that \begin{align}\label{April-106}
\max_{1\leq i\leq n}|\iota_i^{}-\bar \iota_i^{}|=o\left( \frac{\log n}{\tau^2n^{1/2}}G^{1/2}(\chi,N,3)(1-\chi)^{-1/2}\right),\quad a.s.
\end{align}
Define $S_n=[\lambda_{k_n}^2/2,2]$, and
\begin{align}
\delta_n=\inf\{|l-s|; l\in \lambda(\bar L^{}\bar L^{}),l\not \in S_n,s\in S_n\},\\
\delta'_n=\inf\{|l-s|; l\in \lambda_{S_n}(\bar L^{}\bar L^{}),s\not \in S_n\}.
\end{align}
Then $\delta_n=\delta'_n=\lambda_{k_n}^2/2$. The quantity $\delta_n$ measures the distance between spectrum (eigenvalue) of $\bar L\bar L$ outside $S_n$ and $S_n$. $\delta'_n$ measures the how the $S_n$ separates the eigenvalues of $\bar L\bar L$. They are needed for the application of the Davis-Kahan Theorem (\cite{Rohe2011}). Together with \eqref{April-106}, by assumption, the number of elements in $\lambda_{S_n}(\bar L^{}\bar L^{})$ will be equal to the number of elements in $\lambda_{S_n}( L^{} L^{})$. By the definition of $\mathcal M$, as well as the properties (a) \& (b), and the definition of a centroid, we have that
\begin{align}\label{Aug-59}
|\mathcal M|\leq 8P_n\|V_2-ZV_1O\|_F^2=o\left(\frac{P_n\log^2 n}{\lambda_{k_n}^4\tau_n^4n}G^{}(\chi,N,3)(1-\chi)^{-1}\right),\quad a.s.
\end{align}
The last equality follows from the Davis-Kahan theorem, the results of Theorem 2.2, Lemmas 3.1 and 3.2, and the proof of Theorem 3.1 in  \cite{Rohe2011}. \hfill$\Box$

Theorem \ref{Spectralthm3} shows the consistency of the spectral clustering algorithm for composite SBM under regularity conditions. With $\chi$ gets closer to $1$, the dependence between  edge variables becomes stronger and the theoretically guaranteed convergence rate deteriorates. On the other hand, the requirement that $\tau_n^2>M/\log n$ is almost as restrictive as the requirement of at least linearly growing expected degree for all nodes. The following proposition is key to study the tail probability of $\|LL-\bar L\bar L\|_F$ (Proposition \ref{lemmaconsistent}), which controls the mis-clustering rate in Theorem \ref{Spectralthm3}. The proposition is of general interest. It reveals the covariance structure of the latent time-order graph.

\begin{proposition}\label{newprop6}
	Let $C>0$ be a sufficiently large constant. Consider the 
	latent time-order graph with fixed memory parameter $l<\infty$. Let $h\geq 1$ be an integer.
	{\color{black} Define $\{\tilde B_j\}_{1\leq j\leq N}$ to be $\{B_j-\E(B_j|\xi)\}_{1\leq j\leq N}$.}
	Then for any $2l$ different integers $\{i_s,1\leq s\leq 2h\}\in [N]^{2h}$
	and for a sufficiently large positive constant $C$ which may depend on $h$, we have that (a)
	\begin{align}
	|Cov(\Pi_{u=1}^h\tilde B_{i_u}, \Pi_{v=h+1}^{2h}\tilde B_{i_v}|\xi)|=O(\chi^{\Lambda(i_s,1\leq s\leq 2h)/2}(1-\chi)^{-1}),
	\end{align}
	where $\Lambda(i_s,1\leq s\leq 2h)=\min_{1\leq s\leq h}(\min_{h+1\leq u\leq 2h}|i_s-i_u|).$
	In addition, we have (b)
	\begin{align}
	|Cov(\Pi_{u=1}^h\tilde B_{i_u}, \Pi_{v=h+1}^{2h}\tilde B_{i_v}|\xi)|= O(\chi^{\max_{1\leq s\leq 2l}{\iota(i_s)}/2}(1-\chi)^{-1}),
	\end{align}
	where $\iota(i_s)=\min_{1\leq j\leq 2l,j\neq s}(|i_s-i_j|)$,  $Cov(A,B|\xi)=\E(AB|\xi)-\E(A|\xi)\E(B|\xi)$.
\end{proposition}
\begin{remark} 
	Proposition \ref{newprop6} shows that for the composite SBM, the upper bound of the covariance between the product of two groups of  ordered edge variables, $\{B_{j_s},1\leq s\leq l\}$ and  $\{B_{j_s},l+1\leq s\leq 2l\}$, is determined by (a) $\min_{1\leq s\leq l}(\min_{l+1\leq u\leq 2l}|i_s-i_u|)$ and (b) $\max_{1\leq s\leq 2l}\iota(i_s)$. 
	$Cov(\Pi_{u=1}^h\tilde B_{i_u}, \Pi_{v=h+1}^{2h}\tilde B_{i_v}|\xi)$ will become smaller when the terms described by (a) and (b) become larger. The term in (a) is large if the two groups of labels are far away from each other, i.e. , the smallest distance between two labels, one from the $j_s, 1\leq s\leq l$ and the other from $j_s, l+1\leq s\leq 2l$ is large. The quantity (b) is large if there is a label far way from all other labels.
\end{remark}
Recently, many complex models have been proposed based on the SBM to capture additional and important graph structure. For instance, the general SBM proposed in {\color{black} \cite{Caili2015} allows for a portion of arbitrary outliers, where the majority of nodes are generated from a fixed SBM. As a comparison, all nodes from the composite SBM in this paper differ from the SBM when edge variables are conditionally dependent on the latent membership. Another prominent model that can generate arbitrary degree inhomogeneity is the degree corrected stochastic block model (DC-SBM) (see \cite{KN2011}). For this model, consistency of community detection has been studied (see for example \cite{zhao2012consistency}), and corresponding spectral clustering algorithms have been proposed (for example see {\cite{QinRohe2013}}). Also, SBM has been generalized to a mixed membership (for example \cite{Airoldi2013}), and the $K$-median approach (\cite{ZhangL2014}). A tensor approach (\cite{Anandkumar2014}) have been proposed to address the mixed-membership. In this paper, we have built up a general framework for non-exchangeable graphs, and investigate the spectral clustering algorithm for composite SBM in detail.}

\begin{remark}
The key concepts of a ``composite graph'' and composite SBM in Sections \ref{Chap:graphon} and \ref{Chap:spectral}  are  closely related to {\color{black} notion of} composite likelihood. Composite likelihood inference is a popular and successful tool for statistical research when the joint likelihood is hard to evaluate, see \cite{VaRe2011} among others for a comprehensive review. In the literature of network analysis, the idea of analyzing pseudo or approximate likelihood has been proposed to tackle the complex and computational-infeasible joint likelihood of graph models, see for example (\cite{Amini2013}, \cite{Amini2017} and \cite{BiChCh2013} among others.)
\end{remark}

\begin{remark}\label{New.remarkspec6}
We shall assume that $k_n\equiv k$. Assume that there exists $\zeta>0$ such that $|n_{k}|\geq \lf \zeta  n\rf$.
Straightforward calculations show that $G(\chi,N,3)$ is of the order $(1-\chi^{1/2}_n)^{-4}$. As a result, the condition $n^{-1/2}(\log n)^2 G^{1/2}(\chi,N,3)(1-\chi)^{-1/2}=O(\lambda_{k}^2)$ reduces to \begin{align}(1-\chi_n)(1-\chi_n^{1/2})^{-4}\log ^4n=O(n),\label{consistent_condition} \end{align} which yields the weak consistency of clustering in the sense of \cite{zhao2012consistency}, i.e, the mis-clustering rate in Theorem \ref{Spectralthm3} is therefore simplified to  $|\mathcal M|=o(1)$. A straightforward calculation shows that a sufficient condition for \eqref{consistent_condition} is that
\begin{align}
1-\chi_n=\Omega(n^{-1/5}\log n^{-4/5}).
\end{align}
As a comparison, Remark \ref{Consiscondition} shows that the estimation error of Theorem \ref{cconsistent} is negligible if  $1-\chi_n=\Omega(n^{-1/2})$.
\end{remark}

\subsection{Simulation Study of the Mis-clustering Rate}\label{simuspec}
In this section, we examine the performance of the spectral clustering algorithm used on the composite SBM. 
We consider two simulation scenarios: two groups and three groups where the $n$ nodes are partitioned into. A latent order $\omega_1(\{\cdot,\cdot\})$ representing a strong dependence and a latent order $\omega_2(\{\cdot,\cdot\})$ representing a weak dependence are considered. The orders (and these are literal orderings, not orders of magnitude) are constructed in a  way such that their corresponding marginal  edge variables linkage probabilities are the identical. Their forms are deferred to and discussed in detail in Sections \ref{sec4.1} and \ref{sec4.2}. The other detailed parameters of the considered two and  three group composite SBMs can be found in Section \ref{simuspecpara} of the supplemental material. 
\begin{figure}
    \centering
    \includegraphics[width=15cm,height=6cm]{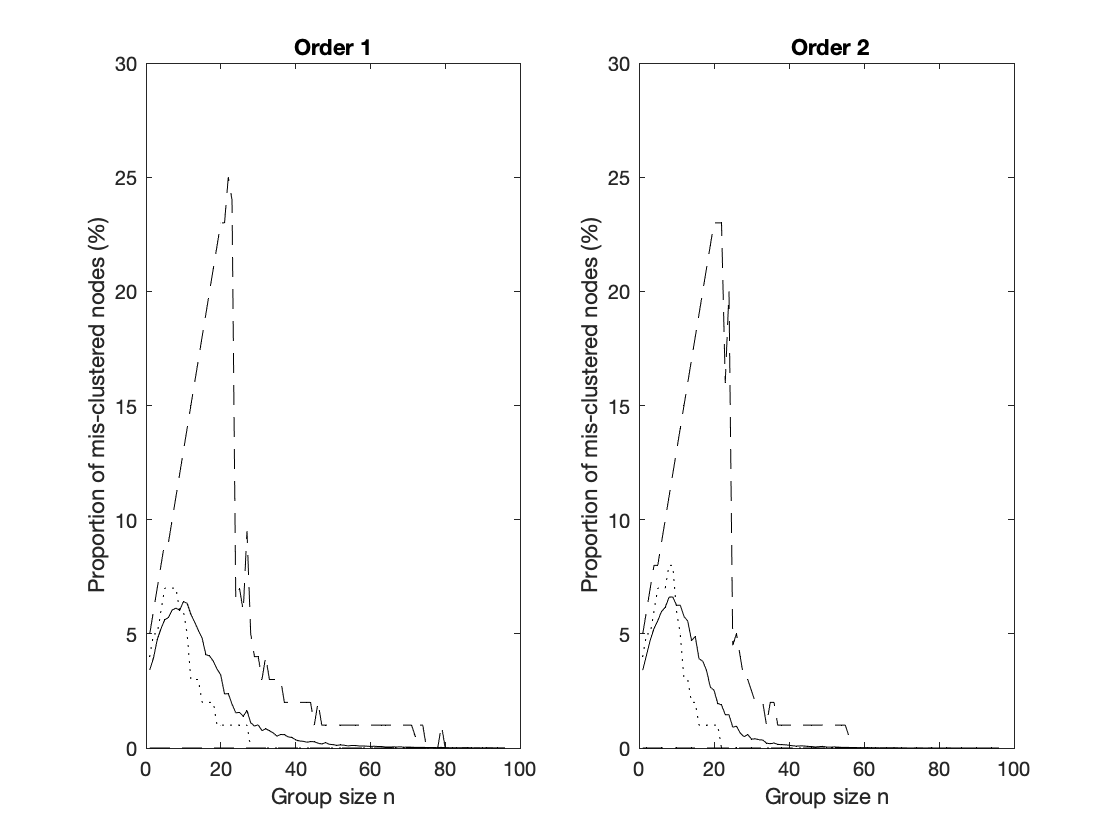}
    \caption{Exploring the mis-clustering under different latent orders for two groups. Left: the number of mis-clustered nodes under order $\omega_1(\{\cdot,\cdot\})$. Right: the number of mis-clustered nodes under order $\omega_2(\{\cdot,\cdot\})$. We display the mean  (solid line), median (dotted line) and upper/lower frequency band (dashed line) of the mis-clustering number over 1000 simulations.  }\label{two_groups}
\end{figure}

The simulation results are displayed in Figures \ref{two_groups} for two groups and which show that there are fewer mis-clustered  nodes under order $\omega_2(\{\cdot,\cdot\})$ than under $\omega_1(\{\cdot,\cdot\})$.  This  reinforces the message of Section \ref{Chap:spectral}, namely that  stronger conditional dependence between  edge variables, as introduced by the order $\omega_1(\{\cdot,\cdot\})$, tends to increase the  mis-clustering rate.   
These simulation results also support the consistency of the spectral clustering algorithm for the composite SBM.  
The clustering results for the three group scenario are similar to that of the two group cases and have been shown in Figure \ref{three_groups} of the supplemental material.
\section{Example: Marginally Edge Constant Latent Time--Order Graph}\label{simpleedge}

In this section, we study the effect of the posited  conditional dependence by studying the given model of the marginally edge constant latent time-order graph sequence model (MECLTG). The MECLTG sequence model is defined as a composite graphon model with $f(\cdot,\cdot)\equiv c_n$ where $c_n$ is a function of $n$, {and $f(\cdot)$ is defined in \eqref{Compositegraphon}}. 
For each fixed $n$, $\p(A_{i,j}=1)=\p(A_{k,l}=1)=c_n$ for every $(i,j),(k,l),i\neq j,k\neq l$. 
It is also a composite SBM with only one group, (as the SBM with one group corresponds to the Erd\"{o}s-R\'{e}nyi model).  Meanwhile, edge variables are correlated with respect to the latent order $\omega_n(\{\cdot,\cdot\})$. Via studying the MECLTG, we can investigate the effect of the (conditional) dependence separately from the effect of inhomogeneous (marginal)  edge variables' linkage probabilities. In the following arguments, for simplicity we omit the subscript $n$ of $\omega_n(\{\cdot,\cdot\})$, $p_{1,n}, p_{0,n}$ and $p_n$ if no confusion arises.

Consider the Markov process of MECLTG, which we write as $CG(V,\omega,p_0,p_1)$:
\begin{align}\label{Model1}
\p(B_i|B_{i-1})=(p_1^{B_i}(1-p_1)^{1-B_i})^{B_{i-1}}(p_0^{B_i}(1-p_0)^{1-B_i})^{1-B_{i-1}},
\\ \text{where}\ \ B_i=A_{\omega^{-1}(i)}.\notag
\end{align}
From \eqref{Model1},  if $B_{i-1}=1$ then $B_i$ is distributed as $Bernoulli(p_1)$, otherwise it is distributed as $Bernoulli(p_0)$.
By the fundamental theorem of Markov Chains, the  edge variables of $CG(V,\omega,p_0,p_1)$ have a limiting distribution
\begin{align}\label{stationp}
p:=\p(B_{\infty}=1)=p_0/(1+p_0-p_1)\in[p_0\wedge p_1, p_0\vee p_1].
\end{align}
We then consider the stationary scenario, i.e.,
\begin{align}\label{stationp1}\p(B_j)=p~~~~\forall 1\leq j\leq N.\end{align} This is because when the total number of the  edge variables is large,
the majority of the  edge variables of $CG(V,\omega,p_0,p_1)$ have marginal linkage probabilities close to $p$. 
\begin{definition}\label{defofCMMG}
We say that the graph $CG(V,\omega,p_0,p_1)$ is a first order homogeneous MECLTG with vertices $V$, driven by the order $\omega(\{\cdot,\cdot\})$ if \eqref{Model1} and \eqref{stationp1} hold. 
\end{definition}
 Notice that when $p_0=p_1=p$, then the  latent structure is not active. As a result, $CG(V,\omega,p_0,p_1)$ reduces to the standard Erd\H{o}s-R\'{e}nyi graph $G(|V|,p)$. The following corollary explicitly calculates the conditional probability of $B_i$ given $B_{i-k}$ for $k\geq 2$:
\begin{corol}\label{coroldependent}
	Consider $CG(V,\omega,p_0,p_1)$. Define $\p_k(a|b)=\p(B_j=a|B_{j-k}=b)$ for $a,b\in \{0,1\}$. Then we have that
	\begin{align*}
	&\p_k(1|0)=\frac{p_0(1-(p_1-p_0)^k)}{1-p_1+p_0};&\p_k(1|1)=\frac{p_0+(1-p_1)(p_1-p_0)^k}{1-p_1+p_0},\notag\\
	&\p_k(0|1)=\frac{(1-p_1)(1-(p_1-p_0)^k)}{1-p_1+p_0};&\p_k(0|0)=\frac{(1-p_1)+p_0(p_1-p_0)^k}{1-p_1+p_0}.
	\end{align*}
	{If we let $p=\frac{p_0}{1-p_1+p_0}$}, {then} we have that 
	\begin{align}
	&\p_k(1|0)=p-p(p_1-p_0)^k,\ \ \p_k(1|1)=p+(1-p)(p_1-p_0)^k\notag\\
	&\p_k(0|1)=(1-p)-(1-p)(p_1-p_0)^k,\ \ \p_k(0|0)=(1-p)+p(p_1-p_0)^k.
	\end{align}
\end{corol}
The results show that the dependence between edge variables, or equivalently
$\Delta_n(k)$ defined in \eqref{dependencemeasure} decays at the geometric rate $(p_1-p_0)^k$. We discuss the phase transition of MECLTG in Section \ref{gai} of the online supplementary
material. In the remaining of the paper, we focus on the degree
distribution of the MECLTG.
\subsection{Degree Distribution}\label{df}
  In previous sections we find that the {ordering}  $\omega_n(\{\cdot,\cdot\})$ has an asymptotically negligible impact on graphon estimation and on community detection under weak dependence (between  edge variables), i.e. $\chi=\chi_n\in (0,1)$. 
  Any $\omega_n(\{\cdot,\cdot\})$ will yield a consistent estimator of the graphon or detection for communities  when $0\leq \chi=\chi_n<1$, where $\chi$ is defined in Proposition \ref{dependenceprop}.
  In this subsection, by investigating simple examples,  we shall see that i) different {ordering} $\omega_n(\{\cdot,\cdot\})$ have different impact on the network structure when $\chi_n\rightarrow 1$, so the impact of missing information of $\omega_n(\{\cdot,\cdot\})$ is no longer asymptotic negligible; and ii) our model is flexible enough to produce networks both { with and without a heavy-tailed degree distribution}.  
  The idea is that when $\p(B_{i}=1|B_{i-1}=1)\rightarrow 1$, we can design  latent orders such that incident edge variables (e.g., $A_{ij}$ and $A_{ik}$) are strongly correlated (or weakly correlated), and hence their summation, or corresponding degrees, cannot (can) be well approximated by sums of independent Bernoulli random variables.

  To illustrate this, assume that $p_0=\lambda_0/n$ and also $p_1\geq p_0$. Recall the homogeneous probability $p=\frac{p_0}{1-p_1+p_0}$. Obviously, $p$ and $p_0$ are of the same order if either $p_1$ is a constant or goes to $0$.  If $p_1=1-n^{-1}g_n$ for a sequence of positive real numbers $g_n$, then $p$ still goes to $0$ as long as $g_n\rightarrow \infty$.

\subsubsection{Examples of MECLTG with heavy tail degree distribution }\label{sec4.1}

Let $\varpi_1(i,j)=n(i-1)-i(i-1)/2+j-i$ for $1\leq i<j\leq n$, and $\omega_1(\{i,j\})=\varpi(i\wedge j, i\vee j)$, $i\neq j$.
Consider the first order  homogeneous process MECLTG $CG(V,\omega_1,p_0,p_1)$.
In particular we choose the ordering $\omega_1(\{i,j\})$ where the  edge variables are generated as follows: 
\begin{align*}
{A_{1,2}, A_{1,3},...,A_{1,n},
A_{2,3},...,A_{2,n}, A_{3,4},..,A_{3,n},...,A_{n-1,n}.}
\end{align*}
 In other words, a node generates its  edge variables after all the node labels before it have generated their  edge variables. 
 We shall see that the considered ordering is able to generate a heavy-tailed degree distribution if we set $p_1\rightarrow 1$, since a connected  edge variables $A_{a,b}$ will lead to a high chance that the next  edge variable $A_{a,b+1}$ is connected, where the two  edge variables have the same vertex. In such a way our model is able to produce a larger number of high degree nodes than the Erd\"os-R\'enyi model.
 
We now study the empirical degree distribution $n^{-1}\sum_{i=1}^nI(d_i=k)$ for $1\leq k\leq n$, where $d_i$ is the degree of node $i$. When the nodes have a homogeneous degree distribution, for example in the Erd\"{o}s-R\'{e}nyi graph, $n^{-1}\sum_{i=1}^nI(d_i=k)$ is an unbiased estimator of $\E d_i$, see for example \cite{Ouadah2015}. Meanwhile, inhomogeneity introduced by strong dependence will distort the empirical degree distribution,
 i.e., as we shall show among a wide range of $k$, the expectation of $n^{-1}\sum_{i=1}^nI(d_i=k)$ of the graph decays with $k$ {at} a polynomial rate. In this way, the graph displays the power law degree distribution. Different from {the Erd\"{o}s-R\'{e}nyi model}, {$\p(d_i=k)$ for MECLTG is heterogeneous in $i$ instead of remaining constant in $i$}. 
\begin{theorem}\label{degree}(Heavy-tailed Degree Distribution)
Consider the first order  homogeneous MECLTG $CG(V,\omega_1,p_0,p_1)$ with $|V|=n$. Suppose $p_0=\frac{\lambda_0}{n}$ with $\lambda_0\geq1$, 
and $p_1=1-\lambda_1n^{-c}$, $0<c<1/2$.
For any $\gamma>1$, $\mu>0$, define $M_\gamma: \sum_{k=1}^n M_\gamma\frac{1}{k^\gamma}=1$, and $M_{\gamma,\mu}: \sum_{k=1}^n M_{\gamma,\mu}\frac{1}{k^\gamma}\exp(-\mu k)=1$.
Let $A_{n,\gamma}=\{k:n^{-1}\sum_{i=1}^n\p(d_i=k)\geq M_\gamma k^{-\gamma}\}$, $B_{n,\gamma,\mu}=\{k:n^{-1}\sum_{i=1}^n\p(d_i=k)\geq M_{\gamma,\mu} k^{-\gamma}\exp(-\mu k)\}$. Then there exist $a_0,b_0,c_0,d_0>0$ (which may depend on $\gamma$), such that
\begin{align}
\left\{k:\lf a_0n^{\frac{2c}{1+\gamma}}\rf\leq k\leq \lf b_0n^c\log^{}n\rf\right\}\subset A_{n,\gamma},\notag\\
\left\{k:\lf c_0\log n\rf\leq k\leq \lf d_0n^{}\rf\right\}\subset B_{n,\gamma,\mu}.
\end{align}
\end{theorem}
{\it Proof.} See section \ref{Sec:4} of Appendix .\hfill $\Box$

 Theorem \ref{degree} shows that, within a wide range of values of  $k$, the tail of the distribution of the degrees of the MECLTG $CG(V,\omega_1,p_0,p_1)$ model behaves similarly to the power law distribution (or to {a} power law {degree} distribution with exponential cutoff, see \cite{Newman2005} ). Consider the usual Erd\"{o}s-R\'{e}nyi graph $G(V,p)$, where $p=\frac{p_0}{1-p_1+p_0}$, so that the marginal linkage probabilities of  edge variables are the same as the first order  homogeneous MECLTG $CG(V,\omega_1,p_0,p_1)$. Let $C_n=\{k: k\geq g(n)\}$ where $g(n)\rightarrow \infty$ arbitrarily slowly. By proposition \ref{possionapprox} (a Poisson approximation) in the supplementary supplement and the large deviation theorem (see the proof of Lemma \ref{homoorder} in the supplementary material), it follows that there exist constants $c,d$ such that both $C_n\cap A_{n,\gamma}$ and $C_n\cap B_{n,\gamma,\mu}$ are subset of $\{k:\lf a n^c\rf\leq k\leq \lf b n^c\rf\}$ when $n$ is large enough.
Thus, the first order  homogeneous MECLTG $CG(V,\omega_1,p_0,p_1)$ has much larger $|C_n\cap A_{n,\gamma}|$ and $|C_n\cap B_{n,\gamma,\mu}|$ than SRG $G(V,p)$.

 \subsubsection{Examples of MECLTG with light-tailed degree distribution}\label{sec4.2}
In this section, we construct a first order  homogeneous MECLTG $CG(V,\omega_2,p_0,p_1)$ which has similar $|C_n\cap A_{n,\gamma}|$ and $|C_n\cap B_{n,\gamma,\mu}|$ to that of Erd\"{o}s-R\'{e}nyi graph $G(V,p=\frac{p_0}{1-p_1+p_0})$. The order $\omega_2$ we consider is $\omega_2(\{i,j\})=\varpi_2(i\wedge j, i\vee j)$, $i\neq j$, where $\varpi_2(i,j)=i+\frac{(2n-(j-i))(j-i-1)}{2}$ for $1\leq i<j\leq n$.  
{In particular, our choices of ordering follow which the  edge variables $A_{i,j}, i<j$ are generated are as follows} 
{\begin{align*}
A_{1,2}, A_{2,3},...,A_{n-1,n},
A_{1,3},...,A_{n-2,n},
A_{1,4},..,A_{n-3,n},
...,
A_{1,n}.
\end{align*}}
 {Observe that the  edge variables are generated in increasing order of $j-i$. Among  edge variables with equal $j-i$, the  edge variables with smaller $i$ are generated earlier.} We now study the expectation of $n^{-1}\sum_{i=1}^nI(d_i=k)$ to { show the characteristic of $A$}.
 \begin{lemma}\label{homoorder}
 Consider the first order  homogeneous MECLTG Graph\\ $CG(V,\omega_2,p_0,p_1)$ where $p_0=\frac{\lambda_0}{n}, \lambda_0\geq 1$ and $p_1=1-\lambda_1n^{-c}, c\in (0,1/2)$ are defined in Theorem \ref{degree}. Let $Y_n$ follow Poisson($p$) for $p=\frac{p_0}{1-p_1+p_0}$.
Let $d_i$ be the degree of node $i$. Let $g(n)$ be a series of real numbers which diverges but may increase at an arbitrarily slow rate.  The we have for some $\iota>c$, $\iota+c<1$,\begin{align}
 &\text{(i)}\ \sum_{k=0}^\infty|\p(d_i=k)-\p(Y_n=k)|=O(n^{\iota+c-1}),\\
 &\text{ (ii)}\ \p(d_i=k)\leq \exp(-0.5(\iota-c)k\log n) \ \text{for $\lf n^\iota g(n)\rf\leq k\leq n$}. \label{new.58-late2016}
\end{align}
 \end{lemma}
 {\it Proof.}  See supplementary material, Appendix Y.\hfill $\Box$

 Lemma \ref{homoorder} shows that the behavior of the MECLTG $CG(V,\omega_2,p_0,p_1)$ is similar to an Erd\"{o}s-R\'{e}nyi graph $G(n,p)$, in the sense that both of their degree distributions can be mimicked by a Poisson($\frac{\lambda_0}{\lambda_1}n^c$) random variable.  Equation \eqref{new.58-late2016} also shows that the tail of the degree distribution decays very {rapidly}.  Together with Theorem \ref{degree}, we find that simply multiplying a scale parameter to the marginal  edge variable linkage probabilities (for example in \cite{BiCh2009} which models the sparse graphon as $\rho_n f(\cdot,\cdot)$) is not able to capture the heteroskedasticity in the probability of linkage, in the sense that $\frac{\max_{i,j} \p(A_{ij}=1)}{\min_{i,j} \p(A_{ij}=1)}$ stays unchanged under this  parameterization. This shows the greater flexibility and rich structure of our model class.
  We illustrate this property in Figure \ref{network} of the supplementary material.
  We discuss the images from left to right in Figure \ref{network}. Figure  \ref{network}  shows typical graphs generated from MECLTGs $CG(V,\omega_1,p_0,p_1)$, $CG(V,\omega_2,p_0,p_1)$ and Erd\H{o}s-R\'{e}nyi $G(V,p)$, {respectively} with $p_0=0.01$, $p_1=1-\frac{1}{n^{1/3}}$, $p=\frac{p_0}{1-p_1+p_0}$  with $|V|=100$. 
  	From the figure, we see that the first network is very inhomogeneous: it has the most hubs among the three networks. 
  	The second network is less inhomogeneous than the first network, but is more inhomogeneous than the third network. Notice that we construct the three networks in such a way that the marginal connection probability is  $n^{-2/3}$, where $n$ is the size of the network. This is larger than the connectivity threshold $\frac{\log n}{n}$.
  	However, all the three networks in Figure \ref{network} have some isolated nodes just like the models of~\cite{BCJHV2017}. For the third network, this is because the sample size is not large enough, so $100^{-2/3}$ is very close to $\log 100/100$. We observed that the first and second network have more components, which is the price we pay for the inhomogeneity. Since the marginal connection probability in our experiment is controlled, the expected total edges of the three networks are fixed. As a result, the structure with more hubs will also tend to have more small degree nodes, and also more isolated nodes. The  edge variables are distributed according to the dependence structures $\omega_1$ and $\omega_2$ in the first and second networks,  and purely randomly distributed in the third network.

\begin{remark}\label{Newremark5}
A  familiar model for networks with power law degree distributions is the preferential attachment (PA) model, where the network is growing {sequentially} node by node. In PA, a node can not affect the relationship among earlier nodes. This shares some features with out model.  Thus, the generating order (or history) of the  edge variables of PA could be written as
  {\begin{align*}
  A_{1,2},
  A_{1,3},A_{2,3},
  A_{1,4},A_{2,4},A_{3,4},
  ...,
  A_{1,n},...,A_{n-1,n}.\end{align*}}
with the associated {ordering} $\omega(\{i,j\})=\varpi(i\wedge j, i\vee j)$, where $\varpi(i,j)=\frac{(j-1)(j-2)}{2}+i$, $1\leq i<j\leq n$. The linkage probability of an edge variable is determined by the popularity of its earlier (more popular) verteces. Hence PA is not a latent time-order graph since the required properties fail to hold for any finite $k$.
  \ The well known heavy-tailed degree distribution of PA is contributed by infinite memory, order $\omega(\{,\})$, and inhomogeneous  edge variables linkage probabilities.

{ Recently \cite{Borgs2014} proposed a class of normalized {\it unbounded} graphon model. Given latent positions $(\xi_i=x_i)_{1\leq i\leq n}$, the edge variables are independently connected with probabilities $\p(A_{ij}=1|\xi_i=x_i,\xi_j=x_j)=\min(1, \rho W(x_i, x_j))$, where $\rho$ is the target density and $W$ is a (possibly) unbounded graphon. For detailed definition of $\rho$ and $W$ we refer to \cite{borgs2020consistent}. Due to the inhomogeneous conditional connection probabilities and the unboundedness of $W$, their model is allowed to have a large portion of high degree nodes and therefore the feature of heavy-tailed degree distribution under some circumstances. } 

In contrast to the aforementioned models, the MECLTG {model} $CG(V,\omega_1, p_1,p_0)$ has homogeneous (marginal) edge variable linkage probabilities. Hence the power law degree distribution {is a consequence of the order} $\omega_1(\{\cdot,\cdot\})$ and the strength of  dependence which is determined by $p_0,p_1$. As a comparison, our construction $CG(V,\omega_2, p_1,p_0)$ does not have a power law distribution though it has the same strength of dependence as $CG(V,\omega_1, p_1,p_0)$. The only difference between the two MECLTGs is the order function. The unobserved ordering $\omega_{i},i=1,2$ introduces correlation, this increasing the probability of an edge variable between nodes adjacent in the ordering.
This reveals the complex nature of the latent time-order graph. We display the adjacency matrix of typical $CG(V,\omega_1,p_0,p_1)$ and $CG(V,\omega_2,p_0,p_1)$  with network size $n=100, p_0=\frac{1}{n}$, $p_1=1-\frac{1}{n^{1/3}}$ in left and right panels of Figure \ref{Adjacent}. Notice that $d_i$ is the sum of $i_{th}$ row of the adjacency matrix. The figure shows that $CG(V,\omega_1,p_0,p_1)$ generates  high degree nodes with greater frequency than  $CG(V,\omega_2,p_0,p_1)$.  Recall that degrees are calculated by averaging along rows or columns, whilst diagonal structure does not aggregate to form larger degrees.
\hspace{-4cm}
\begin{figure}[t]
\centering
\begin{center}
\begin{minipage}[t]{0.495\textwidth}
\includegraphics[width=6.2cm]{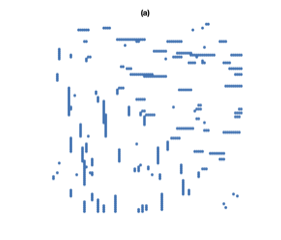}
\end{minipage}
\begin{minipage}[t]{0.495\textwidth}
\includegraphics[width=6.2cm]{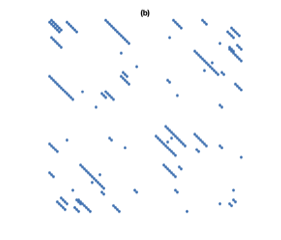}
\end{minipage}
\end{center}
\caption{Adjacency matrix of $CG(V,\omega_1,p_0,p_1)$ (left) and $CG(V,\omega_2,p_0,p_1)$ (right).}
\label{Adjacent}
\end{figure}

\end{remark}

\subsection{Simulation Results of Degree Distributions}\label{simuempricaldegree}

\begin{figure}
 	\centering
 	\begin{subfigure}[b]{0.475\textwidth}
 		\includegraphics[width=\textwidth]{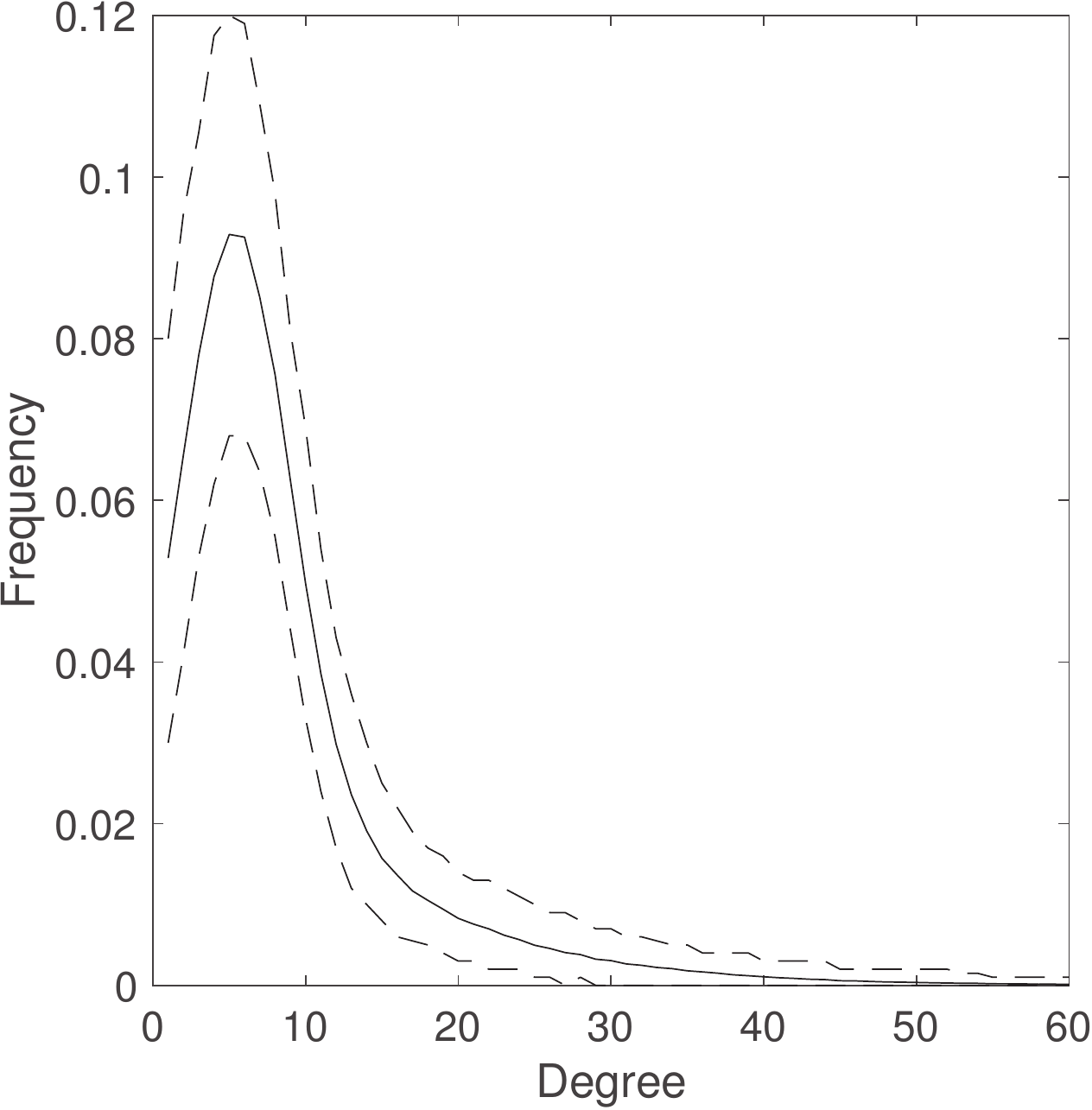}
 		\caption{Mean distribution (solid) and upper/lower frequency band (dashed). }
 	\end{subfigure}
 	\quad
 	~ 
 	\begin{subfigure}[b]{0.46\textwidth}
 		\includegraphics[width=\textwidth]{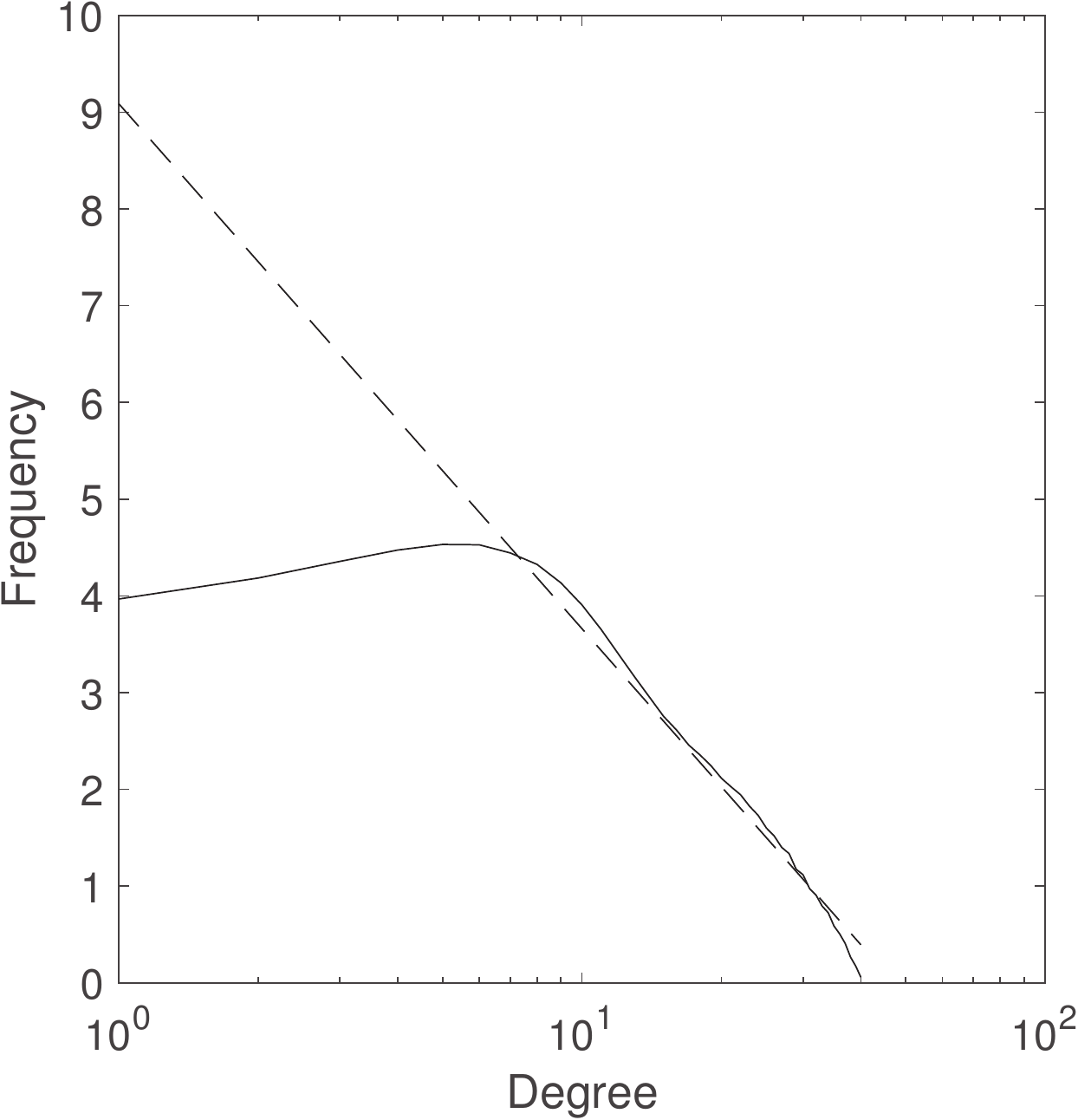}
 		\caption{Powerlaw $\gamma_1=-2.4$(dashed) \& degree distribution (solid).}
 	\end{subfigure}
 	\caption{The degree distribution for $n = 1000$, and  $c = 0.3$.  }\label{gplot1}
 \end{figure}

In this section we generate MECLTGs with $n$ nodes. The orderings are discussed in Section \ref{df}. Let $p_0=\lambda_0/n$, $p_1=1-\lambda_1/n^c$ for $\lambda_0=\lambda_1=1$. We calculate and plot the empirical degree distribution $n^{-1}\sum_{i=1}^nI(d_i=k)$, $k=0,...n$, and fit a power law to the degree distribution as follows: let $k'=\argmax_{0\leq k\leq \sqrt n} n^{-1}\sum_{i=1}^nI(d_i=k)$, and $k''=\max\{k: n^{-1}\sum_{i=1}^nI(d_i=k)>0\}$. We then fit the regression 
\begin{align}\label{fitpowerlaw}\log {\sum_{i=1}^nI(d_i=k)}=\gamma_0+\gamma_1\log k,
\end{align}
 for $k'\leq k\leq k''$ to estimate $\gamma_1$, and use  $\hat \gamma_1$ as the estimate of the power law index. We draw the power line together with the degree distribution in the log scale plot. The empirical distribution is generated by $1000$ replications in each simulation study, and the corresponding 95\% confidence interval is provided by the simulated $0.025$ and $0.975$ quantiles of the simulated samples, respectively. In figure \ref{gplot1} we show the degree distribution for $n=1000$, $c=0.3$  with latent order $\omega_1$.
 We also examine scenarios with various $c's$, $n's$ for latent order $\omega_1$ and $\omega_2$ in figures \ref{gplot2}--\ref{gplot8} shown in the supplementary material. Those figures indicate that $\hat \gamma_1$ increases as $c$ increases.  Also we observe that
 \[\hat \gamma_1(G(V,\omega_1,p_0,p_1))>\hat \gamma_1(G(V,\omega_2,p_0,p_1))\approx \hat \gamma_1(SRG).\]
 
  When modelling the network via either the composite graphon model or the composite SBM, we have demonstrated in Sections \ref{Chap:graphon} and \ref{Chap:spectral} that the usual methods are {still} valid when the dependence is not strong, and the consequences of the inhomogeneous degree pattern is negligible. Our simulation results support this, which coincides with the conclusion in \cite{zhao2012consistency}: comparing with the SBM, the estimation from DC-SBM only improves a little  when the variation of degrees among nodes is not large. On the other hand, when the dependence is strong and the consequences of inhomogeneous degree pattern of network is significant, more advanced approaches are required for modelling and statistically analysing the network.
  
  \section{Discussion}
 As social media data sets, and other types of relational observations (networks) have become prevalent, so unsurprisingly the mathematical treatment of data taking the form of relationships between entities has become increasingly important. The analysis of networks has been the focus of considerable efforts where
the properties of estimators for popular models have now been established. Following on from the understanding of correctly specified parametric models is the usage of non-parametric and incorrectly specified models. For example, our understanding of classical approaches can be found to extend when considering dense exchangeable arrays, see for example~\cite{wolfe2013nonparametric,OlhedeWolfe2012,choi2012stochastic,li2018hierarchical}.

Unfortunately the world contains many data sets that cannot be assumed to be exchangeable, despite how innocuous the assumption may seem, rather like stationarity for time series. For that reason we introduce the composite graphon model, and finite memory latent time-order graphs. By focusing on the latent variables in the model directly we can build a continuum of types of networks that are exchangeable, or strongly non-exchangeable, all tuned explicitly in terms of the dependence strength. This helps us to understand data of this form, and when we can apply regular network tools to novel types of data, and understand the consequences of that choice.

Non-exchangeable networks produce many challenges. The presence of strong powerlaws in the degrees and further heterogeneity in the graphon function itself are still challenging researchers. It is not unreasonable to believe that these features reflect how the network was formed. By assuming that the network formed sequentially we are able to both define a parameter that tunes its degree of exchangeability, and thus we may understand standard tools when applied to such data. Our understanding of this mechanism simultaneously give glimpses into the formation of non-exchangeability, and provide a gray-scale understanding of networks, letting us see how the mechanism allows us to gradually ``dial away'' from exchangeability as a consequence of evolution and growth.

A number of developments have sought to understand greater heterogeneity by modelling edge variables directly rather than relationships between nodes~\cite{crane2018edge}, this allowing a more natural and direct treatment of edge sparsity than some competing models. Others have concerned developing the practical application of work by Kallenberg's constructions~\cite{Kallengberg2005}, such as~\cite{BCJHV2017,borgs2017sparse}. The two key aspects of the latter construction is to use a latent Poisson construction and a latent time. We also used a latent variable which is uniform rather than Poisson. We correlate the latent uniforms directly, and show how the correlation of the latent variables drive the degree of non-exchangeability directly and quantitatively. The advantage of our framework is that it naturally straddles the model space between strong heterogeneity to the standard exchangeable graph model, with a direct tuning of its degree of non-regularity. If the correlation is not too strong, then standard methods apply for estimating the graphon model, rather like in time-series analysis with stationary errors when estimating polynomial trends. As the correlation becomes very strong, the observations exhibit more strong heterogeneity, and standard tools like the stochastic blockmodel approximation of the underlying graphon model will become increasingly problematic.

A number of questions remain unanswered. In parts this falls back to the difficulty of understanding a non-property, which has already haunted both non- stationary and non-linear time series (there are many ways to be non-stationary or non-linear, but only one to be stationary). In parts it falls back to understanding non-exchangeability itself, as one property rather than several real-life observed consequences thereof. By providing this framework, we can better see the limitations of exchangeable models, and how exchangeability can fail to materialize as a consequence of dependence.
  
\section*{Acknowledgements} 
This work was supported by the European Research Council under Grant CoG 2015-682172NETS, within the Seventh European Union Framework Program, and NSFC Young program (No.11901337), and SCO and PJW acknowledge the Isaac Newton Institute for Mathematical Sciences, Cambridge, for support and hospitality during the programme Statistical Network Analysis where work on this paper was
undertaken. This work was therefore also supported by EPSRC grant no. (EP/K032208/1).

\def\theequation{\Alph{section}.\arabic{equation}}
\renewcommand{\thefigure}{\Alph{section}.\arabic{figure}}
\renewcommand{\thelemma}{\Alph{section}.\arabic{lemma}}
\renewcommand{\thetheorem}{\Alph{section}.\arabic{theorem}}
\renewcommand{\theremark}{\Alph{section}.\arabic{remark}}
\renewcommand{\thecorol}{\Alph{section}.\arabic{corol}}
\renewcommand{\theproposition}{\Alph{section}.\arabic{proposition}}
\renewcommand{\thedefinition}{\Alph{section}.\arabic{definition}}
\setcounter{equation}{0}
\setcounter{section}{0}
\renewcommand{\thesection}{\Alph{section}} 
{\begin{center}{\bf Supplemental material for `Tractably Modelling  Dependence in Networks Beyond Exchangeability' }\end{center}}

\begin{abstract}The structure of the supplementary material is organized as follows. Section \ref{simuspecpara} provided arameters in simulation models for spectral clustering algorithm in Section \ref{simuspec} of the main article.
Section \ref{supple-figure} contains additional simulation results for section \ref{simpleedge} of the main article.  Section \ref{gai} contains a discussion of the phase transition for connectivity and giant component of MECLTG. Section \ref{Sec:2} provides the detailed proof of results in Section \ref{Chap:graphon} of the main article. Section \ref{Sec:3} provides the detailed proof of results in Section \ref{Chap:spectral} of the main article. Finally,  Section \ref{Sec:4} provides the detailed proof of results in Section \ref{simpleedge} of the main article. Notice that $N=n(n-1)/2$.
\end{abstract}

\section{Parameters in simulation models for spectral clustering algorithm in Section \ref{simuspec} of the main article}\label{simuspecpara}
For two groups case, we choose the parameter listed in Tables \ref{parameter2-1} and \ref{parameter2-2}, respectively, where we refer the notation $\varrho$ to Section \ref{notationtheta} of the main article. The corresponding results are shown in Figure \ref{two_groups} in the main article.
\begin{table}[htbp]
	\centering
	\caption{The tables for parameters $\varrho^{0,index}_{index}$}
	\begin{tabular}{lrrr}
		\hline
		index & 11    & 12    & 22 \\
		\hline
		value & 0.1   & 0.01  & 0.2 \\
		\hline
	\end{tabular}%
	\label{parameter2-1}%
\end{table}%

\begin{figure}
	\centering
	\includegraphics[width=14cm]{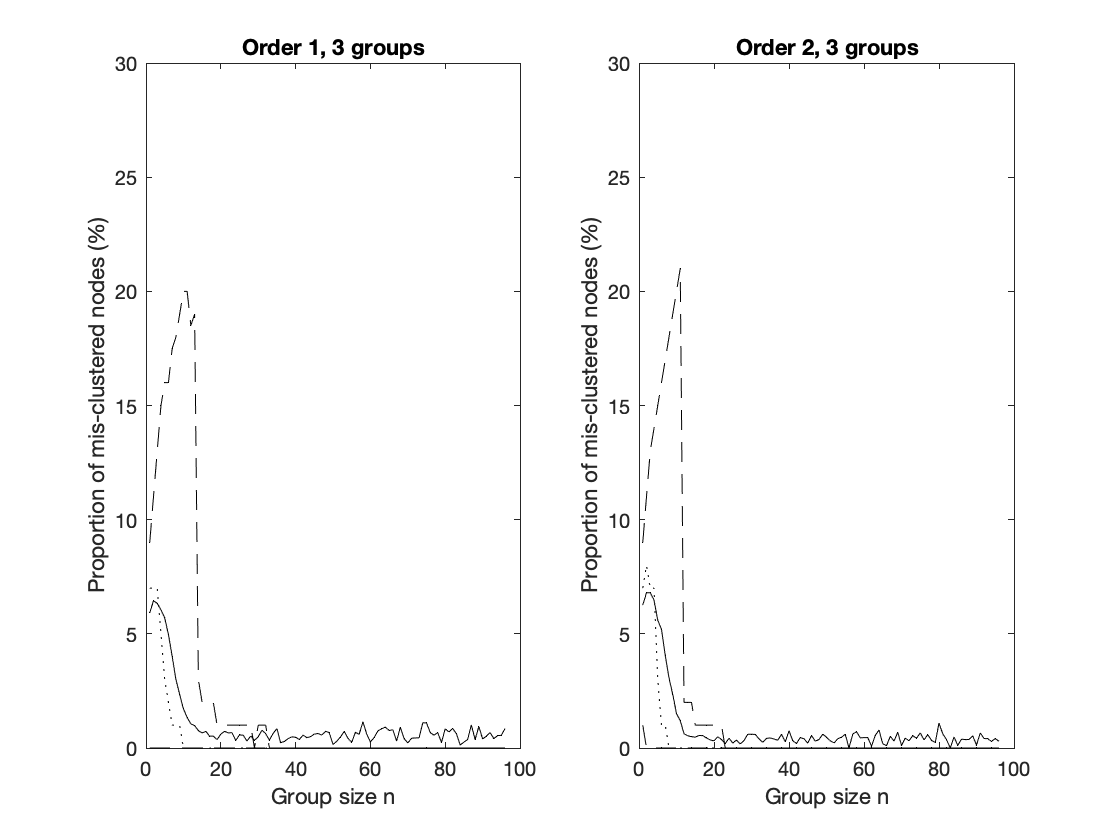}
	\caption{Exploring the mis-clustering  under different latent orders and three groups.  Left: The number of mis-clustered nodes under order $\omega_1(\{\cdot,\cdot\})$.  Right: The number of mis-clustered nodes under order $\omega_2(\{\cdot,\cdot\})$.  We display the Mean  (solid), median (dot) and upper/lower frequency band (dashed) of the mis-clustering number over 1000 simulations. }\label{three_groups}
\end{figure}

\begin{table}[htbp]
	\centering
	\caption{The tables for parameters $\varrho^{1,up}_{down}$}
	\begin{tabular}{c|r|rrr}
		\hline
		&       & \multicolumn{3}{c}{down} \\\hline
		&       & 11    & 12    & 22 \\\hline
		\multirow{3}[0]{*}{up} & 11    & 0.4   & 0.05  & 0.3 \\
		& 12    & 0.3   & 0.1   & 0.1 \\
		& 22    & 0.2   & 0.03  & 0.6 \\\hline
	\end{tabular}%
	\label{parameter2-2}%
\end{table}%
We now show parameters for three groups case in Tables \ref{parameter3-1} and \ref{parameter3-2}, respectively. The clustering result is presented in Figure \ref{three_groups}.

\begin{table}[htbp]
	\centering
	\caption{The tables for parameters $\varrho^{0,index}_{index}$}
	\begin{tabular}{lrrrrrr}
		\hline
		index & 11    & 12    & 13    & 22    & 23    & 33 \\
		\hline
		value & 0.3   & 0.01  & 0.02  & 0.3   & 0.06  & 0.3 \\
		\hline
	\end{tabular}%
	\label{parameter3-1}%
\end{table}%
\begin{table}[htbp]
	\centering
	\caption{The tables for parameters $\varrho^{1,up}_{down}$}
	\begin{tabular}{c|r|rrrrrr}
		\hline
		&       & \multicolumn{6}{c}{down} \\\hline
		&       & 11    & 12    & 13    & 22    & 23    & 33 \\
		\multirow{6}[0]{*}{up} & 11    & 0.5   & 0.02  & 0.05  & 0.3   & 0.02  & 0.04 \\
		& 12    & 0.3   & 0.1   & 0.05  & 0.2   & 0.02  & 0.04 \\
		& 13    & 0.2   & 0.1   & 0.08  & 0.05  & 0.02  & 0.04 \\
		& 22    & 0.15  & 0.02  & 0.02  & 0.6   & 0.02  & 0.04 \\
		& 23    & 0.15  & 0.1   & 0.05  & 0.01  & 0.1   & 0.04 \\
		& 33    & 0.2   & 0.01  & 0.02  & 0.05  & 0.02  & 0.7 \\
		\hline
	\end{tabular}%
	\label{parameter3-2}%
\end{table}%

\section{Simulation results for degree distribution}\label{supple-figure}

\newpage
\begin{sidewaysfigure}
	\centering
	\centering
	\includegraphics[width=19cm,height=11cm]{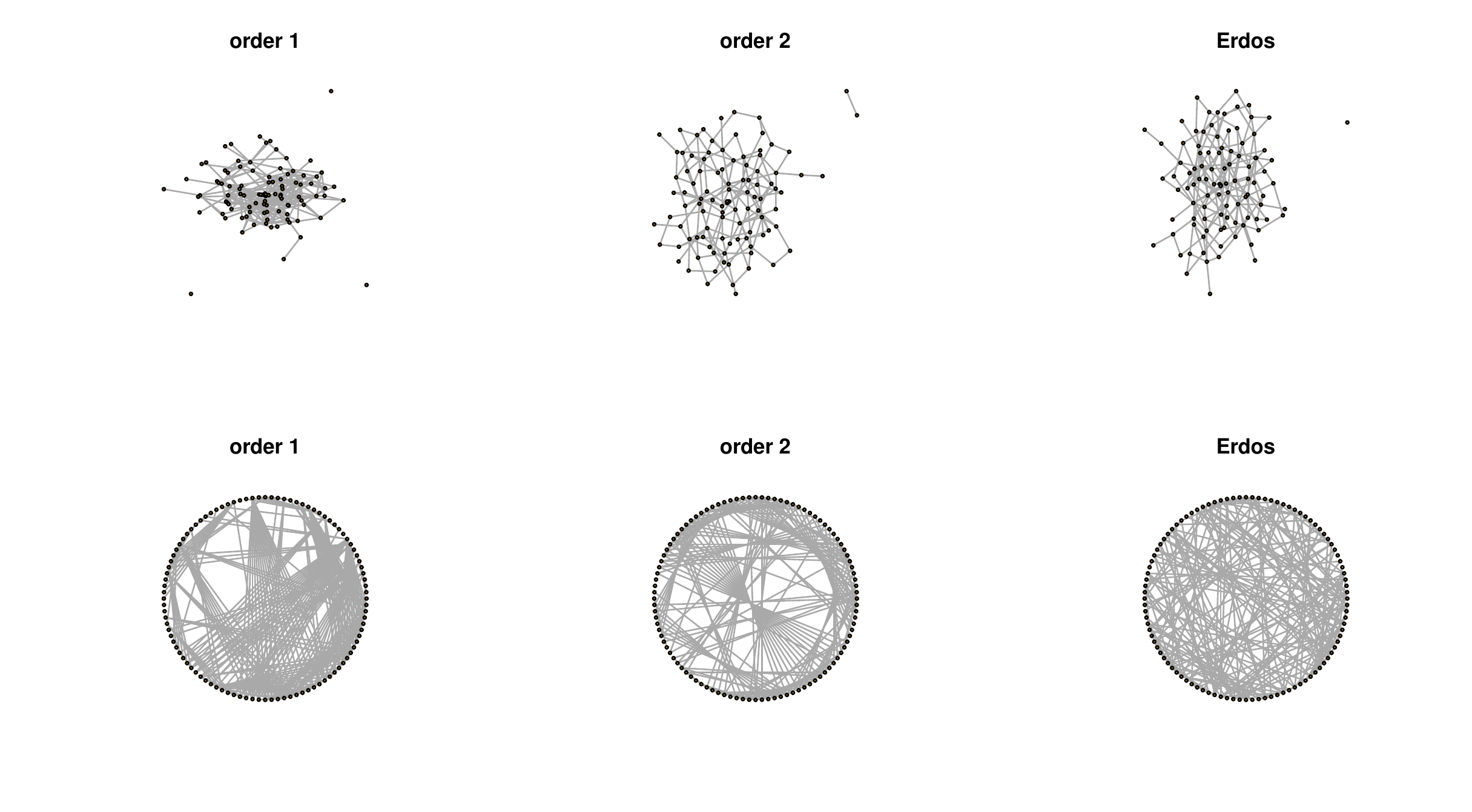}
	\caption{Left and Middle: typical visualization of the size  $n=100$ first order  homogeneous MECLTG with {ordering} $\omega_1$ and $\omega_2$, and
		$p_0=0.01$, $p_1=1-\frac{1}{n^{1/3}}$, {respectively}. Right: Size $n=100$ Erd\H{o}s-R\'{e}nyi Graph with the $p=\frac{p_0}{1-p_1+p_0}$. }\label{network}
\end{sidewaysfigure}
\begin{figure}
	\centering
	\begin{subfigure}[b]{0.475\textwidth}
		\includegraphics[width=\textwidth]{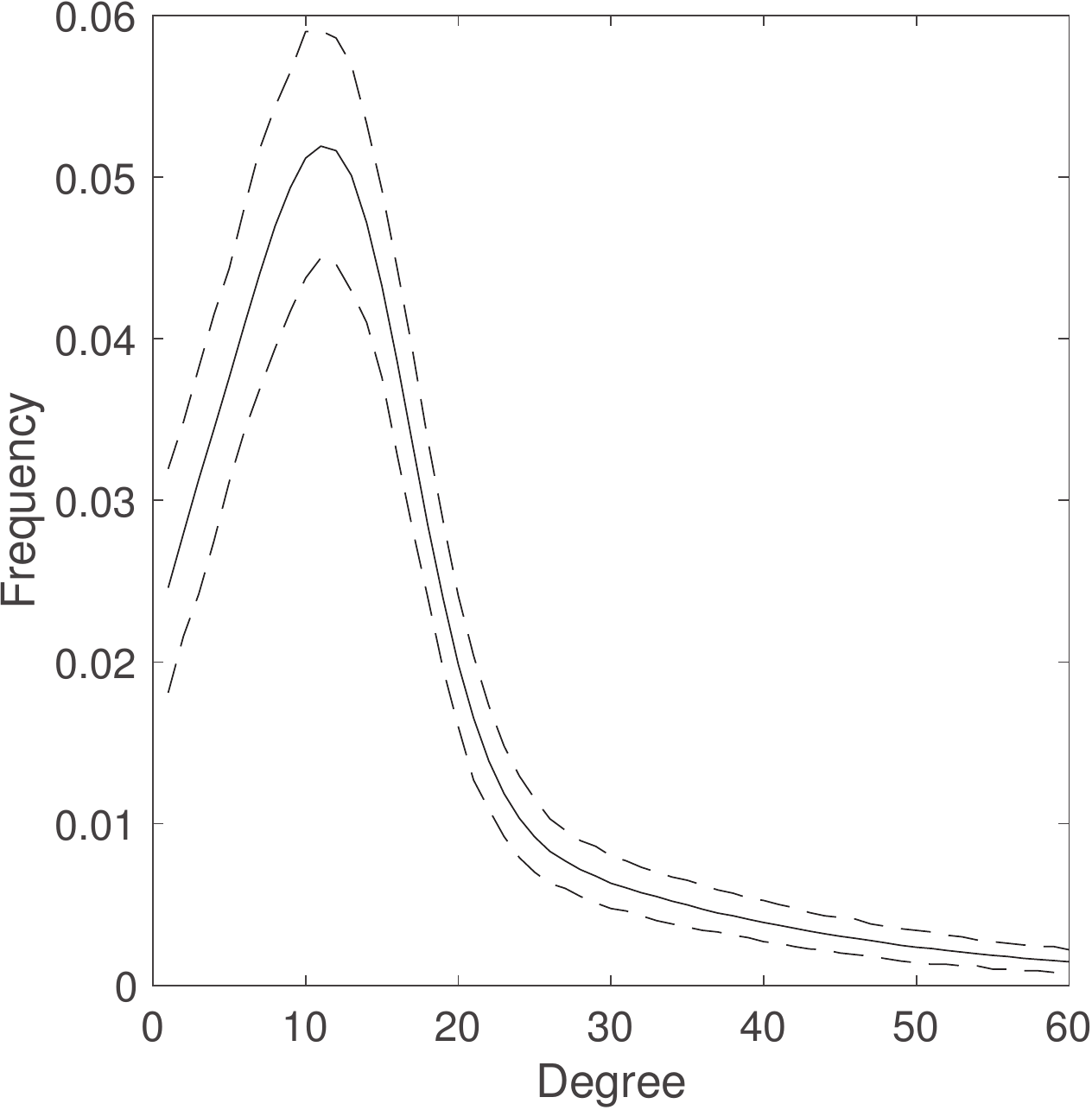}
		\caption{Mean distribution (solid) and upper/lower frequency band (dashed). }
	\end{subfigure}
	\quad
	~ 
	\begin{subfigure}[b]{0.46\textwidth}
		\includegraphics[width=\textwidth]{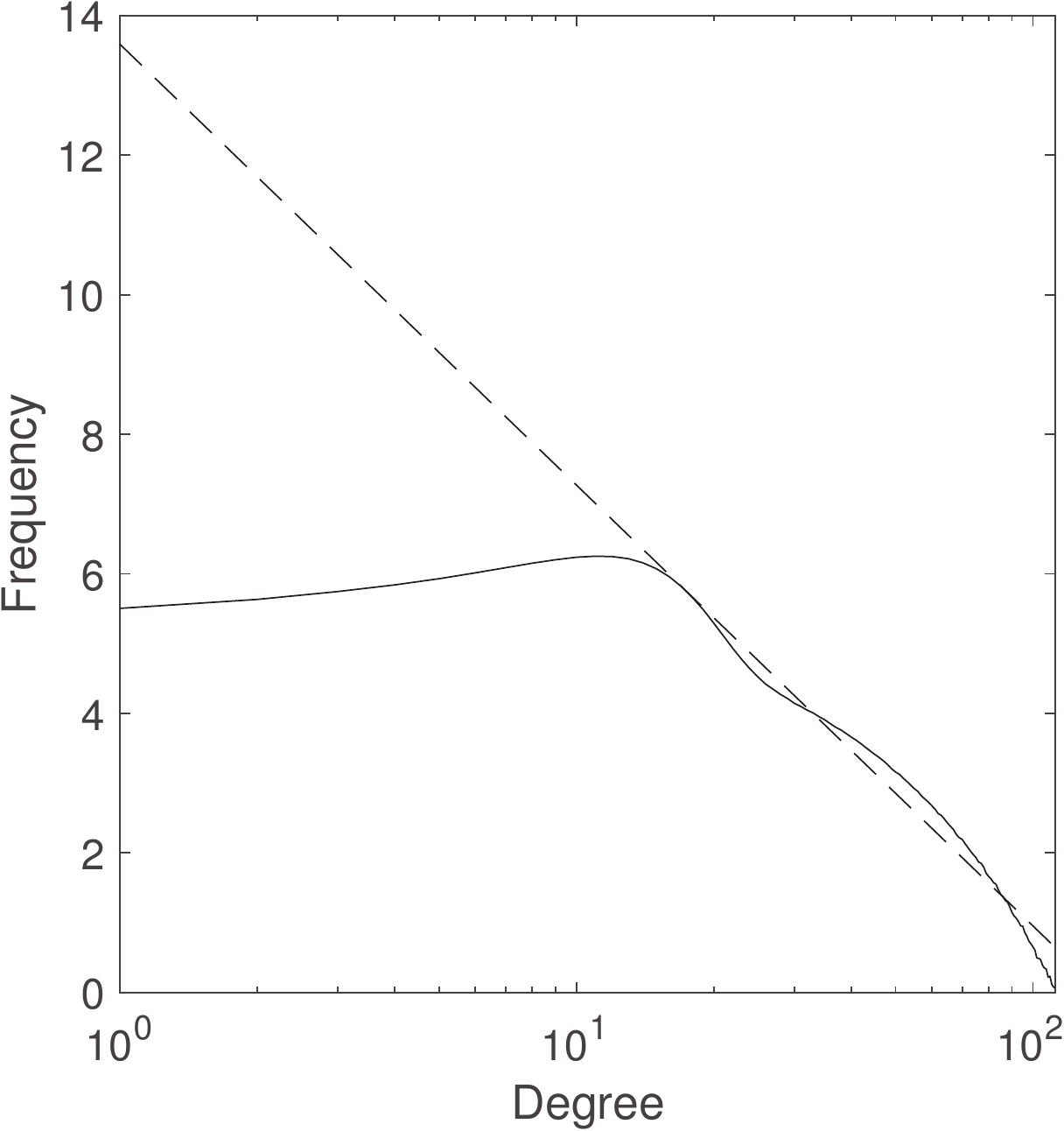}
		\caption{Powerlaw $\gamma_1=-2.74$(dashed) \& degree distribution (solid).}
	\end{subfigure}
	\caption{The degree distribution for $n = 10,000$, and  $c = 0.3$. with latent order $\omega_1$.  }\label{gplot2}
\end{figure}
\begin{figure}
	\centering
	\begin{subfigure}[b]{0.475\textwidth}
		\includegraphics[width=\textwidth]{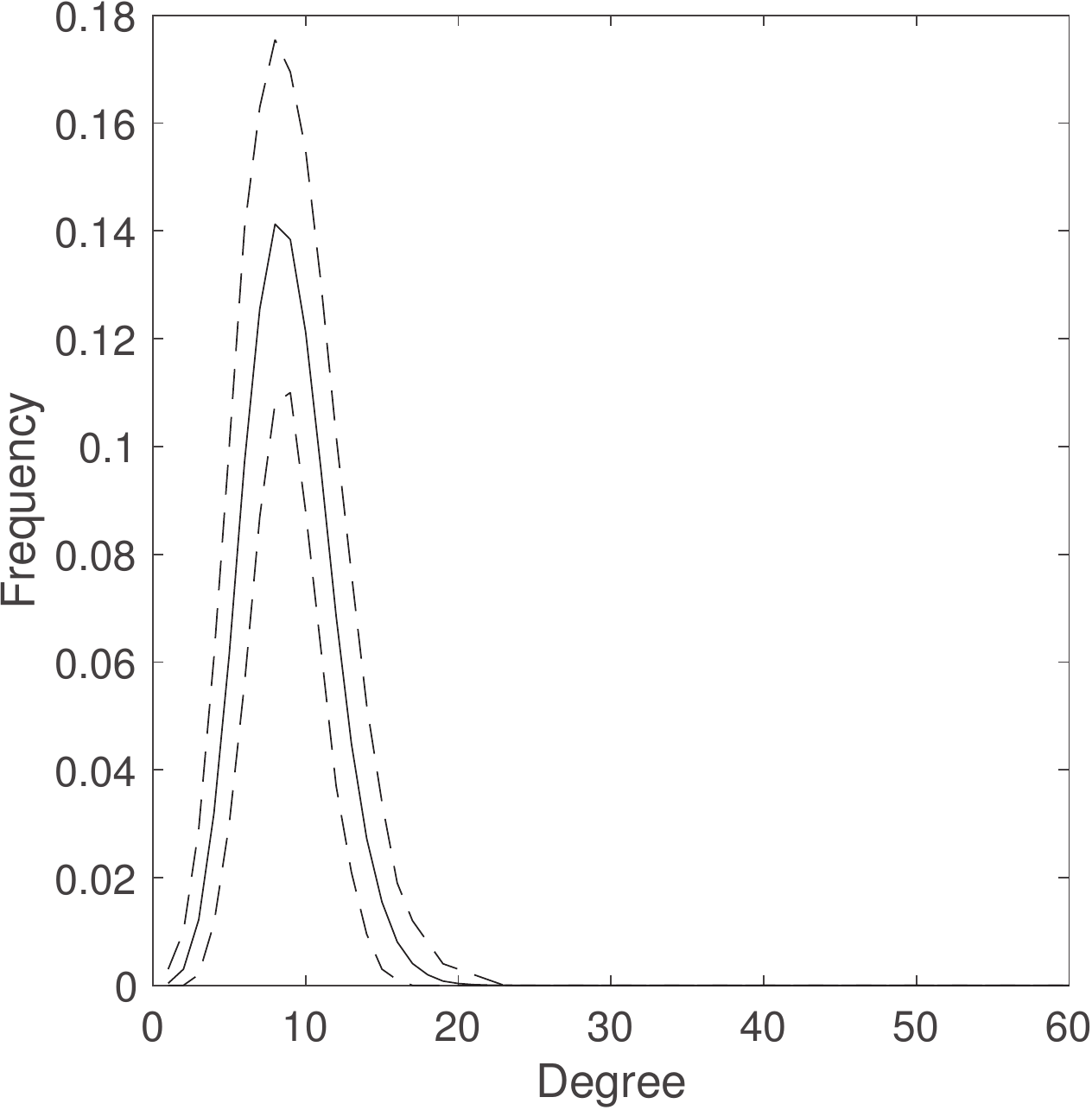}
		\caption{Mean distribution (solid) and upper/lower frequency band (dashed).}
	\end{subfigure}
	\quad
	~ 
	\begin{subfigure}[b]{0.46\textwidth}
		\includegraphics[width=\textwidth]{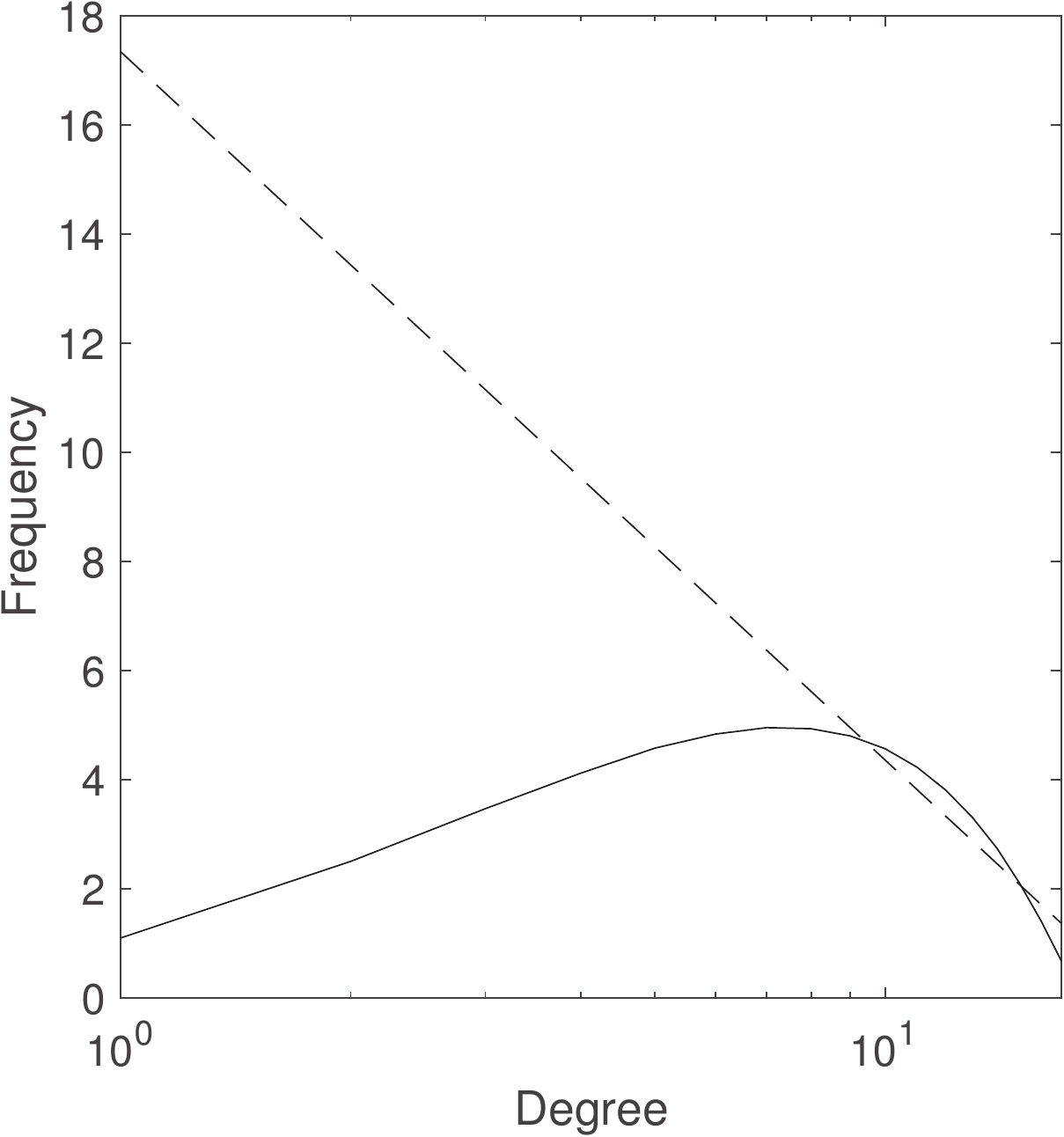}
		\caption{Powerlaw $\gamma_1=-5.6$(dashed) \& degree distribution (solid).}
	\end{subfigure}
	\caption{The degree distribution for $n = 1,000$, and  $c = 0.3$, with latent order $\omega_2$. }\label{gplot3}
\end{figure}
\begin{figure}
	\centering
	\begin{subfigure}[b]{0.475\textwidth}
		\includegraphics[width=\textwidth]{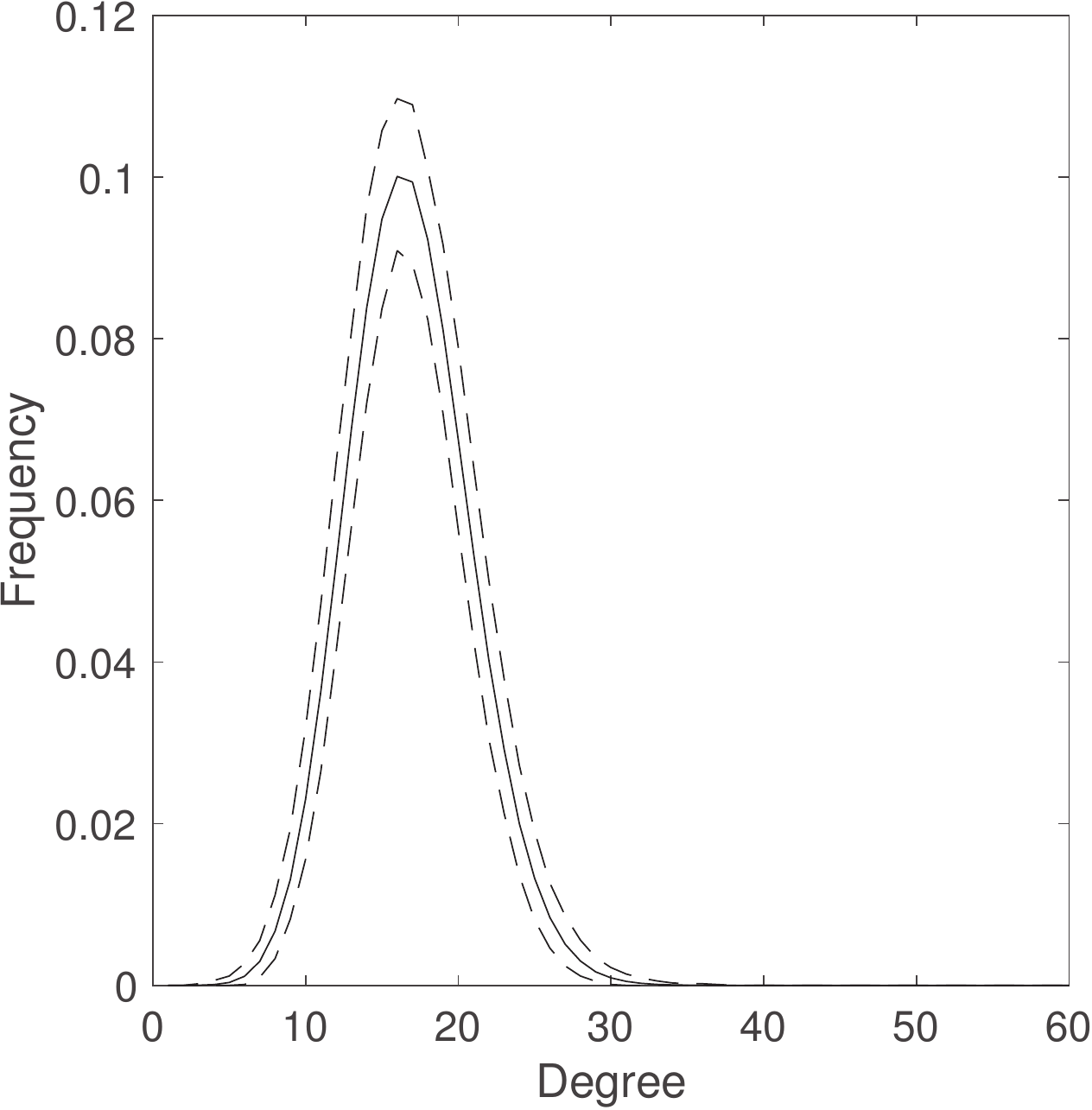}
		\caption{Mean distribution (solid) and upper/lower frequency band (dashed). }
	\end{subfigure}
	\quad
	~ 
	\begin{subfigure}[b]{0.46\textwidth}
		\includegraphics[width=\textwidth]{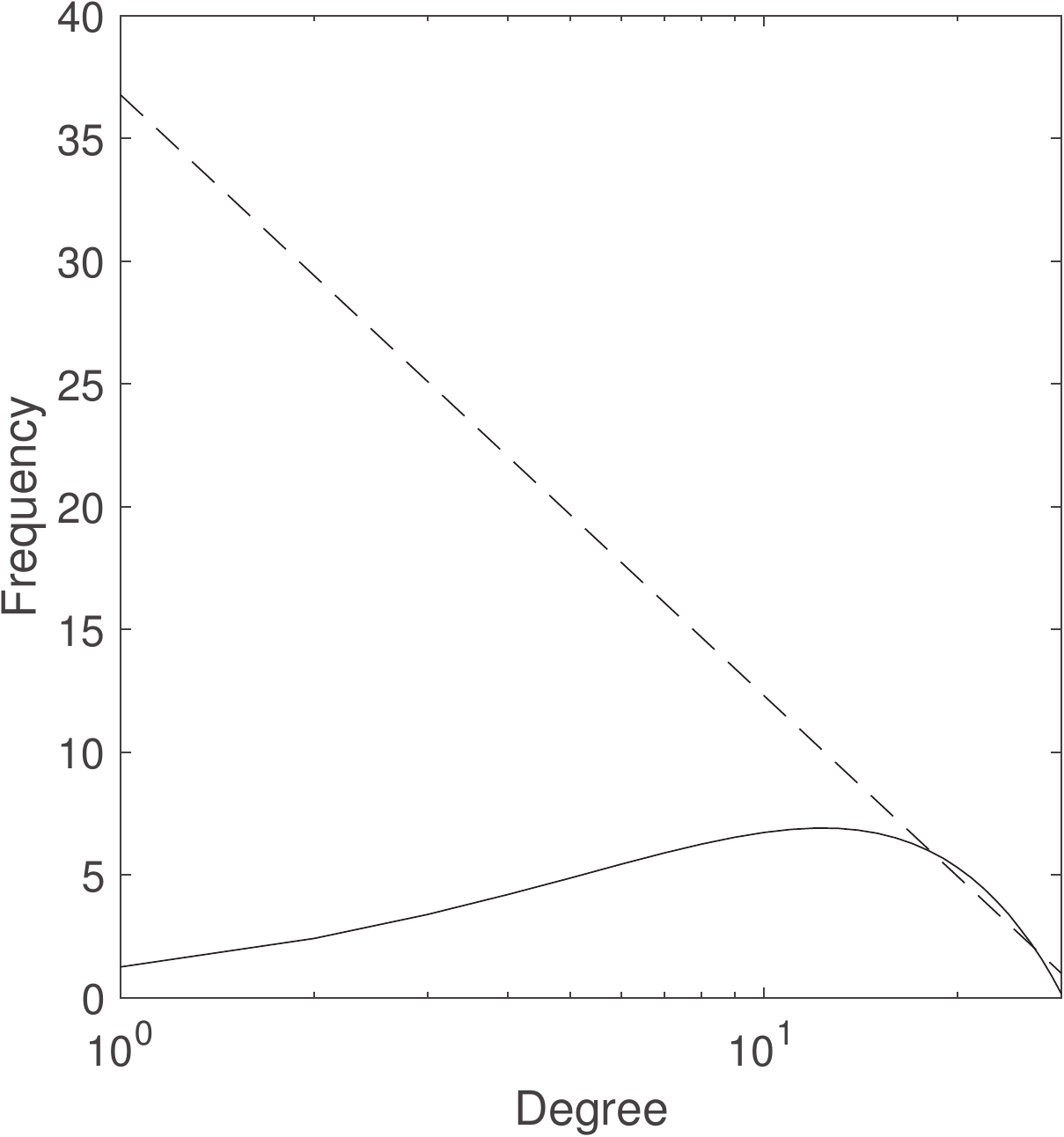}
		\caption{Powerlaw $\gamma_1=-10.2$(dashed) \& degree distribution (solid).}
	\end{subfigure}
	\caption{The degree distribution for $n = 10,000$, and  $c = 0.3$, with latent order $\omega_2$. }\label{gplot4}
\end{figure}
\begin{figure}
	\centering
	\begin{subfigure}[b]{0.475\textwidth}
		\includegraphics[width=\textwidth]{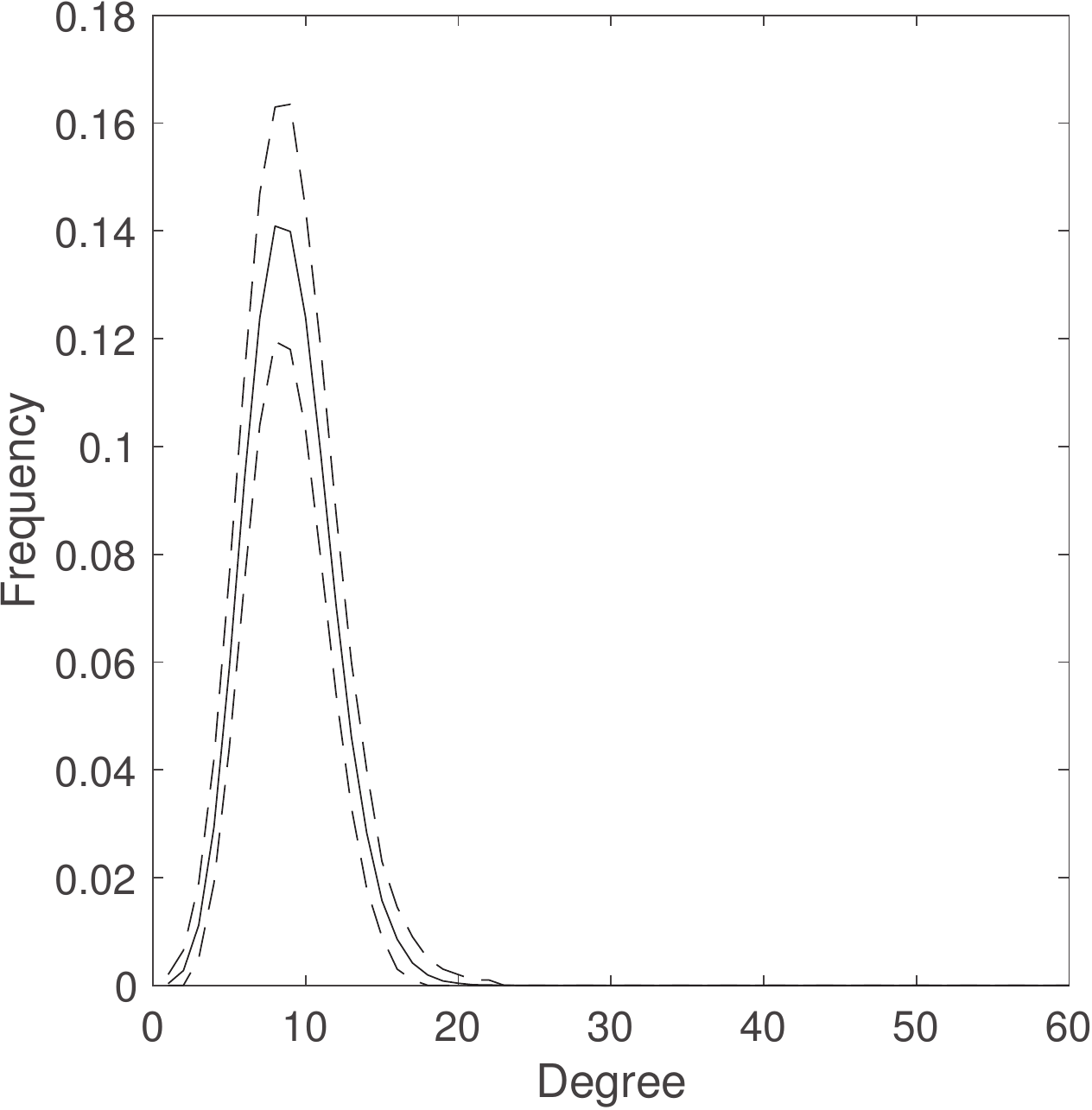}
		\caption{Mean distribution (solid) and upper/lower frequency band (dashed). }
	\end{subfigure}
	\quad
	~
	\begin{subfigure}[b]{0.46\textwidth}
		\includegraphics[width=\textwidth]{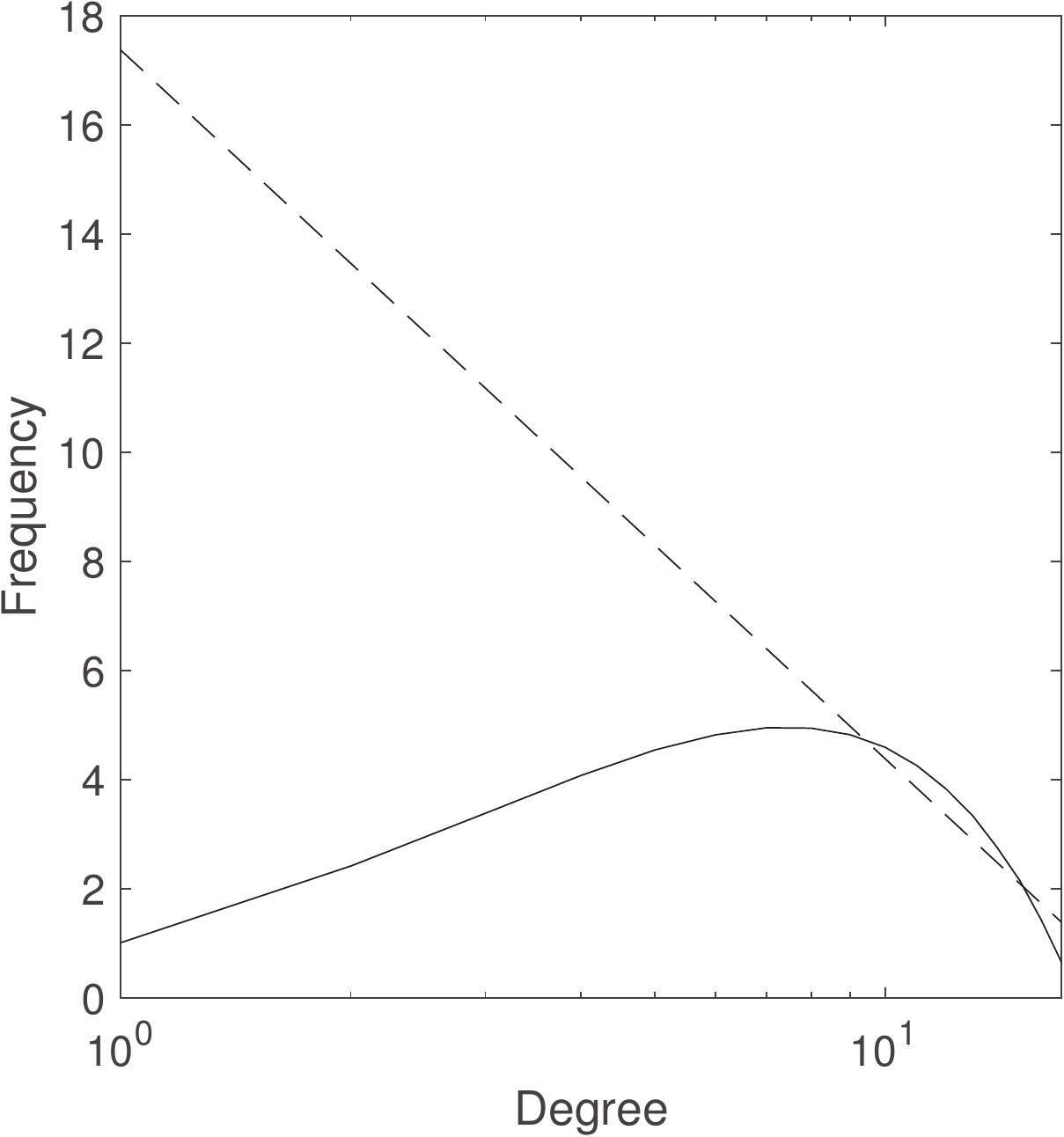}
		\caption{Powerlaw $\gamma_1=-5.6$(dashed) \& degree distribution (solid).}
	\end{subfigure}
	\caption{The degree distribution for $n = 1,000$, and  $c = 0.3$, for simple random graph.  }\label{gplot5}
\end{figure}
\begin{figure}
	\centering
	\begin{subfigure}[b]{0.475\textwidth}
		\includegraphics[width=\textwidth]{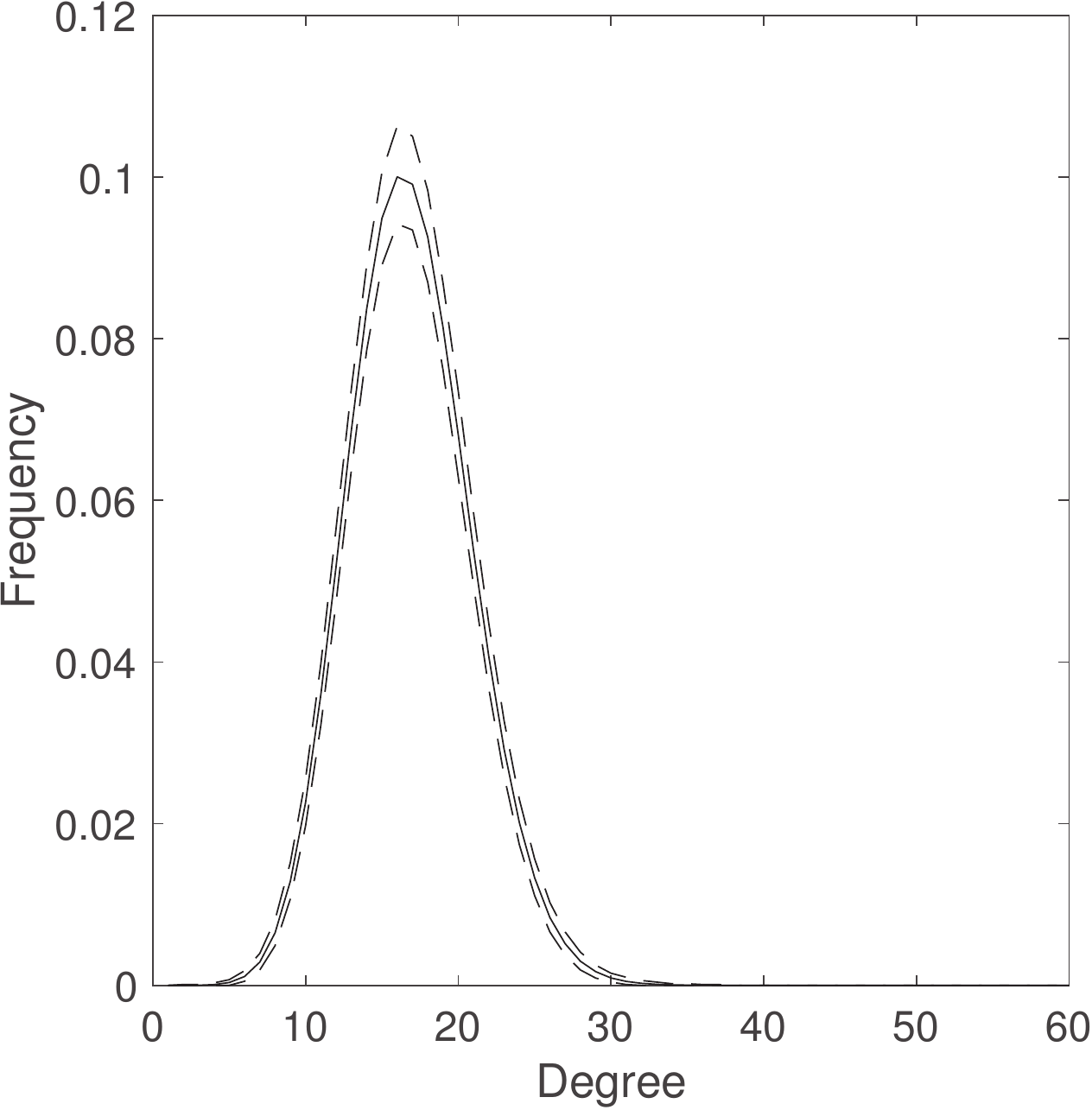}
		\caption{Mean distribution (solid) and upper/lower frequency band (dashed).}
	\end{subfigure}
	\quad
	~ 
	\begin{subfigure}[b]{0.46\textwidth}
		\includegraphics[width=\textwidth]{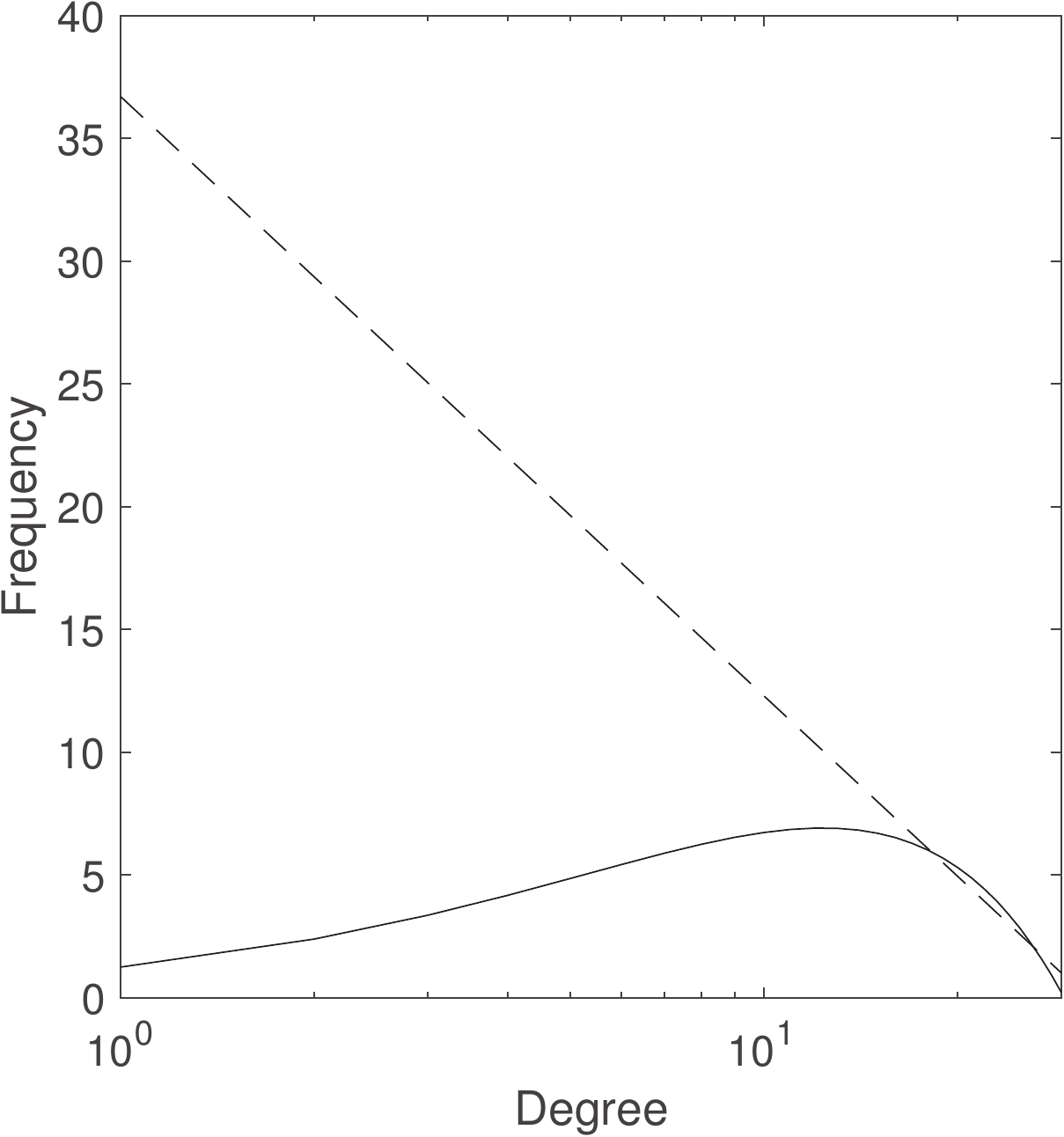}
		\caption{Powerlaw $\gamma_1=-10.6$(dashed) \& degree distribution (solid).}
	\end{subfigure}
	\caption{The degree distribution for $n = 10,000$, and  $c = 0.3$,  for simple random graph.  }\label{gplot6}
\end{figure}
\begin{figure}
	\centering
	\begin{subfigure}[b]{0.475\textwidth}
		\includegraphics[width=\textwidth]{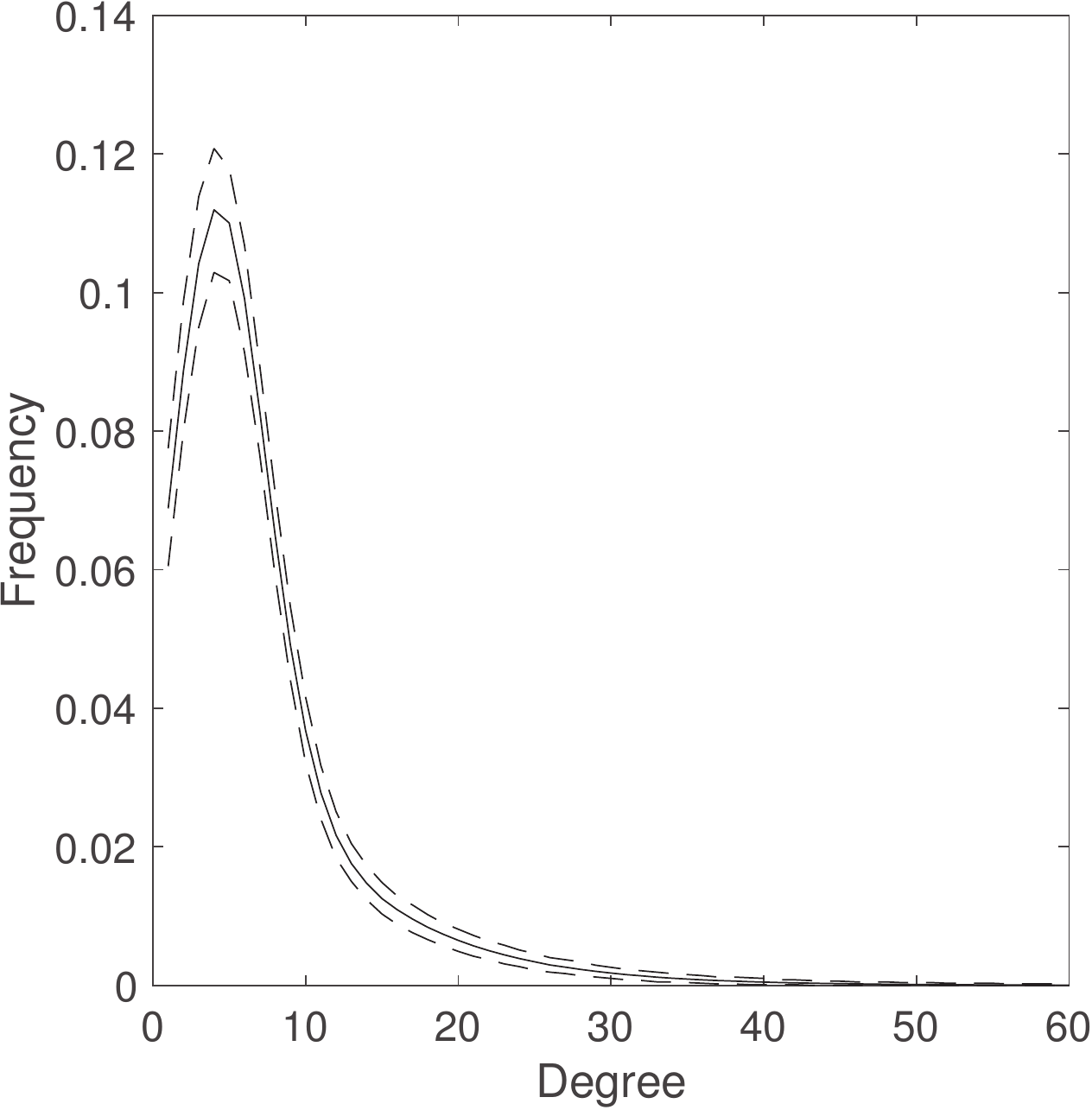}
		\caption{Mean distribution (solid) and upper/lower frequency band (dashed). }
	\end{subfigure}
	\quad
	~ 
	\begin{subfigure}[b]{0.46\textwidth}
		\includegraphics[width=\textwidth]{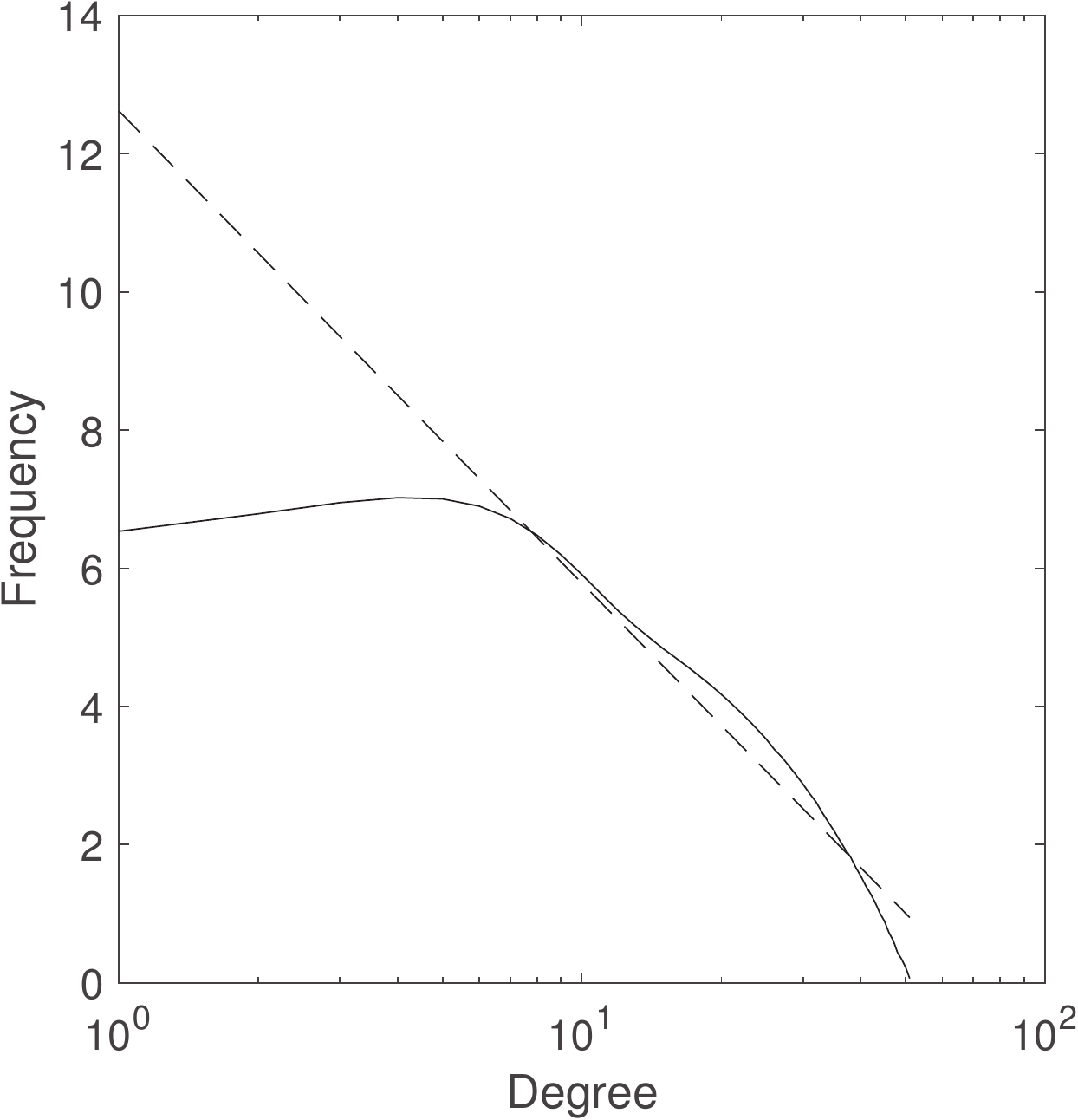}
		\caption{Powerlaw $\gamma_1=-2.97$(dashed) \& degree distribution (solid).}
	\end{subfigure}
	\caption{The degree distribution for $n = 10,000$, and  $c = 0.2$, for latent order $\omega_1$.  }\label{gplot7}
\end{figure}
\begin{figure}
	\centering
	\begin{subfigure}[b]{0.475\textwidth}
		\includegraphics[width=\textwidth]{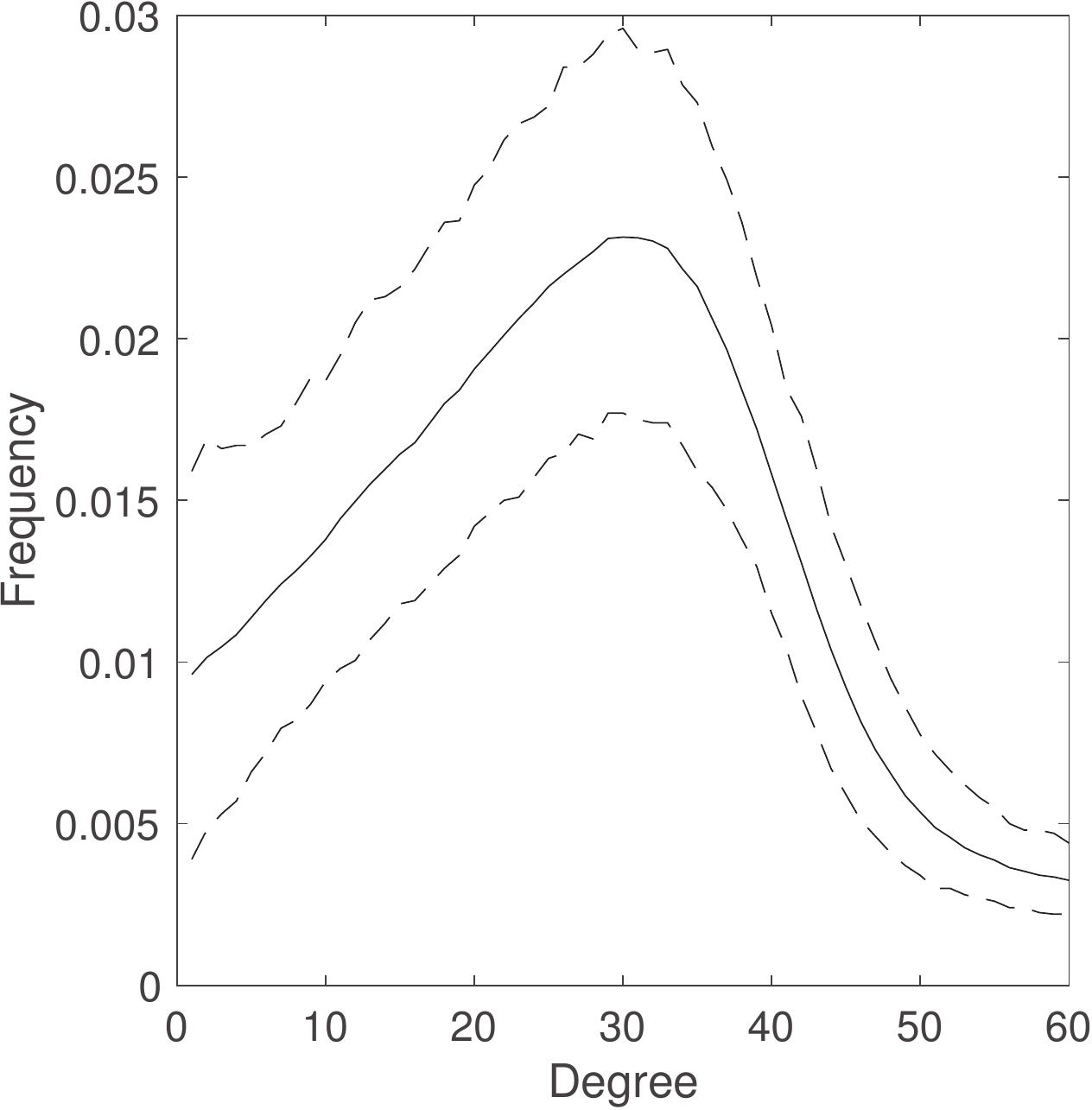}
		\caption{Mean distribution (solid) and upper/lower frequency band (dashed). }
	\end{subfigure}
	\quad
	~
	\begin{subfigure}[b]{0.46\textwidth}
		\includegraphics[width=\textwidth]{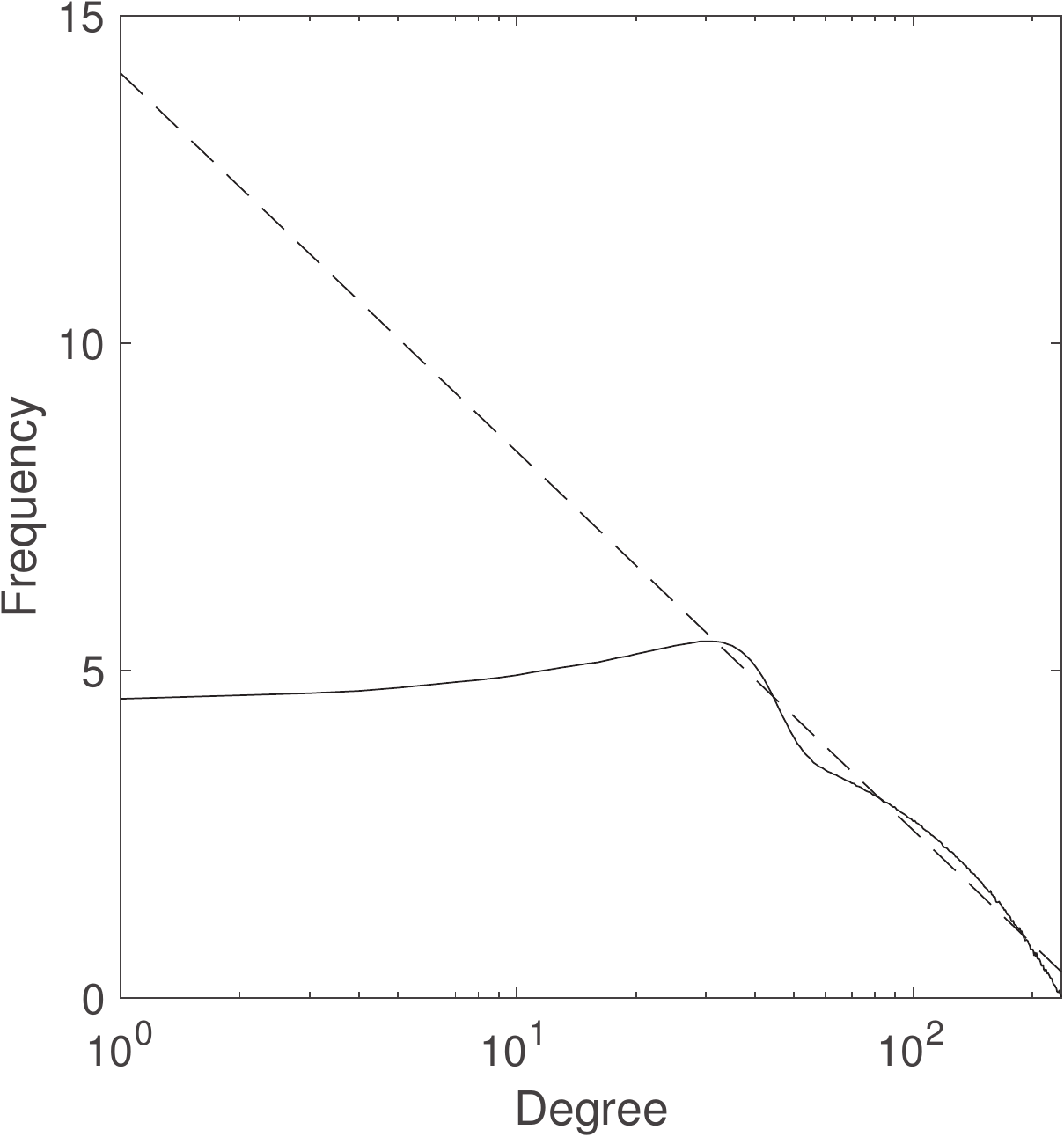}
		\caption{Powerlaw $\gamma_1=-2.52$(dashed) \& degree distribution (solid).}
	\end{subfigure}
	\caption{The degree distribution for $n = 10,000$, and  $c = 0.4$, for latent order $\omega_1$.  }\label{gplot8}
\end{figure}

\newpage
In this section we  display  Figure  \ref{network} which  shows typical graphs generated from MECLTGs $CG(V,\omega_1,p_0,p_1)$, $CG(V,\omega_2,p_0,p_1)$ and Erd\H{o}s-R\'{e}nyi $G(V,p)$, {respectively} with $p_0=0.01$, $p_1=1-\frac{1}{n^{1/3}}$, $p=\frac{p_0}{1-p_1+p_0}$  with $|V|=100$. Then we show the degree distribution of various MECLTGs in Section \ref{simuempricaldegree} of the main article with different sizes, different values of $c$ and two choices of latent order $\omega$, i.e., $\omega_1$ and $\omega_2$. Results are shown in Figures \ref{gplot2}--\ref{gplot8}. The related analysis are concluded in Section \ref{simpleedge}  of the main article.

\section{Phase Transition for Connectivity and Giant Component}\label{gai}
In this section we discuss the phase transition of connectivity and of the giant component for the first order  homogeneous MECLTG defined in Section \ref{simpleedge} of the main article. It is well known that the threshold probability for connectivity in simple random graph (SRG) $G(n,p)$ is $p=\frac{\log n}{n}$ (e.g.\cite{durrett2007random}).  The threshold probability for the emergence of a giant component is $p=\frac{1}{n}$ {(e.g. \cite{durrett2007random})}. However, the traditional method of finding the threshold probability relies heavily on the assumptions of edge variables independence that does not hold for  MECLTG. For example, to calculate the threshold probabilities for connectivity, many traditional methods need to evaluate $Cov (I_i, I_j)$, where $I_i=1$ if node $i$ is isolated and $0$ otherwise. To calculate the threshold for the emergence of giant component in the SRG, techniques based on branching process, random walk, or depth first search algorithm are proposed. However, due to the complicated dependence structure and the unobserved latent order $\omega_n(\{\cdot,\cdot\})$, these techniques are not directly applicable to the first order  homogeneous MECLTG. The threshold probabilities of the connectivity and of the emergence of a giant component are discussed for example
\cite{durrett2007random}, \cite{KrSu2013}.
In the following theorem, we show by construction that under some circumstances, the first order  homogeneous MECLTG possess the threshold properties similar to SRG.
\begin{theorem}\label{connectiviy} Let $G_n=CG(V_n,\omega,p_{0,n},p_{1,n})$ be a series of first order  homogeneous MECLTG, where $V_n$ is the set of vertices such that $|V_n|=n$, and $p_{w,n}=\frac{\lambda_w\log n}{n}$, with $w\in\{0,1\}$ and $\lambda_0,\lambda_1>0$ are positive constants. 
	\begin{description}
		\item (a)	Denote by $\mathcal A_n$ the event that $G_n$ is connected. 
		Then we have i) 
		if $\max\{\lambda_0,\lambda_1\}<1$, then $\lim_{n\rightarrow \infty}\p(\mathcal A_n)= 0$; ii) if $\min\{\lambda_0,\lambda_1\}>1$ then  $\lim_{n\rightarrow \infty}\p(\mathcal A_n)= 1$.
		\item (b)  Let $|S_i|$ be the size of the component that contains $i$. As $n\rightarrow \infty$, we have: i) if $\max\{\lambda_0,\lambda_1\}<1$, then for all (sufficiently large) $a>0$, $\lim_{n\rightarrow \infty}\p(\max_{1\leq i\leq n}|S_i|\geq a\log n)= 0$;  ii) if $\min\{\lambda_0,\lambda_1\}>1$, then 
		there exists $c>0$, such that \begin{align}\p(\text{There exists a component with size} \geq c\sqrt n)\rightarrow 1 \end{align}
		as $n\rightarrow \infty$.
	\end{description}
\end{theorem}
To avoid  tedious computation due to the latent order $\omega_n(\{\cdot,\cdot\})$, and comply with the complex dependence structure, we prove the theorem via studying an algorithm that generates the SRG and the first order MECLTG simultaneously. The algorithm reveals the connection between the two constructions.

{\it Proof}.
Let $p_a$ be a pre-specified number, and $\omega$ be the latent order. Consider the following algorithm, which construct MECLTG by a series of random variables $\tilde A_i$ and $\tilde B_i$, $1\leq i\leq N.$
\begin{description}
	\item (a) Generate $U\sim $Uniform$ (0,1)$. Set $\tilde B_1=1$ if $U\leq p$, where $p$ is a function of $p_0=p_{0,n}$, $p_1=p_{1,n}$ The function is defined in equation \eqref{stationp} of the main article. Set $\tilde A_1=1$ if $U\leq p_a$.
	\item (b) At steps $i\geq 2$, generate independently a new $U\sim $Uniform$ (0,1)$. If $\tilde B_{i-1}=1$, set $\tilde B_i=1$ if $U\leq p_1$. If $\tilde B_{i-1}=0$, set $\tilde B_i=1$ if $U\leq p_0$. Set $\tilde A_i=1$ if $U\leq p_a$.
\end{description}
Thus, $\{B_{\omega^{-1}(i)},1\leq i\leq N\}=\{\tilde B_i,{1\leq i\leq N}\}$ forms the ordered edge variables of a  first order  homogeneous MECLTG  $CG(V,\omega,p_0,p_1)$, while $\{A_{\omega^{-1}(i)},1\leq i\leq N\}=\{\tilde A_{i}, 1\leq i\leq N\}$ forms the edge variables of a first order  homogeneous MECLTG $CG(V,\omega,p_a,p_a)$. By our definition, $CG(V,\omega,p_a,p_a)$ is the simple random graph $G(|V|,p_a)$. Write $\bar p_1=\max\{p_1, p_0\}$, $\bar p_0=\min \{p_1, p_0\}$. Note that $p\in [\bar p_0, \bar p_1]$.

From the construction, when $\bar p_1\leq p_a$, then in each step $i$, $1\leq i\leq N$, $A_{\omega^{-1}(i)}=0$ implies $B_{\omega^{-1}(i)}=0$. This is because if $A_{\omega^{-1}(i)}=0$ (this means $\tilde A_{i}=0$),
\ then $U>p_a\geq \bar p_1$ in step (i).
This will make $\tilde B_i=0$, and consequently $B_{\omega_n^{-1}(i)}=0$. The above fact shows that when $\bar p_1\leq p_a$, if node $i$ is isolated in the constructed $CG(V,\omega,p_a,p_a)$, the corresponding node is also isolated in the constructed $CG(V,\omega,p_0,p_1)$. As a result, we have that $|S_i^A|\geq |S_i^B|$ where $|S_i^A|$ is the size of the component contains node $i$ in $CG(V,\omega,p_a,p_a)$, and $|S_i^B|$ is the size of the component contains node $i$ in $CG(V,\omega,p_0,p_1)$.
Since $CG(V,\omega,p_a,p_a)$ is the Erd\"{o}s-R\'{e}nyi Graph SRG $G(|V|,p_a)$,  when $\bar p_1\leq p_a$ the theorem follows from the well-known results of the threshold probabilities of the connectivity and the emergence of a giant component for Erd\"{o}s-R\'{e}nyi Graph (e.g.\cite{durrett2007random}).  When $\bar p_0\geq p_a$, the theorem follows from  a similar argument, and the proof is completed. 
\hfill $\Box$

\section{Proof of auxiliary results for Theorem \ref{cconsistent} and Theorem  \ref{Marginalgraphon} in the main article}\ \\\label{Sec:2}

\noindent{\bf Proof of Proposition  \ref{dependenceprop}.}
By assumption $\min_{a,b\in \mathbf \chi\times \mathbf \chi} \p(U_{i}=a|U_{i-l}=b)\geq \alpha'>0$, 
By definition, $\alpha'\leq \frac{1}{2}$.
The key to proof the proposition is to show that for $k\geq l$,
\begin{align}\label{2018-New-S2}
\max_b \p(U_i=u_i|U_{i-k}=b)-
&\min_b \p(U_i=u_i|U_{i-k}=b)\notag \\&\leq (1-2\alpha')[\max_b \p(U_i=u_i|U_{i-k+l}=b)-\min_b \p(U_i=u_i|U_{i-k+l}=b)].
\end{align} If equation \eqref{2018-New-S2} holds
then for $k=lq+s$,  $ 0\leq s\leq l-1$, we have by iteratively applying \eqref{2018-New-S2}, the following inequality
\begin{align}\label{2018-New-S2-new}
\max_b \p(U_i=u_i|U_{i-k}=b)-
&\min_b \p(U_i=u_i|U_{i-k}=b)\notag \\&\leq (1-2\alpha')^{q}[\max_b \p(U_i=u_i|U_{i-s}=b)-\min_b \p(U_i=u_i|U_{i-s}=b)]\notag\\
&\leq M(1-2\alpha'
)^{k/l}\tilde p
\end{align} holds,
where $M=(1-2\alpha')^{-\frac{l-1}{l}}\leq (1-2\alpha')^{-1}$. 
Then the proposition follows from equation \eqref{2018-New-S2-new}.
It remains to show \eqref{2018-New-S2}. Let $b_0=\argmin_b\p(U_i=u_i|U_{i-k+l}=b)$. Notice that\begin{align}\label{New.S3}
&\p(U_i=u_i|U_{i-k}=u_{i-k})\notag \\&=\sum_{b\in \mathbb R^l}\p(U_i=u_i|U_{i-k+l}=b)\p(U_{i-k+l}=b|U_{i-k}=u_{i-k})\notag
\\&\leq\p(U_i=u_i|U_{i-k+l}=b_0)\p(U_{i-k+l}=b_0|U_{i-k}=u_{i-k})\\&\ \ \ \ \ \ \ \ \ \  +\sum_{b\neq b_0}(\max_{b}\p(U_i=u_i|U_{i-k+l}=b))\p(U_{i-k+l}=b|U_{i-k}=u_{i-k})\notag\\
&=\max_{b}\p(U_i=u_i|U_{i-k+l}=b)\notag\\&-\p(U_{i-k+l}=b_0|U_{i-k}=u_{i-k})\times\left(\max_{b}\p(U_i=u_i|U_{i-k+l}=b)-\min_{b}\p(U_i=u_i|U_{i-k+l}=b)\right)
\notag \\&\leq \max_{b}\p(U_i=u_i|U_{i-k+l}=b)-\alpha'\left(\max_{b}\p(U_i=u_i|U_{i-k+l}=b)-\min_{b}\p(U_i=u_i|U_{i-k+l}=b)\right)
\end{align}
By taking maximum on both side of \eqref{New.S3}, we have\begin{align}\label{2018-New-S4}
\max_b \p(U_i=u_i|U_{i-k}=b)&\leq \max_{b}\p(U_i=u_i|U_{i-k+l}=b)\notag\\&-\alpha'\left(\max_{b}\p(U_i=u_i|U_{i-k+l}=b)-\min_{b}\p(U_i=u_i|U_{i-k+l}=b)\right)
\end{align}
Similarly
\begin{align}\label{2018-New-S5}
\min_b \p(U_i=u_i|U_{i-k}=b)&\geq \min_{b}\p(U_i=u_i|U_{i-k+l}=b)\notag\\&+\alpha'\left(\max_{b}\p(U_i=u_i|U_{i-k+l}=b)-\min_{b}\p(U_i=u_i|U_{i-k+l}=b)\right)
\end{align}
The \eqref{2018-New-S2} follows from \eqref{2018-New-S4} and \eqref{2018-New-S5}.\hfill $\Box$\\\ \\
Recall $\chi$ of Proposition \ref{dependenceprop}
of the main article  in the following arguments.
\begin{lemma}\label{Newprop3}
	Assume that the conditions of Theorem \ref{cconsistent} holds. Recall $Z_{n,k}$ defined in Section \ref{Sec:Minimax}. Then for all $z\in\mathcal Z_{n,k}$, there exist constant $C_1$ (which does not depend on $n$) such that
	for $a,b\in [k]$,\begin{align}
	\E\left(\exp\left(\left(\frac{C_1\left(1-\chi\right)}{\sqrt {n_an_b}}\sum_{i\in \bar{z}^{-1}(a),j\in \bar{z}^{-1}(b)}(A_{i,j}-\theta_{i,j})\right)^2\right)\right)\leq 2.
	\end{align}
\end{lemma}
{\it Proof.} Consider the case of $a\neq b$. The case of $a=b$ follows {\it mutatis mutandis}.  Recall the map $\omega(\cdot,\cdot)$ and define $B_i,i\in [N]$ such that $B_{\omega(\{i,j\})}=A_{i,j}$. Define a series of integers $u_j$, $1\leq j\leq n_an_b$ as
\begin{align}\label{defbaromega}\{u_j,1\leq j\leq n_an_b\}=\{\omega(\{i,j\}),i\in \bar{z}^{-1}(a),j\in \bar{z}^{-1}(b)\}\end{align}
such that $u_s<u_l$ for $s<l$. 
\ Consider the filtration $\FF_i=(\{B_{u_j},j\leq i\}, z_1,...,z_n)$, where $z:=\{z_i,1\leq i\leq n\}$ are latent $U(0,1)$ variables. Let 
$B_{u_s}=B_s$ for $s\leq 0$.
Define the projection operator for $j\in \mathbb Z$,
\begin{align}
\pp_{j}(\cdot)=\E(\cdot|\FF_{j},z)-\E(\cdot|\FF_{j-1},z).
\end{align}
By our construction of $B_i$, $i\leq 0$ in Definition \ref{defmarkovgraph} of the main article, we get
\begin{align}
\pp_{j}(B_k)=0\  \text{for\ } j\leq 0,\ \ \ k>0.
\end{align}
It follows from the fact that $\theta_{i,j}=\E(A_{i,j}|z)$,
\begin{align}\label{new.eq67}
\sum_{i\in \bar{z}^{-1}(a),j\in \bar{z}^{-1}(b)}\left\{A_{i,j}-\E(A_{i,j}|z)\right\}=\sum_{i=1}^{n_an_b}\sum_{s=0}^{\infty}\pp_{i-s}B_{u_i}=\sum_{s=0}^{\infty}\sum_{i=1}^{n_an_b}\pp_{i-s}B_{u_i}
.\end{align}
Note that $\pp_{i-s}B_{u_i}$ forms a martingale difference w.r.t $\{\FF_{i-1},z_1,...,z_n\}$.
By Burkholder inequality, for $v>0$,
we have that
\begin{align}
\|\sum_{i=1}^{n_an_b}\pp_{i-s}B_{u_i}\|^2_{\mathcal L_v }\leq C v\sum_{i=1}^{n_an_b}\|\pp_{i-s}B_{u_i}\|^2_{\mathcal L_v },
\end{align}
where $C$ is a constant independent of $v$, $n_a$, $n_b$ and $n$, and $\|X\|_{\mathcal L_v}:=\{\E(|X|^v)\}^{1/v}$. By Corollary \ref{Coroldependence} in the main article and the fact that $u_i-u_{i-s}\geq s$, we have
\begin{align}
|\pp_{i-s}B_{u_i}|\leq C'\chi^s.
\end{align}
Here the above equation does not depend on $i$ due to Corollary \ref{Coroldependence} which only involves the distant of indices.
Thus we obtain \begin{align}
\left\|\sum_{i=1}^{n_an_b}\pp_{i-s}B_{u_i}\right\|^2_{\mathcal L_v }\leq C'' vn_an_b\chi^{2s}.
\end{align}
Combining with \eqref{new.eq67} and by triangle inequality, we have that
\begin{align}\label{2018-newS.15}
\left\|\sum_{i\in \bar{z}^{-1}(a),j\in \bar{z}^{-1}(b)}(A_{i,j}-\theta_{i,j})\right\|_{\mathcal L_v }\leq  \frac{C_1\sqrt{n_an_bv}}{1-\chi},
\end{align}
where $C'$, $C''$ and $C_1$ are constants independent of $n$ and $v$. By using Taylor expansion, equation \eqref{2018-newS.15} shows that there exists a small $\eta\in(0,1/2)$ such that
\begin{align}\label{New.74-2016}
\E\left(\exp\left(\eta\left(\frac{\left(1-\chi\right)}{C_1\sqrt {n_an_b}}\sum_{i\in \bar{z}^{-1}(a),j\in \bar{z}^{-1}(b)}(A_{i,j}-\theta_{i,j})\right)^2\right)\right) \leq 1+\sum_{u=1}^\infty\frac{u^u\eta^u}{u!}\leq 2.
\end{align}
The last inequality use the fact that $u!>(\frac{u}{e})^u$ and we take $\eta$ small such that  $\eta e\leq 1/2$.
\hfill $\Box$
\begin{corol}\label{corol2}Assume that the conditions of Theorem \ref{cconsistent} holds.
	Consider real numbers $\{\gamma_{i,j}\}_{1\leq i\leq n,1\leq j\leq n, i\neq j}$ satisfies $\sum_{i,j}\gamma_{i,j}^2=1$.
	Then we have that there exists a constant $C_2$ such that\begin{align}
	\E\left(\exp\left(\left(C_2(1-\chi)\sum_{i,j}\gamma_{i,j}(A_{i,j}-\theta_{i,j})\right)^2\right)\right)\leq 2.
	\end{align}
\end{corol}
{\it Proof.} Recall the proof of Lemma \ref{Newprop3} and the filtration $\FF_i=(\{B_{j},j\leq i\},z)$ defined there, associate with the projection operator 
$\pp_{j}(\cdot)=\E(\cdot|\FF_{j},z)-\E(\cdot|\FF_{j-1},z)$ in Lemma \ref{Newprop3}. Define a series $\gamma'$ such that $\gamma'_{ \omega\{i,j\}}=\gamma_{i,j}$ for $i<j$ such that $\sum_{j=1}^N\gamma'^2_j=1/2$. Thus we have
\begin{align}
\sum_{i>j}\gamma_{i,j}\left(A_{i,j}-\theta_{i,j}\right)=\sum_{s=0}^{\infty}\sum_{i=1}^{N}\gamma'_i\pp_{i-s}B_{i}.
\end{align}
By Burkholder's inequality and  the triangle inequality, a similar argument to that of Lemma \ref{Newprop3} yields that for $v\geq 0$, \begin{align}
\left\|\sum_{i>j}\gamma_{i,j}\left(A_{i,j}-\theta_{i,j}\right)\right\|_{\mathcal L_v }\leq C_2\sqrt v/(1-\chi),
\end{align}
where $C_2$ is a constant independent of $n$ and $v$.
The corollary follows from the same argument of \eqref{New.74-2016}.\hfill $\Box$
\begin{lemma}\label{New.lemma_7} Assume that the conditions of Theorem \ref{cconsistent} hold.
	Then for any $C'>0$, there exists  $C>0$
	such that \begin{align}
	\left|\left\langle \frac{\tilde \theta-\theta}{\|\tilde \theta-\theta\|},A-\theta\right\rangle\right|\leq C\sqrt{n\log k}(1-\chi)^{-1}
	\end{align}
	with probability $1-\exp(-C'n\log k)$.
\end{lemma}
{\it Proof.} By lemma 4.2 of Gao et al. (2015), we have that
\begin{align}
\left|\left\langle \frac{\tilde \theta-\theta}{\|\tilde \theta-\theta\|},A-\theta\right\rangle\right|\leq \max_{z\in \mathcal Z_{n,k}}\left|\sum_{i,j}\gamma_{i,j}(z)(A_{i,j}-\theta_{i,j})\right|
\end{align}
where $\gamma_{i,j}({z})\propto \sum_{a,b\in [k]}\bar \theta_{a,b}({z})\mathbf 1((i,j)\in {z}^{-1}(a)\times {z}^{-1}(b))-\theta_{i,j}$, satisfying $\sum_{i,j}\gamma_{i,j}^2({z})=1$.
Then Markov's inequality, union bound and Corollary \ref{corol2} lead to 
\begin{align}\label{S22}
\p\left((1-\chi)\max_{z\in \mathcal Z_{n,k}}|\sum_{i,j}\gamma_{i,j}(z)(A_{i,j}-\theta_{i,j})|>t\right)\leq C_1\exp(-C_0t^2+n\log k).
\end{align}
Then the lemma follows from letting $t=M\sqrt{n\log k}$ in \eqref{S22} for some large constant $M$.
\hfill $\Box$
\begin{lemma}\label{New.lemma_8}Assume that the conditions of Theorem \ref{cconsistent} hold.
	Then for any constant $C'>0$, there exists a constant $C>0$ only depending on $C'$, such that
	\begin{align}
	\left|\left\langle \frac{\hat \theta-\tilde \theta}{\|\hat \theta-\tilde \theta\|},A-\theta\right\rangle\right|\leq C\sqrt{k^2+n\log k}(1-\chi)^{-1}
	\end{align}
	with probability at least $1-\exp(-C'n\log k).$
\end{lemma}
{\it Proof.} Define an $n-$dimensional ball $\mathcal B\subset\left\{a\in \mathbb R^{n\times n}:\sum_{ij}a_{ij}^2\leq 1\right\}$ with the following property \begin{align}\label{new.property1}
\text{if}\  a,b\in \mathcal B,\  \text{then}\  \frac{a-b}{\|a-b\|}\in\mathcal B.
\end{align}
Let $\mathcal B'$ be a $1/2$-net of $\mathcal B$ such that $|\mathcal B'|\leq \mathcal N\left(1/2,\mathcal B,\|\cdot\|\right)$ where  $\mathcal N\left(1/2,\mathcal B,\|\cdot\|\right)$ is the covering number of $1/2$ net of set $\mathcal B$.
By Lemma A.1. of Gao et al. (2015), we have that:
\begin{align}
\sup_{a\in\mathcal B}\left|\sum_{i,j}a_{i,j}(A_{i,j}-\theta_{i,j})\right|\leq 2\max_{b\in\mathcal B'}\left|\sum_{i,j}b_{i,j}(A_{i,j}-\theta_{i,j})\right|.
\end{align}
As a consequence, by Corollary \ref{corol2} and the union bound, we have
\begin{align}\label{2018-new.2.26}
&\p\left((1-\chi)\sup_{a\in\mathcal B}\left|\sum_{i,j}a_{i,j}(A_{i,j}-\theta_{i,j})\right|\geq t\right)\notag\\
&\leq\p\left( 2(1-\chi)\max_{b\in\mathcal B'}\left|\sum_{i,j}b_{i,j}(A_{i,j}-\theta_{i,j})\right|\geq t\right)\notag
\\&\leq |\mathcal B'| C_1\exp(-C_0 t^{2}),
\end{align}
for some constants $C_0\in \mathbb R^+$ and $C_1\in \mathbb R^+$.
For $\bar{z}\in \mathcal Z_{n,k}$, define
\begin{align}
\mathcal B_{\bar z}=\large\{\{c_{i,j}\}:c_{i,j}=Q_{ab} \ \text{if}\  (i,j)\in \bar{z}^{-1}(a)\times \bar{z}^{-1}(b)\ \text{for some}\ Q_{ab}, \notag
\ \\ \text{and}\ \sum_{i,j}c_{i,j}^2\leq 1, c_{i,j}=c_{j,i}, c_{i,i}=0, 1\leq i\leq n\}.
\end{align}
Notice that $\mathcal B_{\bar z}$ satisfy \eqref{new.property1}, and $\mathcal N\left(1/2,\mathcal B_{\bar z},\|\cdot\|\right)\leq \exp(C_2k^2)$ for some constant $C_2\in \mathbb R^+$.
By the proof of Corollary \ref{corol2} and the union bound, we have that
\begin{align}\label{S.28}
\p\left((1-\chi)\left|\left\langle \frac{\hat \theta-\tilde \theta}{\|\hat \theta-\tilde \theta\|},A-\theta\right\rangle\right|>t\right)\\\notag
\leq
\p\left((1-\chi)\max_{z\in \mathcal Z_{n,k}}\sup_{c\in \mathcal B_z}\left|\sum_{ij}c_{i,j}\left(A_{i,j}-\theta_{i,j}\right)\right|>t\right)\notag\\
\leq C_0\exp\left(-C_0t^2+C_1n\log k+C_2k^2\right)\notag
\end{align}
for some positive constants $C_0$, $C_1$, $C_2$. The last inequality requires \eqref{2018-new.2.26}. Then the lemma follows from taking $t=M\sqrt{n\log k+k^2}$ for sufficiently large $M$ in \eqref{S.28}.
\hfill$\Box$
\section{Proof of auxiliary results for Theorem \ref{Spectralthm3}}\label{Sec:3}
For latent variable $z$ and random variable $X$, let $\|X|z\|_{\mathcal L_2}=\sqrt{\E(X^2|z)}$. We need the following proposition to show Proposition \ref{newprop6} of the main article.
\begin{proposition}\label{prop5}
	Consider the latent undirected Markov model with memory parameter $h<\infty$. Let $\FF_i$ be the filtration generated by $\{B_s\}_{s\leq i}$, where  $\{B_s\}_{s\leq i}$ is defined in Definition \ref{defmarkovgraph} in the main article.
	Let $\delta(k)$: $\delta(k)=(1-2\alpha')^{-1/2}\chi^{k/2}$ for $k\geq 0$, $\chi$ and $\alpha'$ defined in Proposition \ref{dependenceprop} in the main article, and $0$ if $k\leq -1$.
	Then for any integers
	$i_1\leq..\leq i_l$ for some fixed constant $l$, we have that
	\begin{align}
	\|\pp_a(\Pi_{s=1}^lB_{i_s})|z\|_{\mathcal L_2}\leq 2\sum_{s=1}^l\delta(i_s-a)\mathbf 1(i_{s-1}<a),
	\end{align}
	where $\pp_j(\cdot)=\E(\cdot|\FF_j,z)-\E(\cdot|\FF_{j-1},z)$, $\FF_i=(B_{-\infty},...,B_i)$ and $i_0=-\infty$.
\end{proposition}
{\it Proof.} Without loss of generality, let $i_s\leq i_{s+1}$ for $s\leq s+1$. Then by direct calculations,
we obtain that
\begin{align}
\|\pp_a(\Pi_{s=1}^lB_{i_s})|z\|^2_{\mathcal L_2}&=\E\left(
\E\left(\pp^2_a(\Pi_{s=1}^lB_{i_s})|\FF_{a-1},z\right)|z\right)\notag
\\&=\E(\E^2(\Pi_{s=1}^lB_{i_s}|\FF_a,z)|z)-\E(\E^2(\Pi_{s=1}^lB_{i_s}|\FF_{a-1},z)|z):=\E(AB|z),
\end{align}
where \begin{align}
A=\E(\Pi_{s=1}^lB_{i_s}|\FF_a,z)+\E(\Pi_{s=1}^lB_{i_s}|\FF_{a-1},z),\\ B=\E(\Pi_{s=1}^lB_{i_s}|\FF_a,z)-\E(\Pi_{s=1}^lB_{i_s}|\FF_{a-1},z).
\end{align}
As a result we have that 
\begin{align}
\|\pp_a(\Pi_{s=1}^lB_{i_s})|z\|_{\mathcal L_2}\leq (\|A|z\|_{\mathcal L_v}\|B|z\|_{\mathcal L_v})^{1/2}.
\end{align}
By the boundedness of Bernoulli random variable, we have that $|A|\leq 2$.
For $a\geq i_l+1$, it is easy to see that $\|\pp_a(\Pi_{s=1}^lB_{i_s})|z\|=0$.
Recall the definition of $U_i=(B_i,...,B_{i-h+1})'$ below the Definition \ref{defmarkovgraph} in the main article.
Let $u_i=(b_i,...,b_{i-h+1})'\in\{0,1\}^{h}$.
Note that for $a\leq i_1-1$ and $v=a, a-1$,
\begin{align}
\notag&\E(\Pi_{s=1}^lB_{i_s}|\FF_v,z)=\\
&\sum_{ \substack{u_i=(b_i,...,b_{i-h+1})',,\\b_j=\{0,1\},i_1\leq j\leq i_l+h,\\ j\neq i_s,1\leq s\leq l}}
\Pi_{i=i_1+h}^{i_l+h}\p(U_i=u_i|U_{i-1}=u_{i-1},z)\p(U_{i_1+h-1}=u_{i_1+h-1}|\FF_v,z)\big|_{b_{i_s}=1,s=1,...,l}.
\end{align}
Notice that
\begin{align}
\sum_{ \substack{u_i=(b_i,...,b_{i-h+1})',,\\b_j=\{0,1\},i_1\leq j\leq i_l+h,\\ j\neq i_s,1\leq s\leq l}}
\Pi_{i=i_1+h}^{i_l+h}\p(U_i=u_i|U_{i-1}=u_{i-1},z)\big|_{b_{i_s}=1,s=1,...,l}\\\notag
=\E(\Pi_{s=j}^lB_{i_s}|\FF_{i_1+h-1},z), \text{\ $j$ satisfies $i_{j-1}\leq i_1+h-1\leq i_j$}.
\end{align}
Using this fact,  Proposition \ref{dependenceprop} in the main article and the upper bound of $M$ in equation \eqref{2018-New-S2-new} we have that $B\leq 2(1-2\alpha')^{-1}\chi^{i_1-a}$. Consequently, we have that 
\begin{align}\|\pp_a(\Pi_{s=1}^lB_{i_s})|z\|_{\mathcal L_2}\leq 2\delta(i_1-a).\label{new.106}
\end{align}
Observe that \eqref{new.106} still holds for $a= i_1$. For $i_u+1\leq a\leq i_{u+1}$, $1\leq u\leq l-1$, we have that for $v=a,a-1$,
\begin{align}
\E(\Pi_{s=1}^lB_{i_s}|\FF_v,z)=(\Pi_{w=1}^uB_{i_w})\E(\Pi_{s=u+1}^lB_{i_s}|\FF_v,z).
\end{align}
Therefore we have that\begin{align}
\pp_a(\Pi_{s=1}^lB_{i_s})=(\Pi_{w=1}^uB_{i_w})\pp_a(\Pi_{s=u+1}^lB_{i_s}|\FF_v,z).
\end{align}
It follows from the boundedness of $B_i's$ and a similar argument to the  $a<i_1$ case  that
\begin{align}
\|\pp_a(\Pi_{s=1}^lB_{i_s})|z\|_{\mathcal L_2}\leq 2\delta(i_{u+1}-a),
\end{align}
from which the proposition follows. \hfill $\Box$\\
\ \\\ \\\
\noindent{\bf Proof of Proposition \ref{newprop6} in the main article.}
Recall the definition of $\delta(k)$ in Proposition \ref{prop5}. Notice that by the orthogonality of $\pp_k's$ defined in Proposition \ref{prop5} and Fubini's theorem, we have that
\begin{align}\label{newproof111}
&|Cov(\Pi_{u=1}^l\tilde B_{i_u}, \Pi_{v=l+1}^{2l}\tilde B_{i_v}|z)|\notag\\&=|Cov(\Pi_{u=1}^l\tilde B_{i_u}-\E(\Pi_{u=1}^l\tilde B_{i_u}|z), \Pi_{v=l+1}^{2l}\tilde B_{i_v}-E( \Pi_{v=l+1}^{2l}\tilde B_{i_v}|z)|z)|\notag\\
&=\left|\E\left(\sum_{a=-\infty}^N\pp_a(\Pi_{u=1}^l\tilde B_{i_u})
\sum_{b=-\infty}^N\pp_b(\Pi_{v=l+1}^{2l}\tilde B_{i_v})\bigg|z\right)\right|\notag
\\&=\left|\E\left(\sum_{a=-\infty}^N\pp_a(\Pi_{u=1}^l\tilde B_{i_u})
\pp_a(\Pi_{v=l+1}^{2l}\tilde B_{i_v})\bigg|z\right)\right|
\notag\\&\leq \sum_{a=-\infty}^N\|\pp_a(\Pi_{u=1}^l\tilde B_{i_u})|z\|_{\mathcal L_2}
\|\pp_a(\Pi_{v=l+1}^{2l}\tilde B_{i_v})|z\|_{\mathcal L_2}.
\end{align}
By Proposition \ref{prop5}, we have
\begin{align}
\eqref{newproof111}\leq 4\sum_{a=-\infty}^N\left(\sum_{s=1}^l\delta(i_s-a)\mathbf 1(i_{s-1}<a)\sum_{u=l+1}^{2l}\delta(i_u-a)\mathbf 1(i_{u-1}<a)\right)
\end{align}
which shows (a) by straightforward calculations.
For (b), we shall show that \begin{align}\label{new.113}
|Cov(\Pi_{u=1}^l\tilde B_{i_u}, \Pi_{v=l+1}^{2l}\tilde B_{i_v}|z)|\leq 12l(1-2\alpha')^{-1}\chi^{{\iota(i_s)}/2}(1-\chi)^{-1},
\end{align}
for $s=1,...,2l$. Without loss of generality, we only show the case for $s=1$. The other cases follow similar arguments. Direct calculations show that
\begin{align}
&|Cov(\Pi_{u=1}^l\tilde B_{i_u}, \Pi_{v=l+1}^{2l}\tilde B_{i_v}|z)|\notag\\
&=|\E(\Pi_{u=1}^{2l}\tilde B_{i_u}|z)-\E\left(\Pi_{u=1}^l\tilde B_{i_u}|z\right)\E\left( \Pi_{v=l+1}^{2l}\tilde B_{i_v}|z\right)|\leq I+II,\end{align}
where 
\begin{align}
&I:= |\E(\Pi_{u=1}^{2l}\tilde B_{i_u}|z)|=|Cov(\tilde B_{i_1},\Pi_{u=2}^{2l}\tilde B_{i_u}|z)|
\\
&II:=|\E\left(\Pi_{u=1}^l\tilde B_{i_u}|z\right)\E\left( \Pi_{v=l+1}^{2l}\tilde B_{i_v}|z\right)|\leq |Cov(\tilde B_{i_1},\Pi_{u=2}^{l}\tilde B_{i_u}|z)|.\label{new.114}
\end{align}
The second equality of $I$ is due to $\E \tilde B_{i_1}=0$, and the second equality of $II$ is due to $\E \tilde B_{i_1}=0$,$\E \tilde B_{i_{l+1}}=0$.
By similar arguments to equation \ref{newproof111}. we have that for a sufficiently large constant $C>0$,
\begin{align}
&I\leq 4\sum_{a=-\infty}^N\left(\delta(i_1-a)\sum_{u=2}^{2l}\delta(i_u-a)\right),\label{new.116}\\
&II\leq 4\sum_{a=-\infty}^N\left(\delta(i_1-a)\sum_{u=2}^l\delta(i_u-a)\right)
\end{align}
Then \eqref{new.113} follows from \eqref{new.114} and straightforward calculations, which finishes the proof.
\hfill $\Box$\ \\\ \\

\noindent{\bf Proof of Proposition \ref{lemmaconsistent} in the main article}\\\ \\
Let $C$ be a sufficiently large positive generic constant which varies from line to line. Recall the re-parameterization $\theta$ in Section \eqref{MarginalSBM} of the main article. Let $\theta_{i,j}$ be the $(i,j)$-entry of $\theta$.
Define normalized Laplacian $\tilde L=\bar D^{-1/2}A\bar D^{-1/2}$. We shall show (i): there exists a sufficiently large  positive constant $\eta_0$ such that
\begin{align}
\p\left(\|\tilde L\tilde L-\bar L\bar L\|_F\geq\frac{\log n}{\tau^2n^{1/2}}G^{1/2}(\chi,N,3)(1-\chi)^{-1/2}\right)\leq \frac{\eta_0}{n\log ^4 n},
\end{align}
and (ii): there exists a set $E_n$ with its complement $E_n^c$ and  sufficiently large positive constants $\eta_1$ and $M'$, such that when $n\geq M'$
\begin{align}
\p\left(\|\tilde L\tilde L- LL\|_F\geq\frac{\log n}{\tau^2n^{1/2}}G^{1/2}(\chi,N,3)(1-\chi)^{-1/2},E_n \right)=0,\label{variate123}\\
\p(E^c_n)\leq \eta_1n^{-2}.\label{variate124}
\end{align}
where $C_0$ is a sufficiently small positive constant. The LHS of equation \eqref{variate123} means the probability of intercept of events $E_n$ and $\{\|\tilde L\tilde L- LL\|_F\geq\frac{M\log n}{\tau^2n^{1/2}}G^{1/2}(\chi,N,3)(1-\chi)^{-1/2}\}$.
For (i), direct calculation shows that
\begin{align}
|\tilde L\tilde L-\bar L\bar L|_{ij}=\frac{1}{n^2\sqrt{c_ic_j}}\sum_{k=1}^n\left(\frac{A_{i,k}A_{k,j}-\theta_{i,k}\theta_{k,j}}{c_k}\right).
\end{align}
By definition of $\tau$,
we only need to show that
\begin{align}\label{new.129}
\p\left(\sum_{i,j}\left(\sum_{k=1}^n\left(\frac{A_{i,k}A_{k,j}-\theta_{i,k}\theta_{k,j}}{c_k}\right)\right)^2\geq \frac{n^3\log^2n}{\tau^2}G(\chi,N,3)(1-\chi)^{-1}\right)\notag\\\leq \frac{\eta_0}{n\log ^4 n}.
\end{align}
Write $\tilde A_{i,j}=A_{i,j}-\theta_{i,j}$ for the centered adjacency matrix. Then we have
\begin{align}
A_{i,k}A_{k,j}-\theta_{i,k}\theta_{k,j}=\tilde A_{i,k}\tilde A_{k,j}+\theta_{k,j}\tilde A_{i,k}+\theta_{i,k}\tilde A_{k,j}.
\end{align}
As a result, we get
\begin{align}
&\sum_{i,j}\left(\sum_{k=1}^n\left(\frac{A_{i,k}A_{k,j}-\theta_{i,k}\theta_{k,j}}{c_k}\right)\right)^2\notag\\
&=\sum_{i,j,k,l}\left(\left(\frac{\tilde A_{i,k}\tilde A_{k,j}+\theta_{k,j}\tilde A_{i,k}+\theta_{i,k}\tilde A_{k,j}}{c_k}\right)\times\left( \frac{\tilde A_{i,l}\tilde A_{l,j}+\theta_{l,j}\tilde A_{i,l}+\theta_{i,l}\tilde A_{l,j}}{c_l}\right)\right)\notag\\
&=I_1+I_2+I_3+I_4+I_5+I_6+I_7+I_8+I_9,
\end{align}
where \begin{align}
I_1=\sum_{i,j,k,l}\frac{\tilde A_{i,k}\tilde A_{k,j}\tilde A_{i,l}\tilde A_{l,j}}{c_kc_l},
I_2=\sum_{i,j,k,l}\frac{\tilde A_{i,k}\tilde A_{k,j}\tilde A_{i,l}\theta_{l,j}}{c_kc_l},
I_3=\sum_{i,j,k,l}\frac{\tilde A_{i,k}\tilde A_{k,j}\theta_{i,l}\tilde A_{l,j}}{c_kc_l},\\
I_4=\sum_{i,j,k,l}\frac{\tilde A_{i,k}\theta_{k,j}\tilde A_{i,l}\tilde A_{l,j}}{c_kc_l},
I_5=\sum_{i,j,k,l}\frac{\tilde A_{i,k}\theta_{k,j}\tilde A_{i,l}\theta_{l,j}}{c_kc_l},
I_6=\sum_{i,j,k,l}\frac{\tilde A_{i,k}\theta_{k,j}\theta_{i,l}\tilde A_{l,j}}{c_kc_l},\\
I_7=\sum_{i,j,k,l}\frac{\theta_{i,k}\tilde A_{k,j}\tilde A_{i,l}\tilde A_{l,j}}{c_kc_l},
I_8=\sum_{i,j,k,l}\frac{\theta_{i,k}\tilde A_{k,j}\tilde A_{i,l}\theta_{l,j}}{c_kc_l},
I_9=\sum_{i,j,k,l}\frac{\theta_{i,k}\tilde A_{k,j}\theta_{i,l}\tilde A_{l,j}}{c_kc_l}.
\end{align}
We will show that, there exist constants $M_s$, $1\leq s\leq 9$ such that
\begin{align}\label{new.135Is}
\p\left(I_s\geq \frac{n^3\log^2n}{\tau^2}G(\chi,N,3)(1-\chi)^{-1}\right)\leq \frac{\eta_{0,s}}{n\log ^4 n}.
\end{align}
By taking $\eta_0=\sum_{s=1}^9\eta_{0,s}$, we shall show \eqref{new.129}. For the sake of brevity, we only show the case that $s=1$, and the situation that $2\leq s\leq 9$ follows from a similar argument. Further define $\bar{A}_{(i,k),(k,j)}=\tilde A_{i,k}\tilde A_{k,j}-\E(\tilde A_{i,k}\tilde A_{k,j}).$
Notice that
$$\E\left(\bar {A}_{(i,k)(k,j)}\tilde A_{i,l}\tilde A_{l,j}\right)=Cov(\bar {A}_{(i,k)(k,j)}\tilde A_{i,l},\tilde A_{l,j})=Cov(\bar {A}_{(i,k)(k,j)},\tilde A_{i,l}\tilde A_{l,j}).$$
By the similar arguments to the proof of Proposition \ref{newprop6}, we have that
\begin{align}
\left|\E\left(\frac{\bar A_{(i,k)(k,j)}\tilde A_{i,l}\tilde A_{l,j}}{c_kc_l}\right)\right|\leq \frac{40(1-2\alpha')^{-1}}{\tau^2}\chi^{\Xi(i,k,j,l, \omega)/2}(1-\chi)^{-1},
\end{align}
where we have defined $$\Xi(i,k,j,l,\omega)=\max\{\max_{3\leq s\leq 4}\min_{1\leq j\leq 4,j\neq s}(|i_s-i_j|),\min(|i_1-i_3|,|i_1-i_4|,|i_2-i_3|,|i_2-i_4|)\},$$ $\chi$ and $\alpha'$ are defined in defined in Proposition \ref{dependenceprop} in the main article, and integers $i_1=\omega(\{i,k\}), i_2=\omega(\{k,j\}),i_3=\omega(\{i,l\}),i_4=\omega(\{l,j\})$. We then argue that for any $\omega(\{\cdot,\cdot\})$ and  given $i,k,j$ and $r$, the possible number of $l$ such that
$\Xi(i,k,j,l,\omega)=r$ is at most $8r$.
We also use the following argument and its analog. If $i_5=\omega(\{u,v\})$ has $k$ possible different values, then the number of  different  values of $u$ and $v$ should be both smaller than
$k+1$. This can be easily seen by contradiction. 
When $i,j,k$ are fixed, $i_1,i_2$ are fixed. For any real number $a$ and integer $r>0$, denote by $B_a(r)$ the interval $[a-r,a+r]$. Then one of the $i_3$ and $i_4$ must fall into $B_{i_1}(r)\cup B_{i_2}(r)$. Otherwise $\min(|i_1-i_3|,|i_1-i_4|,|i_2-i_3|,|i_2-i_4|)>r$, which implies
$\Xi(i,k,j,l,\omega)>r$. Since $\omega(\{\cdot,\cdot\})$ is a $1-1$ map, when $i, j, k$ are fixed, the number of possible choices of $l$ is less than or equal to both the numbers of possible choices of $i_3$ and $i_4$.  The total number of integers in the interval $B_{i_1}(r)\cup B_{i_2}(r)$ are is most $(2r+1)\times 2=4r+2$, which implies that the possible number of $l$ such that $\Xi(i,k,j,l,\omega)=r$ is at most $8r+5$. Thus,
\begin{align}\label{137-centered}
\sum_{i,j,k,l}\left|\E\left(\frac{\bar A_{(i,k)(k,j)}\tilde A_{i,l}\tilde A_{l,j}}{c_kc_l}\right)\right|\leq \frac{40(1-2\alpha')^{-1}}{\tau^2}(1-\chi)^{-1}\sum_{i,j,k}\sum_{r=0}^{N}(8r+5)\chi^{r/2}\notag
\\ \leq \frac{Cn^3}{\tau^2}G(\chi,N,1)(1-\chi)^{-1}.
\end{align}
Similarly, we get
\begin{align}
\left|\E\left(\frac{\E(\tilde A_{i,k}\tilde A_{k,j})\tilde A_{i,l}\tilde A_{l,j}}{c_kc_l}\right)\right|\leq \frac{16(1-2\alpha')^{-1}}{\tau^2}\chi^{| \omega(\{i,l\})- \omega(\{j,l\})|/2}(1-\chi)^{-1},
\end{align}
Then similar but easier arguments of \eqref{137-centered} (by fixing $i$, $l$, $k$) lead to
\begin{align}
\sum_{i,j,k,l}\left|\E\left(\frac{\E(\tilde A_{i,k}\tilde A_{k,j})\tilde A_{i,l}\tilde A_{l,j}}{c_kc_l}\right)\right|\leq \frac{Cn^3}{\tau^2}G(\chi,N,1)(1-\chi)^{-1},
\end{align}
which together with \eqref{137-centered} shows that
\begin{align}\label{new.140-exp}
\E I_1\leq \frac{Cn^3}{\tau^2}G(\chi,N,1)(1-\chi)^{-1}.
\end{align}
On the other hand, Proposition \ref{newprop6} leads to that
\begin{align}
Cov\left(\sum_{i,j,k,l}\frac{\tilde A_{i,k}\tilde A_{k,j}\tilde A_{i,l}\tilde A_{l,j}}{c_kc_l},\sum_{u,w,v,p}\frac{\tilde A_{u,v}\tilde A_{w,v}\tilde A_{u,p}\tilde A_{w,p}}{c_vc_p}\right)\notag
\notag
\\ \leq \frac{112(1-2\alpha')^{-1}\sum_{i,j,k,l}\sum_{u,w,v,p}}{\tau^4}\chi^{\frac{1}{2}\max(\Lambda(i_s,1\leq s\leq 8),\max_{1\leq s\leq 8}{\iota(i_s)})}(1-\chi)^{-1},
\end{align}
where $i_1=\omega(\{i,k\})$, $i_2=\omega(\{i,l\})$, $i_3= \omega(\{j,k\})$, $i_4=\omega(\{j,l\})$,
$i_5=\omega(\{u,v\})$, $i_6=\omega(\{w,v\})$, $i_7=\omega(\{u,p\})$, $i_8= \omega(\{w,p\})$.
For given $i,j,k,l$ and $r$, we shall show that the total number of possible pairs of integers $u,w,v,p$ such that
$\max(\Lambda(i_s,1\leq s\leq 8),\max_{1\leq s\leq 8}{\iota(i_s)})=r$ is at most $(8r +5)^2(12r+7)n$. Recall the definition of $\Lambda(c\dot)$ and $\iota(\cdot)$. On one hand, $\Lambda(i_s,1\leq s\leq 2l)=\min_{1\leq s\leq l}(\min_{l+1\leq u\leq 2l}|i_s-i_u|)$ is the smallest distance between 2 indices, of which one is from $i_1,..,i_l$ and the other is from $i_{l+1},..,i_{2l}$. On the other hand, $\iota(i_s)=\min_{1\leq j\leq 2l,j\neq s}(|i_s-i_j|)$ measures the smallest distance between index $i_s$ and other indices.  We now calculate the total number of possible pairs of integers $u,w,v,p$ such that $\max(\Lambda(i_s,1\leq s\leq 8),\max_{1\leq s\leq 8}{\iota(i_s)})=r$.

First note that at least one of $i_5,i_6,i_7, i_8$ should belong to the interval $B_{i_1}(r)\cup B_{i_2}(r)\cup B_{i_3}(r)\cup B_{i_4}(r)$. Otherwise all distances between two indices, such that one indices is from $i_1,i_2,i_3, i_4$ and the other is from the group of  $i_5,i_6,i_7, i_8$ are large than $r$, so $\Lambda(i_s,1\leq s\leq 8)>r$. Without loss of generality, consider the case that $i_5=\omega(\{u,v\})\in B_{i_1}(r)\cup B_{i_2}(r)\cup B_{i_3}(r)\cup B_{i_4}(r)$, then the possible numbers of $i_5$ is at most $4\times(2r+1)=8r+4$, which implies the possible numbers of different values of $u,v$ are both at most $8r+5$. For each pair of $(u,v)$, define $I=B_{i_1}(r)\cup B_{i_2}(r)\cup B_{i_3}(r)\cup B_{i_4}(r)\cup B_{r_5}(r)$.  Then if (a) $i_6\in I$, the total possible number of different values of $w$ is no larger than one plus the possible numbers of $i_6$, which is  $1+(2r+1)\times 5=10 r+6$. Meanwhile, the total possible number of $p$ is $n$. So for (a) the total possible number of pairs $(u,w,v,p)$ is $(8r+5)^2(10r+6)n$.

If (b) $i_6\not\in I$, then $w$ has  at most $n$ possibilities. However at least one of $i_7=\omega(\{u,p\})$ and $i_8=\omega(\{w,p\})$ should fall into $B_{r_6}(r)$, otherwise $\iota(i_6)>r$.  Since $i_7=\omega(\{u,p\})$ and $i_8=\omega(\{w,p\})$, the possible number of different values of $p$ is bounded by one plus the minimal possible numbers of $i_7$ and $i_8$  Thus, $p$ has at most $2r+1$ choices, and for case (b) the  total possible number of pairs $(u,w,v,p)$ is $(8r +5)^2(2r+1)n$. Combining (a), (b), we find that the total possible number of pairs $(u,w,v,p)$ is $(8r +5)^2(12r+7)n$. The explanation of (b) is in the scanned figure.

In conclusion, we have that
\begin{align}
&\frac{112(1-2\alpha')^{-1}\sum_{i,j,k,l}\sum_{u,w,v,p}}{\tau^4}\chi^{\frac{1}{2}\max(\Lambda(i_s,1\leq s\leq 8),\max_{1\leq s\leq 8}{\iota(i_s)})}(1-\chi)^{-1}\notag\\
\leq  &\frac{112(1-2\alpha')^{-1}n^4}{\tau^4} \sum_{r=0}^N(8r +5)^2(12r+7)n\chi^{r/2}(1-\chi)^{-1}\notag\\\leq &\frac{Cn^5}{\tau^4}G(\chi,N,3)(1-\chi)^{-1}.
\end{align}
Equation \eqref{new.140-exp} with Markov's inequality show \eqref{new.135Is} (with $s=1$). As a result, (i) follows from similar arguments to $I_v, 2\leq v\leq 9$. For (ii), define \begin{align} \label{En-def}E_{n}=\cap_{i}\{D_{ii}\in\bar D_{ii}(l(n),u(n))\}\end{align} with
$l(n)=1-b(n)$, $u(n)=1+b(n)$, where we have chosen \begin{align}b(n)=\frac{\log n}{16 \sqrt n}G^{1/2}(\chi,N,3)(1-\chi)^{-1/2}.\end{align} We take $\mathfrak D$ such that
\begin{align}
\mathfrak D:=|D_{ss}^{-1}(D_{ii}D_{jj})^{-1/2}-\bar D_{ss}^{-1}(\bar D_{ii}\bar D_{jj})^{-1/2}|\leq \frac{1}{\tau^2n^2}\left|\frac{\bar D_{ss}^{}(\bar D_{ii}\bar D_{jj})^{1/2}}{D_{ss}^{}(D_{ii}D_{jj})^{1/2}}-1\right|.
\end{align}
On the other hand, on set $E$ such that 
for all $i$, $D_{ii}\in \bar D_{ii}(l(n),u(n))$, 
\begin{align}
\frac{\bar D_{ss}^{}(\bar D_{ii}\bar D_{jj})^{1/2}}{D_{ss}^{}(D_{ii}D_{jj})^{1/2}}-1\in ((1+b(n))^{-2}-1, (1-b(n))^{-2}-1).
\end{align}
Notice that $b(n)\rightarrow 0$. Let $M'$ be the integer that $b(n)<0.5$ if $n\geq M'$,
then mean value theorem leads to  \begin{align}
\mathfrak D\leq \frac{16b(n)}{\tau^2n^2}.
\end{align}
By the boundedness of a Bernoulli random variable, we have that when $n\geq M'$,
\begin{align}
|LL-\tilde L\tilde L|_{i,j}\leq \sum_{s=1}^n|A_{i,s}A_{s,j}||D_{ss}^{-1}(D_{ii}D_{jj})^{-1/2}-\bar D_{ss}^{-1}(\bar D_{ii}\bar D_{jj})^{-1/2}|\\
\leq \frac{16b(n)}{n^2\tau^2}\sum_{s=1}^n|A_{i,s}A_{s,j}|\leq \frac{16(n)}{n\tau^2}.
\end{align}
As a result, we have that when $n\geq M'$
\begin{align}\p\left(\|\tilde L\tilde L- LL\|_F\geq\frac{\log n}{\tau^2n^{1/2}}G^{1/2}(\chi,N,3)(1-\chi)^{-1/2},E_n \right)=0,
\end{align}
which further implies \eqref{variate123}.
It remains to show \eqref{variate124}. The definition of set $E_n$ in equation \eqref{En-def} yields that 
\begin{align}\label{Yu4}
\p(E_{n}^c)\leq \sum_{i}\p\left(\frac{|D_{i,i}-\bar D_{i,i}|}{\sqrt n}\geq \frac{\bar D_{i,i}b_n}{\sqrt n}\right)
\leq \sum_{i} \p\left(\frac{|D_{i,i}-\bar D_{i,i}|}{\sqrt n}\geq \tau \sqrt{n} b_n\right).
\end{align}
Moreover it follows from the proof of Lemma \ref{Newprop3},
that there exists a sufficiently small constant $C_0$ which is independent of $i$ such that
\begin{align}
\E\left(\exp\left(\left(\frac{C_0(1-\chi)}{\sqrt n}(D_{i,i}-\bar D_{i,i})\right)^2\right)\right)\leq 2.
\end{align}
The above expression together with \eqref{Yu4}, Markov's inequality and the fact that $\tau_n^2\log n>M$ for sufficiently large $M$ implies \eqref{variate124}. Thus (ii) follows, which together with (i) completes the proof.
\hfill $\Box$.

\section{Proof of Results in Section \ref{simpleedge}}\label{Sec:4}
We need the following proposition to show Corollary \ref{coroldependent} in the main article.
\begin{proposition}\label{propdependent}
	Define $W_k(B_j)$, $U_k(B_j)$ for $j\geq 0$ iteratively as follows:
	\begin{align}\label{new.A4}
	&W_1(B_j)=p_1^{B_j}(1-p_1)^{1-B_j}, \ \ \ U_1(B_j)=p_0^{B_j}(1-p_0)^{1-B_j},\notag\\
	&W_{k}(B_j)=p_1W_{k-1}(B_j)+(1-p_1)U_{k-1}(B_j),\ \ k\geq 2\notag\\
	&U_{k}(B_j)=p_0W_{k-1}(B_j)+(1-p_0)U_{k-1}(B_j),\ \ k\geq 2.
	\end{align}
	Then the dependence between two ordered edge variables in $CG(V,\omega,p_0,p_1)$ could be represented by $U_k(B_j)$ and $W_k(B_j)$ as follows. For $j-1\geq k\geq 2$, we have \begin{align}\label{induction1}
	\p(B_j|B_{j-k})=W_{k-1}(B_j)p_1^{B_{j-k}}p_0^{1-B_{j-k}}+U_{k-1}(B_j)(1-p_1)^{B_{j-k}}(1-p_0)^{1-B_{j-k}},
	\end{align}
	where $CG(V,\omega,p_0,p_1)$ is defined in definition \ref{defofCMMG} in the main article.
\end{proposition}
{\it Proof.}
We shall use mathematical induction to prove the proposition. For $k=2$, we have that by Markov property of the sequence $\{B_j\}$ and $j\geq 2$, \begin{align}
\p(B_j|B_{j-2})&=\sum_{B_{j-1}=0}^1\p(B_j|B_{j-1})\p(B_{j-1}|B_{j-2})\notag\\
&=p_1^{B_j}(1-p_1)^{1-B_j}p_1^{B_{j-2}}p_0^{1-B_{j-2}}+p_0^{B_j}(1-p_0)^{1-B_j}(1-p_1)^{B_{j-2}}(1-p_0)^{1-B_{j-2}}.
\end{align}
Suppose that for $k=l$, $l\geq 2$, the equation \eqref{induction1} holds. Then for $k=l+1$,
\begin{align}\label{dependent6}
\p(B_j|B_{j-l-1})&=\sum_{B_{j-l}=0}^1\p(B_j|B_{j-l})\p(B_{j-l}|B_{j-l-1})\notag\\
&=(W_{l-1}(B_j)p_1+U_{l-1}(B_j)(1-p_1))p_1^{B_{j-l-1}}p_0^{1-B_{j-l-1}}\notag\\
&+(W_{l-1}(B_j)p_0+U_{l-1}(B_j)(1-p_0))(1-p_1)^{B_{j-l-1}}(1-p_0)^{1-B_{j-l-1}}\notag\\
&=W_l(B_j)p_1^{B_{j-l-1}}p_0^{1-B_{j-l-1}}+U_l(B_j)(1-p_1)^{B_{j-l-1}}(1-p_0)^{1-B_{j-l-1}},
\end{align}
where the equality is due to \eqref{new.A4}. By mathematical induction the proposition follows.\hfill $\Box$\\\ \\
\noindent {\bf Proof of Corollary \ref{coroldependent} in the main article}\\ \\
Note that \eqref{new.A4} in Proposition \ref{propdependent} implies that
\begin{align}
\begin{pmatrix}
W_k(B_j)\\U_k(B_j)
\end{pmatrix}=
\begin{pmatrix}p_1& 1-p_1\\p_0& 1-p_0
\end{pmatrix}^{}\begin{pmatrix}
W_1(B_{j-1})\\U_1(B_{j-1})
\end{pmatrix},
\end{align} which by further  iteration results in
\begin{align}
\begin{pmatrix}
W_k(B_j)\\U_k(B_j)
\end{pmatrix}=
\begin{pmatrix}p_1& 1-p_1\\p_0& 1-p_0
\end{pmatrix}^{k-1}\begin{pmatrix}
W_1(B_j)\\U_1(B_j)
\end{pmatrix}.
\end{align}
Define the matrix $\mathbf P= {\left(\begin{smallmatrix}p_1 & 1-p_1\\p_0 &1-p_0
	\end{smallmatrix}\right)}$. 
Direct calculations show that the eigenvalues of $\mathbf P$ are $\lambda_1=1$ and $\lambda_2=p_1-p_0$. Let $\mathbb{I}$ be the $2\times 2$ identity matrix. The by Cayley-Hamilton Theorem, we have the following formula
Let notation $\mathbf P^k$ represent the $n_{th}$ power of matrix $\mathbf P$.
\begin{align}
\mathbf P^k=\frac{\lambda_1^{k}-\lambda_2^{k}}{\lambda_1-\lambda_2}\mathbf P-\lambda_1\lambda_2\frac{\lambda_1^{k-1}-\lambda_2^{k-1}}{\lambda_1-\lambda_2}\mathbb{I},
\end{align}
which by simple calculations leads to
\begin{align}
\mathbf P^k=\frac{1}{1-p_1+p_0}\begin{pmatrix}p_0+(1-p_1)(p_1-p_0)^k& (1-p_1)(1-(p_1-p_0)^k)\\
p_0(1-(p_1-p_0)^k) &(1-p_1)+p_0(p_1-p_0)^k\end{pmatrix}.
\end{align}
Then the corollary is a consequence of the fact that $W_1(1)=p_1, U_1(1)=p_0, W_1(0)=1-p_1, U_1(0)=1-p_0$, Proposition \ref{propdependent} and straightforward calculations.\hfill $\Box$\\

The results show that the dependence between  edge variables, or equivalently  $\Delta_n(k)$ defined in \eqref{dependencemeasure} decays at the geometric rate $(p_1-p_0)^k$. 
We also observe that by Proposition \ref{propdependent}, for all $k\geq 2$, the following bounds hold for the conditional probabilities that defined in Corollary \ref{coroldependent}: 
\begin{align}\label{WUbound}
&\p_k(0|1)\in [1-\bar p_1, 1-\bar p_0],
&\p_k(0|0)\in [1-\bar p_1, 1- \bar p_0],\\
&\p_k(1|0)\in [\bar p_0,\bar p_1], &\p_k(1|1)\in  [\bar p_0,\bar p_1],\notag
\end{align}
where $\bar p_1=p_0\vee p_1$ and $\bar p_0=p_0\wedge p_1$, where $\vee$ represents `$\max$' and $\wedge$ represents `$\min$'. \ \\
We then introduce an auxiliary proposition about the rate of Poisson approximation to Binomial distribution, which has been used to compare  MECLTG $CG(V,\omega_1,p_0,p_1)$  with Erd\"{o}s-R\'{e}nyi graph $G(V,p)$  below Theorem  \ref{degree} of the main article.
\begin{proposition}\label{possionapprox}
	Suppose 
	$k^2=o(n)$, 
	Let $Z_n$ follow a Bernoulli$(n,\lambda_0/n)$  distribution with $\lambda_0:=\lambda_0(n)$ satisfying $\lambda_0k=o(n)$.
	Let $Y$ follow a  Poisson($\lambda_0$) distribution.
	Then it follows that \begin{align}
	\p(Z_n=k)=\p(Y=k)\left(1+O\left(\frac{k^2+\lambda^2_0}{n}\right)\right).
	\end{align}
\end{proposition}
{\it Proof}. 
The proof proceeds by direct calculation. Note that
\begin{align}
\p(Z_n=k)&={n\choose k}\left(\frac{\lambda_0}{n}\right)^k\left(1-\frac{\lambda_0}{n}\right)^{n-k}\notag\\
&=\frac{\lambda_0^kn!}{k!(n-k)!}(n-\lambda_0)^{-k}\exp(-\lambda_0)\left(1+O\left(\frac{\lambda^2_0}{n}\right)\right)\notag\\
&=\frac{\lambda_0^k}{k!}\exp(-\lambda_0)\left(1+O\left(\frac{\lambda_0^2}{n}\right)\right)U(\lambda_0,n,k),
\end{align}
where $U(\lambda_0,n,k)=\frac{\Pi_{s=1}^k(n-s+1)}{(n-\lambda_0)^k}$. Direct calculations shows that
$$U(\lambda_0,n,k)\in [U_1(\lambda_0,n,k),U_2(\lambda_0,n,k)],$$ where
\begin{align}
U_1(\lambda_0,n,k)=\left(1+\frac{\lambda_0-k+1}{n-\lambda_0}\right)^k=1+O\left\{\frac{k(\lambda_0-k)}{n}\right\},\notag\\
U_2(\lambda_0,n,k)=\left(1+\frac{\lambda_0}{n-\lambda_0}\right)^k=1+O\left(\frac{\lambda_0k}{n}\right).
\end{align}
which completes the proof. The last step uses  that $k^2/n\rightarrow 0$,\begin{align}
\left(1+\frac{k}{n}\right)^k=\exp\left\{k\ln\left(1+\frac{k}{n}\right)\right\}=\exp\left\{\frac{k^2}{n}+O\left(\frac{k^3}{n^2}\right)\right\}=1+O\left(\frac{k^2}{n}\right).
\end{align}
In fact, this proposition implies that the distribution of a Bernoulli random variable with a small probability of a successful trial is close to that of a Possion random variable,  and quantifies the difference between the two distributions by their parameters.  
\hfill $\Box$\\

The following proposition is required for proving Theorem \ref{degree} in the main article.
\begin{proposition}\label{pro-degree}
	Considering the MECLTG $CG(V,\omega_1,p_0=\frac{\lambda_0}{n}, p_1=1-\frac{\lambda_1}{n^c})$ define in Definition \ref{defofCMMG} of the main article with $\lambda_0>1$.
	Let $n''$ be an integer which is smaller than $\lf \epsilon n \rf$ for some $\epsilon>0$, and $D_{n''}=\sum_{j=1}^{n-n''}I(A_{n'',n''+j}=1)$. Then for any integer $0\leq k\leq n-n''-2$, there exist strictly positive constants $C_1, C_2$ independent of $k$ and $n$, such that
	\begin{align}
	\p(D_{n''}=k)\geq C_1n^{-c}\left(\frac{n'-k-1}{n'}\right)\exp\left(-C_2\left(\frac{k}{n^c}\right)\right).
	\end{align}
\end{proposition}
{\it Proof.} Write $n'=n-n''$. To simplify the notation, we set $\binom{n}{k}=0$ if $k<0$ or $k>n$. 
The $k=0$ case is obvious. For $1\leq k\leq n'-2$ we have
\begin{align}\label{net.42}
\p(D_{n''}= k)&= p\sum_{d=2}^k\binom{n'-k-1}{d-2}\binom{k-1}{d-1}p^{k-d}_1p_0^{d-1}(1-p_1)^{d-1}(1-p_0)^{n'-k-d+1}
\notag\\&+(1-p)\sum_{d=1}^k\binom{n'-k-1}{d-1}\binom{k-1}{d-1}p^{k-d}_1p_0^{d}(1-p_1)^{d-1}(1-p_0)^{n'-k-d}
\notag\\&+p\sum_{d=1}^k\binom{n'-k-1}{d-1}\binom{k-1}{d-1}p^{k-d}_1p_0^{d-1}(1-p_1)^{d}(1-p_0)^{n'-k-d}
\notag\\&+(1-p)\sum_{d=1}^k\binom{n'-k-1}{d}\binom{k-1}{d-1}p^{k-d}_1p_0^d(1-p_1)^{d}(1-p_0)^{n'-k-d-1}
\notag\\&:=I+II+III+IV,
\end{align}
where $p=\frac{p_0}{1-p_1+p_0}$, as discussed in \eqref{stationp} of the main article.
The latter equations define the fewer known I,II, III \& IV.
The first term is correspond to the situation that $(A_{n'',n''+1}, A_{n'',n})=(1,1)$. The second term is correspond to the case of $(0,1)$, The third is correspond to $(1,0)$, and the last term is correspond to the case of $(0,0)$. We only discuss the quantity $I$ in \eqref{net.42} since $II$, $III$, $IV$ follow {\it mutatis mutandis}. To this end, we need to calculate for each $k$, \begin{align}\p(A_{n'',n''+1}=a_{1},...,A_{n'',n}=a_{n'})\label{element}\end{align}
where $\{a_{i}\in\{0,1\}, 1\leq i\leq n': \sum_{i=1}^{n'}a_i=k,a_1=1,a_{n'}=1 \}$. If $a_1=1$, $a_{n'}=1$, and in the series of $a_i, 1\leq i\leq n'$, there are $d,2\leq d\leq n'$ pieces of consecutive $1's$, then we have that, starting with $1$, the series will have $k-d$ of $1'$s following $1$, $d-1$ of $0's$ following $1$, $d-1$ of $1's$ following $0$, and
$n'-(k-d+d-1+d-1)-1=n'-k-d+1$ of $0's$ following $0$, where the last term of $-1$ in the LHS of the above equation is due to $a_1=1$. Define $U_{d,k}=\{\{a_{i}, 1\leq i\leq n'\}\in\{0,1\}^{n'}: \sum_{i=1}^{n'}a_i=k, \text{there are $d$ strings of $1's$},a_1=1,a_{n'}=1\}.$  Here ``a strings of $1$'' means a sub-series only containing $1$. Then we have that
\begin{align*}
\p(A_{n'',n''+1}=a_{1},...,A_{n'',n}=a_{n'},\{a_{i}\}_{i=1}^{n'}\in U_{d,k})=pp_1^{k-d}(1-p_1)^{d-1}p_0^{d-1}(1-p_0)^{n'-k-d+1}.
\end{align*}
Now we study $|U_{d,k}|$ as follows. Since there are $n'-k$ of $0's$,
we first choose $d-2$ out of $n'-k-1$ positions to place strings of $1's$, where $d-2$ is due to $a_1=1$ and $a_{n'}=1$.
Then we assign the length of each string of $1$. This leads to choose $d-1$ out of $k-1$. Thus we obtain
\begin{align*}
&\sum_{\{a_{i}, 1\leq i\leq n'\}\in U_{d,k} }\p(A_{n'',n''+1}=a_{1},...,A_{n'',n}=a_{n'})\notag\\&=\binom{n'-k-1}{d-2}\binom{k-1}{d-1}pp_1^{k-d}(1-p_1)^{d-1}p_0^{d-1}(1-p_0)^{n'-k-d+1},
\end{align*}
Also, we have $\p(D_{n''}=k,A_{n'',n''+1}=1,A_{n'',n}=1)=\sum_{d=2}^{k}\sum_{\{a_{i}, 1\leq i\leq n'\}\in U_{d,k} }\p(A_{n'',n''+1}=a_1,...,A_{n'',n}=a_{n'})$. Similar arguments apply to the  $(A_{n'',n''+1}, A_{n'',n})$=$(0,1), (1,0),$ and the $(0,0)$ cases. Notice that $I,II,III$ and $IV$ are all positive. In particular, $IV$ is greater than its first term, which is
\begin{align}\label{lowerboundIV}
IV&\geq (1-p)\binom{n'-k-1}{1}\binom{k-1}{0}p^{k-1}_1p_0^{}(1-p_1)^{}(1-p_0)^{n'-k-2}\notag\\
&=\frac{(1-p_1)^2}{1-p_1+p_0}p_1^{k-1}\binom{n'-k-1}{1}p_0(1-p_0)^{n'-k-2}.
\end{align}
Since $p_0=\frac{\lambda_0}{n}$, it follows from straightforward calculations that for sufficiently large $n'$ there exists a positive constant $\eta_0<1$ such that  \begin{align}\binom{n'-k-1}{1}p_0(1-p_0)^{n'-k-2}\geq \eta_0\frac{(n'-k-1)\lambda_0}{n'}\exp(-\lambda_0).\label{lowerboundp0}\end{align}
Since $p_1=1-\frac{\lambda_1}{n^c}$,  expression \eqref{lowerboundp0} yields that if $n'$ is sufficiently large 
\begin{align}
\eqref{lowerboundIV}&\geq \eta_1\eta_0\frac{n'-k-1}{n'}\lambda_0\exp(-\lambda_0)\frac{\lambda_1}{n^c}(1-\frac{\lambda_1}{n^c})^{k-1}\notag\\&\geq \eta_2\eta_0\frac{n'-k-1}{n'}\lambda_0\exp(-\lambda_0)\frac{\lambda_1}{n^c}\exp(-\frac{(k-1)\lambda_1}{n^c}),
\end{align}
for positive constants $0<\eta_2<\eta_1<1$,
which completes the proof.
\ \\
\noindent{\bf Proof of Theorem \ref{degree} in the main article}\\

Let $n''$ be an integer which is smaller than $\lf\epsilon n\rf$ for some fixed $\epsilon\in (0,1)$. We consider node $n''$. 
Denote by $d_{n''}$ the degree of the node $n''$. Notice that for node $n''$, the  edge variables are generated in the following order
$$A_{1,n''},A_{2,n''},...,A_{n''-2,n''}, A_{n''-1,n''}, A_{n'',n''+1},..., A_{n'',n}.$$ As a result, the degree has the following decomposition 
\begin{align}
d_{n''}=D_{n''}+W_{n''},\ \text{where}\ \label{new.A32}
D_{n''}=\sum_{j=1}^{n-n''}I(A_{n'',n''+j}=1), \quad W_{n''}=\sum_{j=1}^{n''-1}I(A_{j,n''}=1).
\end{align}
For $k\in \mathbb N$, \begin{align}\label{degree-rep1}
\p(d_{n''}=k)=\sum_{s=0}^k\p(D_{n''}=k-s|W_{n''}=s)\p(W_{n''}=s).
\end{align}
Recall the definition of $\omega_1$ in Section \ref{sec4.1} of the main article, i.e., $\omega_1(\{i,j\})=\varpi(i\wedge j, i\vee j)$ for $i\neq j$, where $\varpi_1(i,j)=n(i-1)-i(i-1)/2+j-i$ for $1\leq i<j\leq n$.
Direct calculations based on the form of $\omega_1$ show that $\omega_1(\{n''-k,n''\})-\omega_1(\{n''-k-1,n''\})=n-n''+k$ for $1\leq k\leq n''-2$. Also we have that $\omega_1(\{n'',n''+1\})-\omega_1(\{n''-1,n''\})=n-n''+1$.
Thus it follows from from \eqref{new.A32} that\begin{align}\label{New.s104-2017}
&\p(D_{n''}=k-s|W_{n''}=s)=\notag\\&\sum_{
	\substack{\sum_{l=1}^{n-n''}a_{n'',n''+l}=k-s,\\\{a_{n'',n''+l}\}_{l=1}^{n-n''}\in\{0,1\}^{n-n''}}}\Pi_{l=2}^{n-n''}
\p(A_{n'',n''+l}=a_{n'',n''+l}|A_{n'',n''+l-1}=a_{n'',n''+l-1})\notag\\
&\times\p(A_{n'',n''+1}=a_{n'',n''+1}|W_{n''}=s).
\end{align}
On the other hand, by Markov property and theory of total probability, we obtain that \begin{align}\label{circle3}
&\p(A_{n'',n''+1}=a_{n'',n''+1},W_{n''}=s)\notag
\\&=\sum_{\substack{\sum_{l=1}^{n''-1}a_{l,n''}=s,\notag\\\notag\{a_{l,n''}\}_{l=1}^{n''-1}
		\in\{0,1\}^{n''-1}}}\p(A_{n'',n''+1}=a_{n'',n''+1}|A_{n''-1,n''}=a_{n''-1,n''})
\\&\times \Pi_{l=1}^{n''-2}\p(A_{n''-l,n''}=a_{n''-l,n''}|A_{n''-l-1,n''}=a_{n''-l-1,n''})\p(A_{1,n''}=a_{1,n''}).
\end{align}
By using Corollary \ref{coroldependent} in the main article, we get
\begin{align}\label{NewA.36-2017}
&\p(A_{n'',n''+1}=a_{n'',n''+1}|A_{n''-1,n''}=a_{n''-1,n''})\notag
\\&=\p(A_{n'',n''+1}=a_{n'',n''+1})(1+O(n^{1-c}(p_1-p_0)^{n-n''+1}))\notag
\\&=\p(A_{n'',n''+1}=a_{n'',n''+1})(1+O(n^{1-c}\exp(-n^{\lambda_1(1-c)})))
\end{align}
if $n$ is sufficiently large. Then by \eqref{NewA.36-2017}, \eqref{circle3} can be simplified to equal to\begin{align}
\eqref{circle3}=\p(A_{n'',n''+1}=a_{n'',n''+1})\p(W_{n''}=s)\{1+O(n^{1-c}\exp(-n^{\lambda_1(1-c)}))\},
\end{align}
which further implies that
\begin{align}\label{degree-rep2}
\p(A_{n'',n''+1}=a_{n'',n''+1}|W_{n''}=s)=\p(A_{n'',n''+1}=a_{n'',n''+1})\{1+O(n^{1-c}\exp(-n^{\lambda_1(1-c)}))\}.
\end{align}
Therefore \eqref{degree-rep1}, \eqref{New.s104-2017} and \eqref{degree-rep2} result in\begin{align}\label{degree-rep3} 
\p(d_{n''}=k)=\sum_{s=0}^k\p(D_{n''}=k-s)\p(W_{n''}=s)\{1+O(n^{1-c}\exp(-n^{\lambda_1(1-c)}))\}.
\end{align}
Similar arguments using the fact $\omega_1(\{n''-k,n''\})-\omega_1(\{n''-k-1,n''\})=n-n''+k$ yield that
\begin{align}\label{Y_n-Wn}
\p(W_{n''}=s)&=\left(\sum_{\substack{\sum_{l=1}^{n''-1}a_{l,n''}=s,\\\{a_{l,n''}\}_{l=1}^{n''-1}\in\{0,1\}^{n''-1}}}
\Pi_{l=1}^{n''-1}\p(A_{l,n''}=a_{l,n''})\right)\{1+O(n^{2-c}\exp(-n^{\lambda_1(1-c)}))\}\notag
\\&=\p(Y_{n''}=s)(1+O(n^{2-c}\exp(-n^{\lambda_1(1-c)}))),
\end{align}
where $Y_{n''}\sim Binomial(n''-1,p)$, for $p=\frac{p_0}{1-p_1+p_0}$ as discussed in equation \eqref{stationp} in the main article. By Proposition \ref{pro-degree}, there exist constants positive $C_1$, $C_2$ such that for $0\leq k\leq n-n''-2$
\begin{align}
\p(D_{n''}=k)\geq C_1n^{-c}\left(\frac{n'-k-1}{n'}\right)\exp\left(-C_2\left(\frac{k}{n^c}\right)\right).
\end{align}
Let $k^*=\lf k/p\rf$.
Then we have for 
$0\leq k\leq (n-n''-2)$,
\begin{align}
&\E\left(\frac{\sum_{n''=1}^nI(d_{n''}=k)}{n}\right)\geq \frac{\sum_{n''=1}^{k^* \wedge n} \p(d_{n''}=k)}{n}\notag
\\&=\frac{1}{n}\sum_{n''=1}^{k^* \wedge n}\sum_{s=0}^k\p(W_{n''}=k-s)\p(D_{n''}=s)\{1+O(n^{1-c}\exp(-n^{\lambda(1-c)}))\}\notag
\\&\geq \notag\frac{1}{n}\sum_{n''=1}^{k^* \wedge n}\sum_{s=0}^k\p(W_{n''}=k-s)C_1n^{-c}\exp\left(-C_2\left(\frac{k}{n^c}\right)\right)\{1+O(n^{1-c}\exp(-n^{\lambda(1-c)}))\}
\\&= \frac{1}{n}\sum_{n''=1}^{k^* \wedge n}\sum_{s=0}^k\p(Y_{n''}=k-s)C_1n^{-c}\exp\left(-C_2\left(\frac{k}{n^c}\right)\right)
\notag\\&\times \{1+O(n^{1-c}\exp(-n^{\lambda(1-c)}))\}\{1+O(n^{2-c}\exp(-n^{\lambda(1-c)}))\}
\notag\\&\geq\frac{1}{n}\sum_{n''=1}^{k^* \wedge n}\sum_{s=0}^k\p(Y_{n''}=k-s)C_1n^{-c}\exp\left(-C_2\left(\frac{k}{n^c}\right)\right)
\end{align}
for some positive constant $\tilde C_1$, where the equality is due to equation \eqref{Y_n-Wn}. 
Furthermore, observe that by the basic property of a Binomial random variable 
\begin{align}
\sum_{n''=1}^{\lf k/p \rf\wedge n}\p(Y_{n''}\leq k)\geq  \zeta((k/p)\wedge n)
\end{align}
for some constant $\zeta$,
which leads to 
\begin{align}
&\frac{1}{n}\sum_{n''=1}^{k^* \wedge n}\sum_{s=0}^k\p(Y_{n''}=k-s)C_1n^{-c}\exp\left(-C_2\left(\frac{k}{n^c}\right)\right)
\notag\\ &=\frac{1}{n}\sum_{n''=1}^{\lf k/p\rf}\p(Y_{n''}\leq k)C_1n^{-c}\exp\left(-C_2\left(\frac{k}{n^c}\right)\right)
\notag\\ &\geq \check C_1\frac{1}{n^{c}}\left(\frac{k}{n^c}\wedge 1\right)\exp\left(-C_2\left(\frac{k}{n^c}\right)\right)
\end{align} 
for some sufficiently small positive constant $\check C_1$. 
Define $M_\gamma=(\sum_{k=1}^n\frac{1}{k^\gamma})^{-1}$ for $\gamma>1$,
$M_{\gamma,\mu}=(\sum_{k=1}^n\frac{1}{k^\gamma}\exp(-\mu k))^{-1}$ for $\gamma>1$, $\mu>0$. Denote $J_{n''}=[1,n-n''-2]\cap \mathbb Z$, and \begin{align}
\tilde A_{n,\gamma}:=\left\{k\in J_{n''} :\check C_1\frac{1}{n^{c}}\left(\frac{k}{n^c}\wedge 1\right)\exp\left(-C_2\left(\frac{k}{n^c}\right)\right)\geq M_\gamma \frac{1}{k^\gamma}\right\},\\
\tilde B_{n,\gamma,\mu}=\left\{k\in J_{n''}:\check C_1\frac{1}{n^{c}}\left(\frac{k}{n^c}\wedge 1\right)\exp\left(-C_2\left(\frac{k}{n^c}\right)\right)\geq M_{\gamma,\mu} k^{-\gamma}\exp(-\mu k)\right\}.
\end{align}
Direct calculations show that there exist $a_0,b_0,c_0,d_0$ such that for sufficiently large $n$,
\begin{align}
\left\{k:\lf a_0n^{\frac{2c}{1+\gamma}}\rf\leq k\leq \lf b_0n^c\log^{}n\rf\right\}\subset \tilde A_{n,\gamma}\subset A_{n,\gamma},\\
\left\{k:\lf  c_0\log n\rf\leq k\leq \lf d_0n \rf\right\}\subset \tilde B_{n,\gamma,\mu}\subset B_{n,\gamma,\mu},
\end{align}
which completes the proof. \hfill$\Box$\\\ \\

\noindent{\bf Proof of Lemma \ref{homoorder} in the main article}\\ \\

{\it Proof.} Without loss of generality, we assume that $n$ is even. The case of $n$ odd can be shown {\it mutatis mutandis}.
Consider node $k$ such that $2\leq k \leq n/2$. Note that for node $k$, the order of the corresponding edge variables are generated in the following order:
\begin{align*}\mathcal A_k:=(A_{k-1,k}, A_{k,k+1}, A_{k-2,k}, A_{k,k+2},...,A_{k-a,k},\\A_{k,k+a},... A_{1,k},A_{k,2k-1},A_{k,2k},A_{k,2k+1},...,A_{k,n}),\end{align*}
where $\mathcal A_k$ is a $(n-1)$--vector, with $\mathcal A_k^{(l)}$  its $l_{th}$ element.
Recall that $\omega_2(\{i,j\})=i+\frac{(2n-(j-i))(j-i-1)}{2}$. Straightforward calculations show that i): if $|a-b|<|c-d|$, $\omega_2(\{a,b\})< \omega_2(\{c,d\})$  and ii): if $|a-b|=|c-d|$, then  $\omega_2(\{a,b\})< \omega_2(\{c,d\})$ when $\min(a,b)<\min(c,d)$. Let $\mathcal B_i\in \mathbb R^{1\times (n-1)}$ with $\mathcal B_i^{(l)}$ as its $l_{th}$ element. Define for each finite $k$, $1\leq k\leq n/2$, the real series $\{a_i(k)\}_{i=1}^{n-2}$ as
$$
\left\{ \begin{array}{ll}
a_{2s+1}(k)=s+1, &  0\leq s\leq k-2\\
a_{2s}(k)=n-2s-1, & 1\leq s\leq k-2\\
a_l(k)=n+k-1-l, & 2k-2\leq l\leq n-2.
\end{array}
\right.
$$
Let $\mathcal B_k^{(1)}=B_1$, and $\mathcal B_{k}^{(s)}=B_{1+\sum_{u=1}^{s-1}a_u(k)}$ for $s\geq 2$, where $\{B_k\}_{1\leq k\leq N}$ is the ordered edge variables in \eqref{Model1} of the main article with $\p(B_i=1)=\frac{p_0}{1-p_1+p_0}$. Then by the stationary assumption and direct calculation, we can show that
\begin{align}
\mathcal A_k \overset{d}{=} \mathcal B_k
\end{align}
where $\overset{d}{=}  $ means the equivalence in distribution.
For node $k$ with $\frac{n}{2}+1\leq k\leq n$, define similarly
$$
\left\{ \begin{array}{ll}
a_{2s+1}(k)=s+1, &  0\leq s\leq n-k-1\\
a_{2s}(k)=n-2s-1, & 1\leq s\leq n-k\\
a_l(k)=2n-k-1-l, & 2n-2k+1\leq l\leq n-2.
\end{array}
\right.
$$
Let $\mathcal B_{k}^{(s)}=B_{1+\sum_{u=1}^{s-1}a_u(k)}$. By our construction, we have
\begin{align}
\mathcal A_k\overset{d}{=} \mathcal B_k.
\end{align}
Recall $c$ in the definition of $p_1=1-\lambda_1n^{-c}$. Let $\iota$ be a constant such that $c-1+\iota<0$. For the purpose of evaluating $\p(d_i=k)$ we define \begin{align}
d_{i,1}=\sum_{j=1}^{\lf n^\iota\rf}\mathcal B_i^{(j)},
d_{i,2}=\sum_{\lf n^\iota \rf+1}^{n-\lf n^\iota \rf-1}\mathcal B_i^{(j)},d_{i,3}=\sum_{n-\lf n^\iota\rf}^n \mathcal B_i^{(j)}.\end{align}
This therefore splits up $\p(d_i=k)$ into three parts.
Then by similar arguments to the proof of \eqref{New.s104-2017} of Theorem \ref{degree}, we have that
\begin{align}\label{new.14-newdegree}
\p(d_i=k)&=\sum_{s,u}\p(d_{i,3}=u|d_{i,2}=k-s-u)\p(d_{i,2}=k-s-u|d_{i,1}=s)\p(d_{i,1}=s)
\notag\\&=\sum_{s,u}\p(d_{i,3}=u)\p(\bar d_{i,2}=k-s-u)\p(d_{i,1}=s)(1+O(n^{2-c}\exp(-n^{\lambda_1(\iota-c)}))),
\end{align}
where $\bar d_{i,2}$ follows a $Binomial(n-2\lf n^\iota\rf-1,p)$ distribution with $p=\frac{p_0}{1-p_1+p_0}$.
Without loss of generality, we consider node $1$. For node $i,2\leq i\leq n$, the equations follow  {\it mutatis mutandis}.
By similar arguments to the discussion of \eqref{new.14-newdegree}, we then have
\begin{align}
\p(d_1=k)&=\sum_{u}\p(d_{1,3}=u|d_{1,2}+d_{1,1}=k-u)\p(d_{1,2}+d_{1,1}=k-u)\notag
\\&=\sum_{u}\p(d_{1,3}=u)\p(\tilde D_{1,1}=k-u)(1+O(n^{2-c}\exp(-n^{\iota-c}))),\label{degree.new16}
\end{align}
where $\tilde D_{1,1}$ follows a $Binomial (n-\lf n^\iota \rf, p)$ distribution.
Let $V_1$ follow a 
Poisson($n^\iota p$) distribution. 
Let $M$ be a sufficiently large constant.
By Theorem 4.1 of \cite{CR2013} or Theorem 1 of \cite{AGG1989}, we have that for a sufficiently large positive constant $M$,
\begin{align}
\frac{1}{2}\sum_{i=0}^\infty|\p(V_1=i)-\p(d_{1,3}=i)|=M(b_1+b_2+b_3),
\end{align}
where we choose
\begin{align}\label{A.19}
b_1=\sum_{i=1}^{\lf n^\iota\rf}p^2\leq Mn^{\iota+2c-2},b_2=0,\notag\\
b_3=\sum_{i=1}^{\lf n^\iota \rf }\E|\E(B_i|B_j,j\neq i)-p(B_i=1)|\leq 2p\lf n^\iota \rf\leq Mn^{\iota+c-1}.
\end{align}
where for the bound of $b_3$, we have used the property of conditional expectation and the non-negativity of the Binomial random variable.
On the other hand, let $W_i$ follow a $Poisson((n-\lf n^\iota\rf)p)$ distribution. By \cite{chen1974convergence}, we have
\begin{align}\label{A.20}
\sum_{k=0}^\infty|\p(W_i=k)-\p(\tilde D_{1,1}=k)|\leq Mn^{c-1}.
\end{align}
By the property of sum of independent Poisson random variables, \eqref{A.19} and \eqref{A.20}, (i) follows. It remains to show (ii).
As per usual, define the entropy function \begin{align}
H(a,p)=a\log (\frac{a}{p})+(1-a)\log (\frac{1-a}{1-p}).
\end{align}
{\color{black}For any $Binomial (n,p)$ random variable $S_n$, we estimate its tail behavior as follows.} For $0\leq p<a\leq 1$, we have that for $\beta>0$,
\begin{align}\label{June-16-S-126}
\p(S_n\geq an)\leq \exp(-\beta an)\E(\exp(\beta S_n))=\left(\exp\left(-\beta an\right)\left(1-p+p\exp(\beta)\right)\right)^n.
\end{align}
In equation \eqref{June-16-S-126} set $\beta=\log \left(\frac{a(1-p)}{p(1-a)}\right)>0$, then we have that
\begin{align}
\p(S_n\geq an)\leq \exp(-nH(a,p)).\label{A.27-2017}
\end{align}
In the remainder of the proof, we consider node $1$. The other nodes follow {\em mutatis mutandis}.
Let $V$
follow $Binomial(n-\lf n^\iota\rf ,p)$. By \eqref{degree.new16} and the property of Binomial distribution \eqref{A.27-2017}, we have that for $k\geq \lf n^\iota g(n)\rf $, where $g(n)$ is a series of real numbers which go to infinity arbitrarily slowly,  \begin{align}
\p(d_1=k)&\leq \p(V\geq k-\lf n^\iota \rf)(1+O(n^{2-c}\exp(-n^{\lambda_1(\iota-c)})))\notag
\\&\leq  \exp\left(-(n-\lf n^\iota \rf)H\left(a,p\right)\right)(1+O(n^{2-c}\exp(-n^{\lambda_1(\iota-c)}))),
\end{align}
where $a=\frac{k-\lf n^\iota\rf}{n-\lf n^\iota \rf}$.
Since $k\geq \lf n^\iota g(n)\rf$ and $p\asymp n^{c-1}$, when $n$ large enough, $\log (a/p)\geq (\iota-c)\log n$. Therefore it follows from Taylor expansion that
\begin{align}\label{Bound_H}
H(a,p)\geq (\iota-c)a\log n+(1-a)\log(1-a)+O(n^{c-1}).
\end{align}
Notice that \begin{align}\label{limit_a}
\lim_{a\rightarrow 0}\frac{(1-a) \log (1-a)}{a}=-1.
\end{align}
Expression \eqref{Bound_H}, \eqref{limit_a} and the definition of $a$ imply that
\begin{align}
\{n-\lf n^\iota \rf \}H(a,p)\geq \delta(\iota-c)k\log n \ \text{for some $\delta\in(0.5,1)$ }
\end{align}
when $n$ is sufficiently large,
which completes the proof.
\hfill $\Box$
\\\ \\

\bibliographystyle{apalike}
\begin{small}
	\setlength{\bibsep}{2pt}
	\bibliography{lit}
\end{small}

\end{document}